\documentclass[twocolumn]{paper}
\usepackage{amsmath,amssymb,amsfonts}
\usepackage{algorithmic}
\usepackage{graphicx}
\usepackage{textcomp}
\usepackage{subcaption}
\usepackage{tikz}
\usepackage{appendix}
\usetikzlibrary{fit,positioning,shapes,arrows}

\usepackage[a4paper, margin=1.75cm]{geometry}

\usepackage[hidelinks]{hyperref}

\usepackage[acronym]{glossaries}
\newacronym{CF}{CF}{Complementary Filter}
\newacronym{PCF}{PCF}{Passive Complementary Filter}
\newacronym{EKF}{EKF}{Extened Kalman Filter}
\newacronym{IMU}{IMU}{Inertial Measurement Unit}
\newacronym{GNSS}{GNSS}{Global Navigation Satellite System}
\newacronym{ISS}{ISS}{Input-to-State Stable}
\newacronym{UCO}{UCO}{Uniformly Complete Observability}

\newtheorem{problem}{Problem}
\newtheorem{remark}{Remark}
\newtheorem{assumption}{Assumption}
\newtheorem{lemma}{Lemma}
\newtheorem{theorem}{Theorem}

\begin{document}

\tikzstyle{triangle} = [draw, regular polygon, isosceles triangle]
	
\title{High-Performance Motorbike Lean Angle Estimation}
\author{Nicola Mimmo  and Matteo Zanzi
\thanks{N. Mimmo and M. Zanzi are with the
	Department of Electrical, Electronic and Information Engineering "Guglielmo Marconi", University of Bologna,
	Viale del Risorgimento, 2 - 40136 - Bologna (BO) - ITALY (e-mail: \{nicola.mimmo2,matteo.zanzi\}@unibo.it)
}
}
\maketitle
\begin{abstract}
This work deals with the real-time estimation of the lean angle of high-performance motorbikes. The estimate is obtained through measurements provided by an onboard inertial sensor and a GNSS receiver. A two-stage state observer, implementing a kinematic model developed under the novel assumption of coordinated manoeuvre, processes these measurements. A theoretical analysis demonstrates the observer's stability, while a covariance analysis assesses the estimate's accuracy and error bounds. Finally, experimental results obtained on race-track tests and numerical comparisons, with competitive approaches, in simulated realistic scenarios show the superior performance of the proposed estimator.
\end{abstract}


\section{Introduction}
\label{sec:introduction}

\subsection{Motivation}

Lean angle real-time knowledge is crucial for controlling engine and brake power to optimize the motorbike's performance while keeping the biker's safety \cite{Savino2020Active}.
Indeed, the tire-road grip coefficient is a non-linear function of the contact patch shape, which, in turn, depends on the motorbike's leaning \cite{cossalter2014motorcycle}. These non-linearities become critical during high-performance turns when stiff grip variations degrade performance and stability, eventually leading to skidding and highside \cite{BEDOLLA20161872}.


\subsection{State of the art}

The literature extensively investigated the problem of lean angle estimation. The documented solutions can be divided into three main categories: estimation algorithms based on kinematic models (position and velocity), those relying on dynamic models (forces and torques), and image-based. 

Works of the first category present algorithms fed by angular rates, body accelerations, and, eventually, linear speeds and Earth magnetic field measurements. Moreover, the estimation schemes of these works are designed on kinematic models describing the attitude dynamics (commonly Euler's angle dynamics). These models embed gyroscope data as inputs (\textit{e.g.}, for state propagation in Kalman filtering), while non-linear elaborations of accelerometers and linear speeds constitute the output. In detail, \cite{Boniolo2010Motorcycle} and \cite{Maceira2021Roll} propose Complementary Filters (CFs), designed on error frequency separation arguments, which elaborate gyros (motorbike's angular rates) and odometers (wheel speed). Furthermore, \cite{Boniolo2009innovative} proposes strategies using only gyroscope data, whereas \cite{boniolo2014estimation} 
presents a method exploiting only two accelerometers and one gyroscope. These data are successively processed into a CF designed on frequency separation arguments both in \cite{Boniolo2009innovative} and \cite{boniolo2014estimation}. 
Works \cite{Boniolo2010Attitude} and \cite{Corbetta2010Attitude} evaluate the performance of Extended Kalman Filters (EKFs) and unscented Kalman filters applied to the estimation of motorcycles' attitude. These observers rely on the knowledge of the projection on the motorcycle's longitudinal axis of the inertial velocity. Finally, \cite{Sanjurjo2019Roll} provides a scheme for the roll angle estimation relying on a Kalman filter, IMU data, and wheel speed sensors.

Concerning the second category, algorithms rely on dynamic models embedding inertial and geometric data and tire forces descriptions. Moreover, sensor suites comprehend IMU and speed data (as in the first category), potentiometers sensing the steering angle, and torque meters measuring the biker's effort on the handlebar. Commonly, papers in the cited literature assume the knowledge of steering, roll, and yaw angle derivatives. These algorithms focus on estimating a state vector, usually comprehensive of tire forces. The proposed approaches are: Luenberger observers \cite{Damon2019Steering}, EKFs \cite{Teerhuis2012Motorcycle} and \cite{Gabriel2022Accurate},  high-order sliding mode observers \cite{Nehaoua2014Unknown}, \cite{Dabladji2016Unknown}, and \cite{Dabladji2014Estimation}, unknown-input observers \cite{Dabladji2016Unknown} and \cite{Damon2017Lateral},  $H_\infty$ observers \cite{ferjani2018robust}, and adaptive observers \cite{Fouka2022Motorcycle}.

Finally, we report a couple of works belonging to the third category for completeness. In particular, \cite{Damon2018Image} propose using a camera to estimate the motorbike's lean angle. In detail, machine-learning algorithms trained to recognize roll angles from images elaborate onboard camera streams. In addition, \cite{Schlipsing2012Roll} proposes an intriguing comparison between camera-based methods and state observers fed by IMUs.

\subsection{Contribution}

In the context of algorithms based on kinematic models, this paper presents a lean angle estimation approach utilizing standard IMU and GNSS data, such as body accelerations and angular rates (obtained by accelerometers and gyros) and inertial velocities  (from a GPS receiver). 
In particular, we fuse IMU and GNSS data through a novel concept of \textit{coordinated manoeuvre}, which well approximates actual motorbike-plus-biker dynamics.  

The estimator architecture is a cascaded two stages. The first processing level, called \textit{pre-filter}, embeds the coordinated manoeuvre assumption. 
The pre-filter computes a preliminary estimate of the motorbike attitude as a unitary quaternion. The coordinated manoeuvre represents a novel strategy to compensate for the centre-of-gravity displacements due to the biker's movements. This compensation results in a highly accurate estimation, especially when the lean angle data are fundamental, \textit{e.g.}, during high-speed turns. Downstream, an EKF enhances the lean angle estimation by fusing pre-filter and gyroscope outputs.

Theoretical investigations show that the proposed estimator is (locally asymptotically) stable, uniformly on the motorbike's trajectories. Field tests and realistic simulations confirm the good performance of the estimation algorithm proposed in this paper. Finally, a comparison with already existing methods shows the superior performance of the proposed coordinated manoeuvre assumption.

\subsection{Benefits of the proposed approach}

The lean angle estimator designed in this paper has the following benefits.

The overall estimation scheme can be thought as a CF with all the benefits associated to this class of algorithms. In particular, its reduced order (lower than full-order observers with accelerometers and gyroscopes as input and GNSS as a output) lowers the computational burden thus making CFs appealing in real applications.

The proposed estimation scheme does not rely on magnetometers. This improves the estimation accuracy and alleviates the calibration process, as detailed in Remark \ref{rem:Magn}. 

Moreover, the proposed system architecture is more reliable than full-order observers for two reasons. First, 
the proposed algorithm does not suffer from observability issues related to GNSS data unavailability. Second, the CF architecture guarantees estimation stability, although the motorbike does not perform sufficiently exciting trajectories (like on straights). 


\subsection{Notation}

This paper denotes with $\mathbb{R}$ the set of reals and with $\mathbb{N}$
the natural numbers greater than zero. Calligraphic letters, \textit{e.g.}, $\mathcal{X} \subseteq \mathbb{R}^n$, with $n \in \mathbb{N}$, denote subsets. We represent matrices with capital letters, \textit{e.g.}, $X \in \mathbb{R}^{n\times m}$, with $n,m \in \mathbb{N}$. Let $X_i\in\mathbb{R}^{n_i \times m}$ be matrices, with $i = 1,\dots, n$ and $n,\,m,\,n_i \in \mathbb{N}$, then we define $\mathtt{col}\,:\,\mathbb{R}^{n_1 \times m} \times \cdots \times \mathbb{R}^{n_n\times m} \to \mathbb{R}^{(\sum_{i=1}^n n_i)\times m}$ such that
$
\mathtt{col}(X_1,\cdots,X_n) = \left[\begin{array}{ccc}
	X_1^\top& \cdots& X_n^\top
\end{array}\right]^\top
$. 
Symbol $I_n$ denotes identity matrices of size $n \in \mathbb{N}$. Small capital letters, \textit{e.g.}, $x \in \mathbb{R}^n$, with $n \in \mathbb{N}$, denote real vectors of $n$ components. Let $x \in \mathbb{R}^3$ be a vector, then we describe its components with $x_x$,  $x_y$, and $x_z$ such that $x = \mathtt{col}(x_x, x_y, x_z)$. With $\|\cdot\|$, we denote the 2-norm of vectors such that $\|x\|:= \sqrt{x^\top x}$ for any $x \in \mathbb{R}^n$, with $n \in \mathbb{N}$.  Finally, this paper defines $\mathbb{H}$ as the set of unitary-norm quaternions. 

%

\section{Problem formulation and main result}
\label{sec:Problem}


Let $\mathcal{F}_I$ and $\mathcal{F}_B$ be inertial and body reference frames, with the latter rigidly attached to the motorbike. Let $\omega \in \mathbb{R}^3$ be the vector of motorbike angular speeds expressed in $\mathcal{F}_B$. Let $\phi,\theta,\psi \in \mathbb{R}$ be an Euler angle parametrisation for rotation matrices from $\mathcal{F}_I$ to $\mathcal{F}_B$ and define \begin{equation}
	\label{eq:Theta}
	\Theta := \mathtt{col}(\phi,\theta,\psi).
\end{equation}
Then, define $T \,:\, \mathbb{R}^3 \to \text{SO}(3)$ such that $T(\Theta)$ corresponds to the rotation matrix from $\mathcal{F}_I$ to $\mathcal{F}_B$, whose expression is reported in [\cite{titterton2004strapdown}, Eq.(3.63)].
  
Now define $v \in \mathbb{R}^3$ as the motorbike linear speed expressed in $\mathcal{F}_I$.
 Let $\mathtt{v}:=\|v\| \ge 0$ be the inertial speed magnitude, and $\chi,\gamma \in \mathbb{R}$ be the course and the grade angle, then define $$\xi = \mathtt{col}({\mathtt{v}},{\gamma},\chi)$$ such that
\begin{equation}  
	\label{eq:v_aero}
	v = h_v(\xi) := \mathtt{v}\mathtt{col}(\cos\chi\cos\gamma, \sin\chi\,\cos\gamma, -\sin\gamma).
\end{equation}

Denote with $g\in \mathbb{R}^3$ the gravity acceleration expressed in $\mathcal{F}_I$. Then, we made the following assumption with all these quantities at hand.

\begin{assumption}[Sensor Suite]
	\label{hyp:Sensors}
	Assume $\mathcal{F}_B$ be rigidly attached to a combined \acrshort{IMU} and \acrshort{GNSS} board providing
	\begin{equation}
		\label{eq:SensorSuite} 
		\begin{aligned}
			y_a =&\, T(\Theta)  (\dot{v}-g )+ \nu_a(t) && \text{3-axis accelerometer}\\
			\dot{\bar{b}}_g =&\, 0 && \text{3-axis gyro bias}\\
			y_g =&\, \omega + {\nu}_g(t) && \text{3-axis gyroscope}\\
			y_s =&\, h_s(\xi,\tau)+ \nu_s(t) &&\text{GNSS receiver}
		\end{aligned}
	\end{equation} 
	in which, for all $\#\in\{a,g,s\}$,  $y_{\#}$ denotes the sensor output while $\nu_{\#}(t)$ represents bounded measurement errors. More in detail, we define ${\nu}_g(t)=\bar{b}_g(t) + w(t)$ with $w(t)$ an additive error. Let $\overline{\nu}_{a}$, $\overline{\nu}_{g}$, and $\overline{\nu}_{s} > 0$. Then, we assume $\|\nu_{\#}(t)\|_\infty < \overline{\nu}_{\#}$, for all $\#\in\{a,g,s\}$.	\hfill $\square$
\end{assumption}

In agreement with Assumption \ref{hyp:Sensors}, sensors provide  measurements $y_{\#}$ corrupted by  errors $\nu_{\#}$. Moreover, gyroscopes are also affected by the bias $\bar{b}_g$. Finally, the GNSS sampling time, \textit{i.e.}, $\tau > 0$ embedded into $h_s(\cdot,\tau)$, is significant for the application under investigation. In practice, $h_s(\xi,\tau)$ represents a $\tau$-long fixed-period sampling of  $h_v(\xi)$. A description of $h_s(\xi,\tau)$ is given in Section \ref{sec:Pre}, Eq. \eqref{eq:ys}.

\begin{remark}\label{rem:Magn}
The algorithm proposed in this paper does not use data from magnetometers for two main reasons, the distortion of the Earth's magnetic field in the proximity of the motorbike's metal masses and the experimented strong dependence of the magnetometer response on the engine mapping. On the one hand, even if possible for a single test, magnetometer calibrations are time-consuming and too complicated to be carried out during the race weekend. On the other hand, these calibrations require a look-up table to be embedded in the algorithm, thus resulting in a further state dependency, possibly impacting estimation filter stability.  \hfill $\square$
\end{remark}

\begin{problem} [Roll Angle Estimation]
	\label{prob:LeanAngleEstimation}
	Design an algorithm with inputs $y_a(t)$, $y_g(t)$, and $y_s(t)$, state $\hat{x}(t)$, and output $\hat{\phi}(t)$ such that: a) there exists a non empty set of initial conditions, namely $\mathcal{X}_0$, such that $\hat{x}(t)$ is bounded for any $t \ge 0$ and $\hat{x}(0) \in \mathcal{X}_0$; b) there exists $\bar{\phi} > 0$ such that $\limsup_{t \to \infty} \|\hat{\phi}(t)-\phi(t)\| < \bar{\phi}$. \hfill $\square$ 
\end{problem}

\begin{figure}
	\centering	
		\resizebox{\columnwidth}{!}{$
		\begin{tikzpicture}
			\node[draw](Plant){$\begin{aligned}
					\dot{\hat{x}} =&\,  f(\hat{x},t) + K(t)(q_1-\hat{q})\\
					\dot{S} =&\, -SA(\hat{x},t) + \dots\\
				\end{aligned}$};
			\node[draw, shape border rotate=180, triangle, below of = Plant, yshift=-0.5cm](H){$H$};
			\node[draw,circle,left of = Plant, xshift =-1.85cm](sum){$+$};
			\node[draw, right of = Plant, xshift = 2.25cm](h){$h(\hat{x})$};
			\node[ below of = h, yshift = -0.25cm](hatphi){$\hat{\phi}$};
			\node[above of = Plant, yshift = 0.5cm](y2){$y_g(t)$};
			\draw[->] (Plant.east) -| (2.5,-1.5) |- (H.east) ; 
			\draw[->] (H.west) -| (sum.south) node[pos = 0.75, left]{$- \hat{q}$};
			\draw[->] (Plant) -- (h);
			\draw[->] (sum) --  (Plant) ;
			\draw[->] (y2) --  (Plant) ;
			\draw[->] (h) --  (hatphi) ;				
			\node[draw, above of = sum, yshift = 1.25cm](q1){$\bar{h}(f_{\Theta}(y_a,\hat{\xi}_e,\hat{\mathtt{g}}))$};
			\node[draw, right of = q1, xshift = 2.25cm](GNSS){Pre-Filter};
			\node[right of = GNSS, xshift = 0.75cm](yv){$y_s$};
			\node[left of = q1, xshift = -1.25cm](input){$y_a,\,\hat{\mathtt{g}}$};
			\node[above of = Plant, yshift = 0.5cm](y2){$y_g(t)$};
			\draw[->] (q1) --  (sum) node[pos = 0.5, left]{$q_1$};
			\draw[->] (input) -- (q1);
			\draw[->] (GNSS) -- (q1) node[pos = 0.5, above]{$\hat{\xi}_e$};
			\draw[->] (yv) -- (GNSS);
		\end{tikzpicture}
		$}
	\caption{Observer architecture. The attitude propagation through the quaternion dynamics is corrected thanks to the estimation $q_1$. The feedback matrix is $K(t):= S^{-1}(t)H^\top R^{-1}$.}
	\label{fig:CF}
\end{figure}
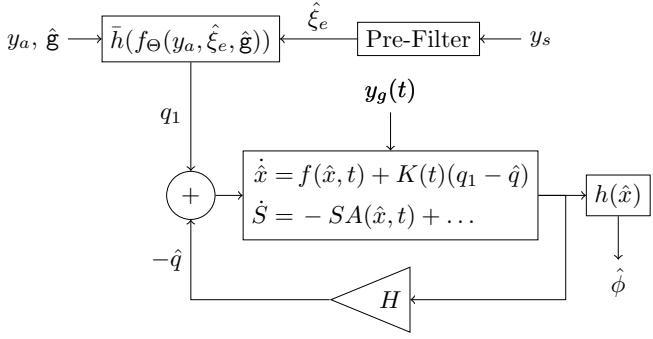

Hereafter, we define some quantities instrumental for introducing the proposed solution, depicted in Figure \ref{fig:CF}.

Let $\bar{h}\,:\,\mathbb{R}^3\to \mathbb{H}$ be such that $q:=\mathtt{col}(q_0,q_x,q_y,q_z) = \bar{h}(\Theta)$ represents the unitary quaternion associated with $\Theta$ (a detailed expression for $\bar{h}(\cdot)$ is reported in [\cite{titterton2004strapdown}, eq. (3.65)]). Moreover, 
the dynamics of $q$ is 
\begin{equation}
	\label{eq:dotTheta}
	\dot{q} = M(q)\, \omega,
\end{equation}
where $M(q)$ is detailed in [\cite{titterton2004strapdown}, eq. (3.61)]. Introduce $\xi_e := \mathtt{col}(\xi,\dot{\xi}),$ define $\mathtt{g} = \|g\| $, and let $\phi_{\text{av}}(\cdot,\cdot,\cdot)\,:\,\mathbb{R}^3\times\mathbb{R}^6\times\mathbb{R} \to \mathbb{R}$ be such that for any $a \in \mathbb{R}^3$, $\xi_e \in \mathbb{R}^6$, and $\mathtt{g} > 0$
\begin{equation}
	\label{eq:phi_av}
	\resizebox{0.89\columnwidth}{!}{$
	\phi_{\text{av}}(a,\xi_e,\mathtt{g}) = \tan^{-1}\Bigg(\dfrac{(\mathtt{g}\cos\gamma -\mathtt{v}\dot{\gamma})a_y-\mathtt{v}\dot{\chi}a_z \cos\gamma}{(\mathtt{g}\cos\gamma -\mathtt{v}\dot{\gamma}) a_z+\mathtt{v}\dot{\chi}a_y\cos\gamma } \Bigg).
	$}
\end{equation}

Define $f_{\Theta}(\cdot,\cdot,\cdot)\,:\,\mathbb{R}^3\times \mathbb{R}^6\times \mathbb{R} \to \mathbb{R}^4$ such that $f_{\Theta}(a,\xi_e,\mathtt{g}) = \mathtt{col}(\phi_{\text{av}}(a, {\xi}_e,\mathtt{g}),\,
{\gamma},\,\chi)$ for each $a \in \mathbb{R}^3$, $\xi_e \in \mathbb{R}^6$, and $\mathtt{g} > 0$. 
Assume $\hat{\xi}_e \in \mathbb{R}^6$ (among whose entries there are $\hat{\gamma}$ and $\hat{\chi}$) be a proxy of $\xi_e$.  Moreover, let $\hat{\mathtt{g}}$ be a proxy of $\mathtt{g}$ and introduce 
\begin{equation}
	\label{eq:hat_theta_e_av}
	\hat{\Theta}_{\text{av}}:= f_{\Theta}(y_a,\hat{\xi}_e,\hat{\mathtt{g}})
\end{equation}
and
\begin{equation} 
	\label{eq:q1}
	q_1 :=\bar{h}(\hat{\Theta}_{\text{av}}).
\end{equation}
  
 Let $x := \mathtt{col}(\bar{b}_g, q)$ 
 and define
\begin{equation}
	\label{eq:f}
f(x,t) = \left[\begin{array}{c}
		0\\
		\begin{array}{r}
		M(q)(y_g(t)-\bar{b}_g)
		\end{array}
	\end{array}
\right]
\end{equation}
and
\[
g(x) = \left[\begin{array}{cccc}
	I & 0& 0 & 0\\
	0 &I & -M(q) & q
\end{array}\right].
\]
Then, to solve Problem \ref{prob:LeanAngleEstimation}, we propose the observer depicted in Figure \ref{fig:CF}. 
In detail, let $\hat{x} := \mathtt{col}(\hat{\bar{b}}_g, \hat{q})$, 
  $\lambda > 0$, $Q = Q^\top \in \mathbb{R}^{12 \times 12}$, $R(\cdot)\,:\,\mathbb{R}\to  \mathbb{R}^{4 \times 4}$ with $R(\cdot) = R^\top(\cdot)$, and $S_0 \in \mathbb{R}^{7 \times 7}$ be such that 
 $Q,\,S_0 \succ 0$, and $R(t)\succ 0$ for all $t \ge 0$. Then, define the following \acrshort{EKF}
\begin{subequations}
	\label{eq:BasicEKF_CF}
	\begin{align}
		\label{eq:EKF_CF}
		\dot{\hat{x}}  =&\, f(\hat{x},t) + 
		S^{-1} {H}^\top 
		R^{-1}(t)(q_1-H \hat{x}) && \hat{x}(0) = \hat{x}_0\\
		\dot{S} =&\,
		 -SA(\hat{x},t)
		-A^\top(\hat{x},t)S \notag\\ 
		&\,-  \lambda  Sg(\hat{x}) Q g^\top(\hat{x}) 
		S   +H^\top R^{-1}(t)H
		&& S(0)=S_0 \label{eq:DRE_CF}\\
			\label{eq:hatphi}
		\hat{\phi} =&\, h(\hat{x})
	\end{align}
\end{subequations}
where 
 $H = \left[\begin{array}{cc}
	0 & I
\end{array}\right]$, $A(x,t) :={\partial f(x,t)}/{\partial x}$, and
\[
h(x):= \left[\begin{array}{ccc}
	1 & 0 & 0 
\end{array}\right] \bar{h}^{-1}(q/\|q\|).
\]

	\begin{remark}
Let $w_1 \in \mathbb{R}^7$ and $w_2 \in \mathbb{R}$, and define $\bar{w} = \mathtt{col}(w_1, \nu_g, w_2)$. Then, we derived $g(x)$ as sum of three contributions, namely $g_1(x)$, $g_2(x)$, and $g_3(x)$, with
\[
\begin{aligned}
	g_1(x) = &\, \left[\begin{array}{cccc}
		0 & 0 & 0 & 0\\
		0 & 0& -M(q) & 0
	\end{array}\right],\, g_2(x) = \left[\begin{array}{cccc}
		I & 0 & 0 & 0\\
		0 & I& 0 & 0
	\end{array}\right]\\
 g_3(x) =&\, \left[\begin{array}{cccc}
		0 & 0 & 0 & 0\\
		0 & 0& 0 & q
	\end{array}\right].
\end{aligned}
\]
The first comes from the linearisation of $\dot{\bar{b}}_g = 0$ and $ \dot{q} = M(q)\omega  = M(q)(y_g - \nu_g)$ with respect to $\bar{w}$, the second keeps $g(\hat{x}) Q g^\top(\hat{x}) \succ 0$ for any $\hat{x} \in \mathbb{R}^7$ (required for the stability of \eqref{eq:BasicEKF_CF}), and the third makes $g(\hat{x}) Q g^\top(\hat{x})$ well conditioned for $\hat{x} \,:\, \|\hat{q}\| \approx 1$ (to improve the performance of \eqref{eq:BasicEKF_CF}). \hfill $\square$
	\end{remark}
	
\begin{remark}
	Observer  \eqref{eq:BasicEKF_CF} provides the estimate  $\hat{q}$, which does not represent a rotation because $\|\hat{q}\|$ is not guaranteed to be unitary. To ensure $\|\hat{q}\| = 1$, one should implement algorithms designed on $\mathbb{H}$, see \cite{Markley2003Attitude,Bonnable2009Invariant} and \cite{Mahony2008Nonlinear}. In \cite{Markley2003Attitude}, the strategy is estimating, through an \acrshort{EKF}, a suitable parametrisation of the attitude (\textit{e.g.}, the Gibbs vector).  The drawback of this approach consists mainly of the non-linearities the \acrshort{EKF} must face. 
	As for \cite{Bonnable2009Invariant}, the observer is composed of a (non-Extended) Kalman Filter designed on a linearisation point. Finally, \cite{Mahony2008Nonlinear} proposes a non-linear \acrshort{CF} whose gain belongs to $\mathbb{R}^2$ (in the case of gyro bias compensation). To the authors' best understanding, the observer gains design is not associated with any physical properties of the sensor suite. 
	
	In this context, the observer proposed in this paper exploits the bi-linear nature of \eqref{eq:dotTheta} to guarantee observability properties, demonstrated in Theorem \ref{Theorem:DREBounds}, which are valid for any observer trajectory (and not only for a linearisation point). Moreover, the design of the observer gain exploits physical features of the selected sensors, thus reducing the number of hand-tuned parameters to one. \hfill $\square$
\end{remark}


Theorem \ref{Theorem:DREBounds}, which summarises the theoretical results of this paper, is valid under the following assumptions.

 \begin{assumption}[Manoeuvres]
		\label{hyp:Non-Ballistic}
		Let $a(t)$ and $\xi_e(t)$ represent the dynamic state of the motorbike at time $t \ge 0$. Then, there exists $\underline{a} > 0$ such that 
		$
			(\mathtt{g}\cos\gamma(t) -\mathtt{v}(t)\dot{\gamma}(t)) a_z(t)
			+\mathtt{v}(t)\dot{\chi}(t)a_y(t)\cos\gamma(t) > \underline{a}
		$
		for all $t\ge 0 $ \hfill $\square$
	\end{assumption}

 \begin{assumption}[Boundedness of $\nu$]
		\label{hyp:nu}
Define $\nu(t)=q_1(t)-q(t)$. Then, there exists $\overline{\nu} > 0$ such that $\|\nu(t)\|_\infty < \overline{\nu}$. \hfill $\square$
		\end{assumption}

\begin{assumption}[Motorbike Pitch]
	\label{hyp:Pitch}
	There exists $\bar{\theta} \in [0,\,\pi/2)$ such that $|\theta(t)|_\infty < \bar{\theta}$.  \hfill $\square$
\end{assumption}

 	Assumption \ref{hyp:Non-Ballistic} ensures $\phi_{\text{av}}(a(t),\xi_e(t),\mathtt{g})$ is well-defined for any $t \ge 0$. 
	In practice, Assumption \ref{hyp:Non-Ballistic} does not represent a limitation. Indeed, for $\underline{a} \ll 1$ and
	with computations similar to those used in the proof of Lemma 1, we can show that the most likely conditions for having Assumption \ref{hyp:Non-Ballistic} not satisfied are ballistic trajectories and turns with extreme roll angles ($|\phi| \approx 90$ deg), which are out of the nominal operating range of on-track race motorbikes.

 Assumption \ref{hyp:nu} is instrumental to assess the local stability of \eqref{eq:BasicEKF_CF}. Section \ref{sec:Rec Err Analysis} deals with the description and analysis of $\nu(t)$.

Assumption \ref{hyp:Pitch} represents a necessary condition to bound roll angle estimation errors. 
 However, in practice, this assumption does not represent a substantial limitation because, in on-track motorsport, motorcycles pitch of few degrees (comprehensive of track grade and wheelie). 

\begin{theorem}
	\label{Theorem:DREBounds}
	Let Assumptions \ref{hyp:Sensors}-\ref{hyp:nu} hold and $S(t)$ be the solution to \eqref{eq:DRE_CF}. Then, there exist $\underline{s},\overline{s} > 0$ such that 
	\[
	\underline{s} I \preceq S(t) \preceq \overline{s} I \qquad \forall t \ge 0. 
	\]
	Moreover, there exist  $\rho > 0$ and a class-$\mathcal{K}$ function $\beta(\cdot)$ such that, for any $\|\hat{x}(0)-x(0)\| \le \rho$, the trajectories of \eqref{eq:EKF_CF} are bounded and $\limsup_{t \to \infty} \|\hat{q}(t)-q(t)\| \le \beta (\|\mathtt{col}({\nu}(t),w(t))\|_\infty)$. { To conclude, let  Assumption \ref{hyp:Pitch} hold. Then, system \eqref{eq:phi_av}-\eqref{eq:hatphi} solves Problem \ref{prob:LeanAngleEstimation}.}\hfill $\square$
\end{theorem}
Theorem \ref{Theorem:DREBounds} is proved in Appendix \ref{sec:ProofOfTheorem1}.

\section{Description of the proposed solution}
\label{sec:Description}

Section \ref{sec:Coordinated} aims to describe the novel concept of \textit{coordinated manoeuvre} modelling complex motion configurations in which the rider's gravity centre is not on the motorbike's plane of symmetry. With the function $\phi_{\text{av}}(\cdot,\cdot,\cdot)$ at hand, the vector $q_1$ is built through \eqref{eq:hat_theta_e_av}-\eqref{eq:q1} where the estimate $\hat{\xi}_e$ is obtained in Section \ref{sec:Pre} via the so called \textit{pre-filter}. Finally, Section \ref{sec:Rec Err Analysis} analyses the \textit{pre-filter} errors.

\subsection{Coordinated manoeuvres}
\label{sec:Coordinated}


Let $v^B:= T(\Theta) v$. 
Then, we define as \textit{coordinated manoeuvres} the set of dynamic states such that $v^B \equiv \mathtt{col}(\mathtt{v},0,0)$.  

\begin{lemma}
	\label{lemma:PhiAv}
	Consider $\phi_{\text{av}}(\cdot,\cdot,\cdot)$ defined in \eqref{eq:phi_av} and $a := T(\Theta)(\dot{v}-g)$, then  
	\begin{equation}
		\label{eq:phi_av2}
		\Theta = \mathtt{col}(\phi_{\text{av}}(a,\xi_e,\mathtt{g}),\gamma,\chi)
	\end{equation}
during coordinated manoeuvres verifying Assumption \ref{hyp:Non-Ballistic}. \hfill $\square$
\end{lemma}

The proof of Lemma \ref{lemma:PhiAv} is reported in Appendix \ref{sec:ProofPhiAv}.	

\begin{figure}
		\centering
		\includegraphics[clip, trim = 7cm 5cm 5cm 2cm,width=1\columnwidth]{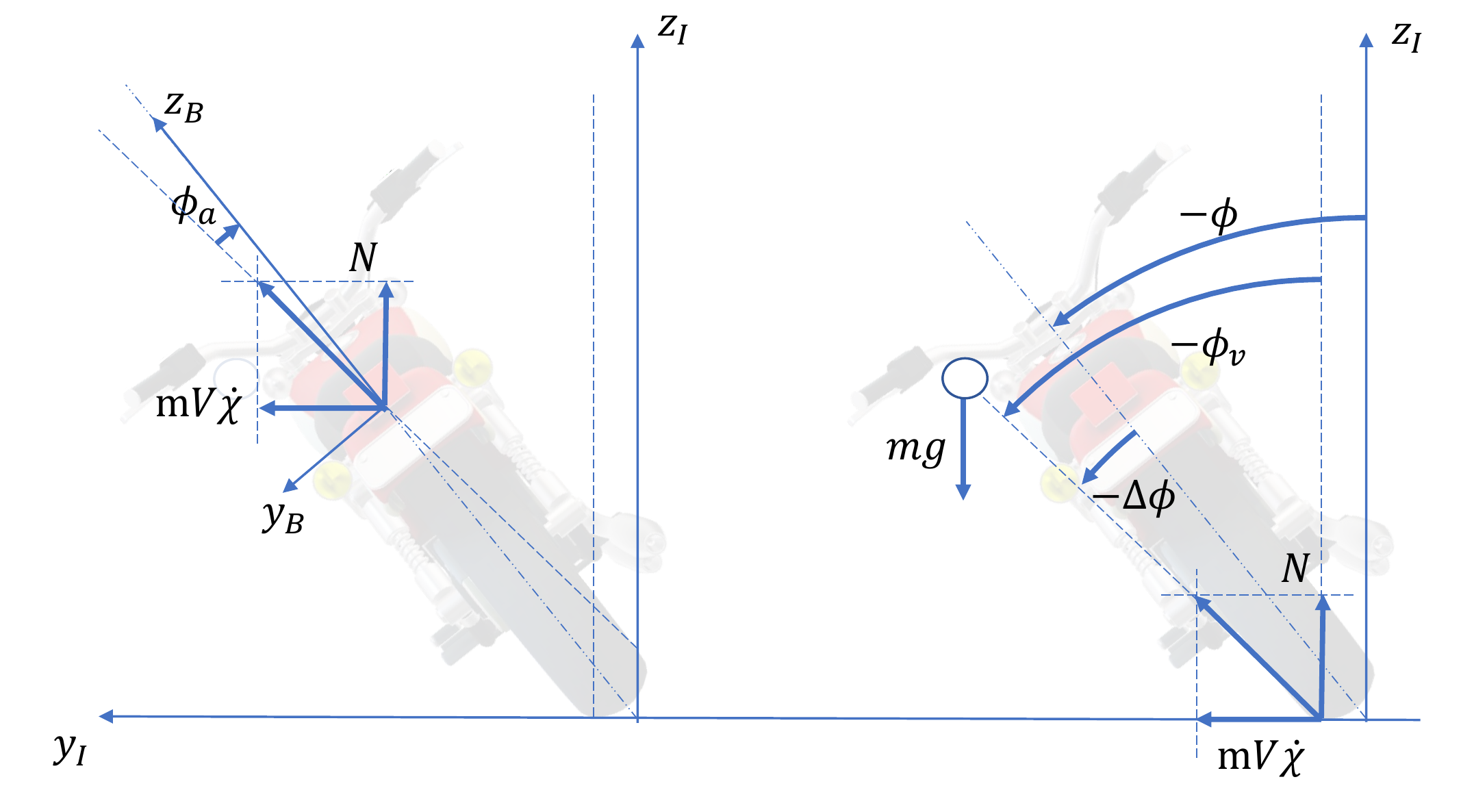}
		\caption{Equilibrium of forces and torques in a flat-coordinated turn assuming the pilot body shifts the gravity centre of the system out of the longitudinal symmetry axis.}
		\label{fig:CoordinatedTurn}
	\end{figure}

In the remaining of this section, we show how the \textit{coordinated manoeuvres} improves the roll angle estimation.
	Let $v := \mathtt{col}(v_x,v_y,v_z)$ and $V := \sqrt{v_x^2 + v_y^2}$. Concerning Figure \ref{fig:CoordinatedTurn}, define a \textit{flat-coordinated turn} as a coordinated manoeuvre performed under the  further constraints $\gamma,\dot{\gamma} = 0$. With these constraints at hand, the system composed of motorbike and  biker is at the equilibrium (translations and rotation) at
	\begin{equation}
		\label{eq:Phiv}
		\phi_v(\xi_e,\mathtt{g}) :=-  \tan^{-1}\left(\dfrac{V\dot{\chi}}{\mathtt{g}}\right).
	\end{equation}
\begin{remark}
		It is worth noting that $\phi_v$ denotes the roll angle that the complete system (motorbike + biker) negotiates to perform a coordinated turn. This angle corresponds to $\phi$ when the biker doesn't move his body out of the motorbike symmetry plane. The difference between $\phi$ and $\phi_v$ due to the tire size and the centre of gravity shift due to the rider movements during flat-coordinated turns, are well-known concepts, as pointed out in [\cite{Cossalter2002Motorcycle}, \S 4.1.2] and recalled in \cite{Boniolo2010Motorcycle}. However, to the author's knowledge, what follows represents the first effective way to compensate for this difference in the context of roll angle estimation.  \hfill $\square$
\end{remark}

	Equation \eqref{eq:SensorSuite} with $\gamma = 0$ and $\psi = \chi$ becomes
	\begin{equation}
		\label{eq:Inclinometer}
		a =  T_1(\phi) T_3(\chi) (\dot{v}-g)
	\end{equation}
where $T_i(s)$ denotes the  matrix associated with a rotation, of magnitude $s$, around the $i$-th axis.	Then, since $\dot{v} = V\dot{\chi}\mathtt{col}(-\sin\chi, \cos\chi, 0)$, and using $a = \mathtt{col}(a_x, a_y, a_z)$, the roll angle is found through \eqref{eq:Inclinometer} as
	\begin{equation}
		\label{eq:TanPhiA}
		\phi_a(a) :=\tan^{-1}\left({a_y}/{a_z}\right).
	\end{equation}
	As detailed in Figure \ref{fig:CoordinatedTurn}, the following equation holds
	\begin{equation}
		\label{eq:CoordinatedTurn}
		\begin{aligned}
			\phi =&\, \phi_v - \Delta\phi
		\end{aligned}
	\end{equation}
	where $\Delta\phi \in \mathbb{R}$ models the effects of the pilot displacement. We can exploit the body accelerations to correct $\phi_v$ in flat-coordinated turns. Indeed, through basic geometric arguments, it results $\phi_a =  -\Delta \phi$. Then, use \eqref{eq:Phiv}, \eqref{eq:TanPhiA} and \eqref{eq:CoordinatedTurn} to verify
		\begin{equation}
		\label{eq:CoordinatedTurnC}
		\begin{aligned}
			\left.\phi_{\text{av}}(a,\xi_e,\mathtt{g})\right|_{\gamma , \dot{\gamma} = 0} =&\, \phi_a(a)+\phi_v(\xi_e,\mathtt{g}). \nonumber
		\end{aligned}
	\end{equation}

\subsection{Pre-Filter}
\label{sec:Pre}

The \textit{pre-filter}, representing the subsystem providing $\hat{\xi}_e$, relies of two subsystems, \textit{i.e.}, the GNSS reconstructor, estimating $v$ and $\dot{v}$, and the continuous-course estimator, providing $\hat{\chi}$.  

\subsubsection{GNSS Reconstructor}
Let $v_e := \mathtt{col}(v,\dot{v})$ and define $f_{\xi_e}\,:\,\mathbb{R}^6\to \mathbb{R}^6$ such that 
\begin{equation}\label{eq:xie}
	\xi_e = f_{\xi_e}(v_e)
	\end{equation} 
with $f_{\xi_e}(\cdot)$ detailed in Appendix \ref{sec:fxie}. This work adopts a second-degree polynomial signal reconstructor that interpolates the $n\in\mathbb{N}$ most recent samples of $y_s$ to estimate $v_e$, namely via $\hat{v}_e$. Moreover, The reconstructor extrapolates the signal values, along the next $\tau$-long interval. Finally, we impose
\begin{equation}
	\label{eq:fhatxie}
		\begin{aligned}
			\hat{\xi}_e =&\, f_{\xi_e}(\hat{v}_e).
		\end{aligned}
\end{equation}

More precisely, let $k:=\lfloor{t}/{\tau}\rfloor$  be the maximum integer not greater than ${t}/{\tau}$ and
 introduce
\begin{equation}
	\label{eq:u}
	u(k,\tilde{t}) := c_2(k)  \tilde{t}^2+c_1(k) \tilde{t}+c_0(k)\quad \tilde{t} \in [-(n-1)\tau, \,\tau)
\end{equation}
where $c_2(\cdot)$, $c_1(\cdot)$, and $c_0(\cdot) \,:\,\mathbb{N} \to \mathbb{R}^3$. 
\\
\begin{remark}
	Function $u(\cdot,\cdot)$ is a vector-valued $2^{\text{nd}}$-degree polynomial describing a trajectory with a piecewise-constant jerk in the time interval $((k-n+1)\tau, \,(k+1)\tau)$. In practice, we define $c_2(\cdot)$, $c_1(\cdot)$, and $c_0(\cdot)$ as a jerk,  acceleration, and speed at time $t = k\tau$ that best describe the last $n$ GNSS data. Finally, we use the same coefficients to preview inertial speed and acceleration within the next $\tau$-long time window.\hfill $\square$
\end{remark}

Use \eqref{eq:u} to define $$\tilde{u}(k,s) = u(k,s-k\tau)-v(s) = u(k,s-k\tau)-h_v(\xi(s))$$ for all $s \in [(k-n+1)\tau, \,(k+1)\tau)$. 
Then, the GNSS receiver provides at time $t=k\tau$
\begin{equation}
	\label{eq:ys}
	y_s(k\tau) = u(k,0)-\tilde{u}(k,k\tau) +\nu_{s}(k\tau).
\end{equation}
In the following, we present an algorithm elaborating the $n$ most recent samples $y_s(k\tau), \dots, y_s((k-n+1)\tau)$ to provide estimations for $c_2(k)$, $c_1(k)$, and $c_0(k)$, namely $\hat{c}_{2k}$, $\hat{c}_{1k}$, and $\hat{c}_{0k}$, respectively. Consequently, we propose to reconstruct $u(\cdot,\cdot)$ as
\begin{equation}
	\label{eq:hatu}
	\hat{u}(k,\tilde{t}) := \hat{c}_{2k}  \tilde{t}^2+\hat{c}_{1k} \tilde{t}+\hat{c}_{0k} \nonumber
\end{equation}
and to use it to approximate $v(t)$ and $\dot{v}(t)$ by
\begin{equation} 
	\label{eq:hat_v_vdot}
	\resizebox{0.89\columnwidth}{!}{$
\begin{aligned}
	\hat{v}(t) =&\, \hat{u}(k,t-k\tau) = \hat{c}_{2k}  (t-k\tau)^2+\hat{c}_{1k} (t-k\tau)+\hat{c}_{0k}\\
	\hat{\dot{v}}(t) =& \left.\dfrac{\partial \hat{u}(k,\tilde{t})}{\partial \tilde{t}}\right|_{\tilde{t} = t-k\tau} 
		=2\hat{c}_{2k}  (t-k\tau)+\hat{c}_{1k}
\end{aligned}
$}
\end{equation}
where $t \in [k\tau,(k+1)\tau)$. Let $R_{\nu_s} =R_{\nu_s}^\top \in \mathbb{R}^{3\times 3}$ with $R_{\nu_s}\succ 0$. Then, we determine $\hat{\zeta}_k := \mathtt{col}(\hat{c}_{2k},\hat{c}_{1k},\hat{c}_{0k})$ through
\begin{equation}
	\label{eq:hatck_prob}
	\hat{\zeta}_k=\underset{x \in {\mathbb R}^9}{\operatorname{argmin}}\{(y_k-Cx)^{\top}(I_n\otimes R_{\nu_s})^{-1}(y_k-Cx)\}
\end{equation}
where $y_k = \mathtt{col}(y_s(k\tau), \dots, y_s((k-n+1)\tau))$,
\[
C = \left[\begin{array}{ccc}
	0 & 0 & 1\\
	\tau^2 & -\tau & 1\\
	\vdots & \vdots & \vdots\\
	(n-1)^2\tau^2 & -(n-1)\tau & 1
\end{array}\right]\otimes I_3,
\]
and $n \ge 3$.
The solution of \eqref{eq:hatck_prob} is (see \cite{barshalom2001estimation})
\begin{equation}
	\label{eq:hatck}
\hat{\zeta}_k = K_s y_k
\end{equation}
with 
\begin{equation}
\label{eq:Least Squares}
K_s=(C^{\top}(I_n\otimes R_{\nu_s})^{-1}C)^{-1}C^{\top}(I_n\otimes R_{\nu_s})^{-1}.
\end{equation}

\subsubsection{Continuous-Course Estimator}
\label{sec:hatchi}
The computation of $\hat{\chi}$ from inertial speeds $\hat{v}_x$ and $\hat{v}_y$, made through whether $\tan^{-1}(\hat{v}_y/\hat{v}_x)$ or $\mathtt{atan}_2(\hat{v}_y,\hat{v}_x)$, is prone to discontinuities, which could induce wrong roll angle estimations. 

This section proposes a novel continuous map $t \mapsto \hat{\chi}(t)$ that solve this issue. With reference to \eqref{eq:hat_v_vdot}, remember that $\hat{v} := \mathtt{col}(\hat{v}_x, \hat{v}_y, \hat{v}_z)$ and $\hat{\dot{v}} := \mathtt{col}(\hat{\dot{v}}_x, \hat{\dot{v}}_y, \hat{\dot{v}}_z)$ and define 
\begin{equation*}
	\resizebox{\columnwidth}{!}{$
\begin{aligned}
	\mathcal{T}_k^{+}=&\, \{t \in[k\tau,(k+1)\tau)\,:\,\hat{v}_y(t) = 0,\,\hat{\dot{v}}_y(t) < 0,\,\hat{v}_x(t) < 0\}\\
	\mathcal{T}_k^{-}=&\, \{t \in[k\tau,(k+1)\tau)\,:\,\hat{v}_y(t) = 0,\,\hat{\dot{v}}_y(t) > 0,\,\hat{v}_x(t) < 0\},\\
\end{aligned}	
$}
\end{equation*} 
 where $\mathcal{T}_k^{+}$ and $\mathcal{T}_k^{-}$ are analytically found thanks to \eqref{eq:hat_v_vdot} being polynomial functions of time. Now, define $\mathcal{T}^{+} = \bigcup_{k \in \mathbb{N}} \mathcal{T}_k^{+}$ and $\mathcal{T}^{-} = \bigcup_{k \in \mathbb{N}} \mathcal{T}_k^{-}$ and use them to feed the lap counters \eqref{eq:LapCounter}. Then, adopt $f_\chi(\cdot,\cdot,\cdot)$ detailed in Appendix \ref{sec:fxie} to introduce 
\begin{equation}
	\label{eq:fchi_new2}
	\begin{aligned}
	\hat{\chi}&\,=f_{\chi}(\hat{v},N_{+},N_{-}). 
		\end{aligned}
\end{equation}
\begin{lemma}
	\label{lemma:ContinuityHatChi}
	Consider \eqref{eq:fchi_new2}, then
	 \[
	 \,\,\,\,\qquad\qquad\hat{\chi}(t)= f_{\chi}(\hat{v}(t),N_{+}(t),N_{-}(t)) \in \mathcal{C}^1. \,\,\,\,\qquad\qquad \square 
	 \]	  
 \end{lemma}Appendix \ref{sec:ProofContinuityHatChi} details the proof of Lemma \ref{lemma:ContinuityHatChi}.

\subsection{Error Boundedness}
\label{sec:Rec Err Analysis}

In this section we investigate the error $\nu = q_1-q$. 
	To this end, let $\tilde{\xi}_e := \hat{\xi}_e-\xi_e$ and 
	\begin{equation} 
		\label{eq:Theta_e_av}
		\Theta_{\text{av}}:= f_{\Theta}(a,\xi_e,\mathtt{g}), 
	\end{equation} 
	define $\tilde{\mathtt{g}} = \hat{\mathtt{g}}-\mathtt{g}$, and remember  $q_1=\bar{h}(\hat{\Theta}_{\text{av}})$. Then,
	\begin{equation}
		\label{eq:q_1}
		\begin{aligned}
		\nu = \,& \bar{h}(\hat{\Theta}_{\text{av}})-\bar{h}(\Theta_{\text{av}})+\bar{h}(\Theta_{\text{av}})-q =\\ 
		\,= \,&\nu_\nu(t,\nu_a,\tilde{\xi}_e,\tilde{\mathtt{g}}) + \nu_{\text{m}}(t)
	\end{aligned}
	\end{equation}
	where
	\begin{equation}
		\label{eq:nu_nu_m}
		\begin{aligned}
						\nu_{\text{m}}(t) &\,:= 	\bar{h}(f_{\Theta}(a(t),\xi_e(t),\mathtt{g}))-q(t)\\
			\nu_\nu(t,\nu_a,\tilde{\xi}_e,\tilde{\mathtt{g}}) 
			&\,:=  \bar{h}(f_{\Theta}(y_a(t), \hat{\xi}_e(t),\hat{\mathtt{g}}))\\
			&-\bar{h}(f_{\Theta}(y_a(t)-\nu_a, \hat{\xi}_e(t)-\tilde{\xi}_e, \hat{\mathtt{g}}-\tilde{\mathtt{g}})).
		\end{aligned}
	\end{equation}
	The first error contribution, \textit{i.e.}, $\nu_\text{m}$, embeds the errors due to the model mismatch, \textit{i.e.}, the difference between the actual motorbike evolution and a coordinated manoeuvre. In detail, $\nu_\text{m}$ is highly dependent on the rider's driving style, mainly due to wheelies and drifts. Thus, with a particular focus on applications like the Grand Prix motorcycle racing, the error $\nu_\text{m}$ is usually negligible except during tail-wagging or corner entries.  However, these represent short-duration driving phases in which the side-slip remains bounded. Thus, we formalise this through the following assumption.
		\begin{assumption}[Boundedness of $\nu_{\text{m}}$]
			\label{hyp:num}
			There exists $\overline{\nu}_{\text{m}} > 0$ such that $\|\nu_{\text{m}}(t)\|_\infty < \overline{\nu}_{\text{m}}$. \hfill $\square$
		\end{assumption}

	As for  $\nu_\nu$, it describes the uncertainties induced by the sensor inaccuracy plus those introduced by the estimator of $\xi_e$. Note that, since $\bar{h}(\Theta)$ and all its derivatives are Lipschitz and bounded  for all $\Theta \in \mathbb{R}^3$, there exists $\overline{\nu}_\nu > 0$ such that $\|\nu_\nu(t,\nu_a,\tilde{\xi}_e,\tilde{\mathtt{g}})\|_\infty < \overline{\nu}_\nu$.
	 \section{Experimental results}

This section presents the results of numerical and field tests. While the former were conducted to check the expected theoretical behaviour, the latter were performed to assess the applicability of the proposed estimation scheme. 

The field test were conducted by AvioRace \cite{AvioRace}, a provider of electronics specialised in motor-sport applications. Since the parties agreed on a data protection policy, sensible data collected during field test are shown without scale. The lack of quantitative evaluations is compensated in Section \ref{sec:Simulation} where the algorithm is tested in a realistic synthetic environment.

\subsection{Pre-Filter Performance Analysis}
\label{sec:PreFilterPerformance}
The boundedness of $\|\nu_\nu(t)\|_\infty$ demonstrated in section \ref{sec:Rec Err Analysis} is necessary for the observer stability proof. In contrast, $\|\nu_\nu(t)\|_\infty$ is too conservative in assessing the pre-filter's performance. Therefore, this section reports a stochastic description of $\nu_\nu$, which is also used to tune \eqref{eq:BasicEKF_CF}.
To this aim, we propose the following process:
\begin{enumerate}
	\item we introduce a stochastic model of sensor noise
\item we estimate how GNSS measurement errors and the GNSS reconstructor impact $\tilde{v}:=\hat{v}-v$ and $\tilde{\dot{v}}:=\hat{\dot{v}}-\dot{v}$
\item we rely on results of point 1) to describe $\tilde{\xi}:= \hat{\xi}-\xi$ and $\tilde{\dot{\xi}}:=\hat{\dot{\xi}}-\dot{\xi}$
\item we use results from point 2) to characterize $\nu_\nu$.
\end{enumerate}
\subsubsection{Stochastic Description of Sensor Noise}

Let  $w_{\#}\,:\,\mathbb{R}\to \mathbb{R}^9$, for all $\#\in\{a,g,s\}$ a stationary stochastic process. Then, in agreement with \cite{Farrell2022IMU}, this paper models the measurement errors appearing in \eqref{eq:SensorSuite} as 
\begin{equation}
	\label{eq:GenSensorModel}
	\begin{aligned}
		\dot{b}_{\#} = &\, A_{\#} b_{\#} + B_{\#}\, w_{\#}(t) && b_{\#}(0) = b_{{\#}0} \\
		\nu_{\#} =& \,C_{\#}\, b_{\#} + D_{\#}\, w_{\#}(t) &&{\#}\in \{a,g,s\} 
	\end{aligned}
\end{equation} 
with $b_{\#} := \mathtt{col}(\bar{b}_{\#},z_{\#})$, $\bar{b}_{\#},z_{\#}\in\mathbb{R}^3$,
\[
\begin{aligned}
	A_{\#} =&\,  \left[\begin{array}{cc}
		0 & 0 \\ 0 & -\tau_{\#}^{-1}
	\end{array}\right]\otimes I_{3}&&B_{\#} =  \left[\begin{array}{ccc}
		1 & 0 & 0\\
		0 & 1  & 0
	\end{array}\right] \otimes I_{3}\\
	C_{\#} =&\, \left[\begin{array}{cc}
		\,1 & \quad 1\quad
	\end{array}\right] \otimes I_{3} && D_{\#} = \left[\begin{array}{ccc}
		0 & 0 & 1\\
	\end{array}\right] \otimes  I_{3},
\end{aligned}
\]
and where $\tau_{\#} > 0$ and $\otimes$ denotes the Kronecker product. 

Usually the Power Spectral Density (PSD) of $w_{\#}$ is constant within the sensor sampling frequency. Hence, for analysis purposes,  $w_{\#}$ can be seen as a white noise by assuming its PSD constant for all the frequencies. Therefore, we assume ${\rm E}[w_{\#}(t)]=0$, and ${\rm E}[w_{\#}(t) w_{\#}^\top(\tau)] = R_{\#}  \delta(t-\tau)$ with $R_{\#} \succ 0$ and block diagonal, for all ${\#} \in \{a,g,s\}$. Lastly, enforce $b_{s0} = 0$ and $B w_s(t) = 0$ for all $t \ge 0$.

Consequently, $\bar{b}_{\#}$ models a biased random walk while $z_{\#}$ represents a coloured noise (with a time constant $\tau_{\#}$). In particular, only a 0-mean white noise affects the GNSS measurement.
\subsubsection{Analysis of $\tilde{v}$ and $\tilde{\dot{v}}$} 
Introduce 
\[
\begin{aligned}
	\zeta_k:= &\,\mathtt{col}(c_{2}(k),c_1(k),c_0(k))	\\
	\tilde{u}_k := &\,\mathtt{col}(\tilde{u}(k,k\tau),\dots,\tilde{u}(k,(k-n+1)\tau))\\
	\nu_k := &\, \mathtt{col}(\nu_s(k\tau),\dots,\nu_s((k-n+1)\tau))\\
\end{aligned}
\]
write $y_k = C\zeta_k -\tilde{u}_k+\nu_k$ and use $R_{\nu_s}:=DR_sD^\top$  into \eqref{eq:Least Squares}.
\begin{remark}
	Usually $R_{\nu_s}=\mathtt{diag}(\sigma_x^2,\sigma_y^2,\sigma_z^2)$ where $\sigma_{x}$, $\sigma_{y}$, and $\sigma_z$ are known and correspond to figures of merit of GNSS, known as \textit{User Range Rate Error}.\hfill $\square$
\end{remark}

Besides, use \eqref{eq:hatck} to write
\begin{equation}
	\label{eq:tilde_ck}
\tilde{\zeta}_k := \hat{\zeta}_k-\zeta_k = K_s y_k - \zeta_k = K_s(\nu_k- \tilde{u}_k).
\end{equation}
Introduce 
\[
\Phi(\tilde{t}) := \left[\begin{array}{ccc} 
\tilde{t}^2& \tilde{t}& 1\\
2\,\tilde{t}& 1& 0\end{array}
\right]\otimes I_3,
\] define $\hat{v}_e := \mathtt{col}(\hat{v},\hat{\dot{v}})$, and  rewrite \eqref{eq:hat_v_vdot} as
\begin{equation}
	\label{eq:hatv_hatvdot}
\hat{v}_e = \Phi(t-k\tau) \hat{\zeta}_k \qquad t \in [k\tau,(k+1)\tau).
\end{equation}
Introduce ${v}_e := \mathtt{col}({v},{\dot{v}})$ and $\tilde{u}_e := \mathtt{col}(\tilde{u},\dot{\tilde{u}})$ and use \eqref{eq:tilde_ck} and \eqref{eq:hatv_hatvdot} to calculate the estimation errors 
\begin{equation}
	\label{eq:tilde_ve}
	\begin{aligned}
	\tilde{v}_e(t)  = &\,\hat{v}_e(t)-{v}_e(t) \\
= &\,
\Phi(t-k\tau) K_s(\nu_k-\tilde{u}_k) + \tilde{u}_e(t)
	\end{aligned}
\end{equation}
for all $t \in [k\tau,\,(k+1)\tau)$. Use \eqref{eq:tilde_ve} to compute
\begin{equation}
		\label{eq:Eve}
	\begin{aligned}
		{\rm E}[\tilde{v}_e(t)] = &\,	{\rm E}[\Phi(t-k\tau) K_s(\nu_k-\tilde{u}_k) + \tilde{u}_e(t)]\\
		 = &\,	-\Phi(t-k\tau) K_s {\rm E}[\tilde{u}_k]+	{\rm E}[\tilde{u}_e(t)] \nonumber
	\end{aligned}
\end{equation}
in which we have exploited ${\rm E}[\nu_k] = 0$.

\begin{remark}
The quantity ${\rm E}[\tilde{v}_e]$ represents the expected velocity and acceleration estimate errors. Roughly, ${\rm E}[\tilde{v}_e(t)] = 0$ because it can be demonstrated that $\tilde{u}_k=\tilde{u}_e(t)=0$ if the motorbike travels at a statistically constant jerk into $n\tau$-long time intervals. We assumed $\tilde{v}_e(t)$ as an ergodic process to test if the expected value is near zero in a real-world scenario. Therefore, we computed the time average by using the GNSS samples reported in the experimental tests of Section \ref{sec:FieldTest}. The results show that, in practice, ${\rm E}[\tilde{v}_e] \approx 0$. 
\hfill $\square$
\end{remark}

Finally, compute the covariance of the estimation error
\begin{equation}
\nonumber
		R_v(t):=	\, 	{\rm E}[(\tilde{v}_e(t)-{\rm E}[\tilde{v}_e(t)])(\tilde{v}_e(t)-{\rm E}[\tilde{v}_e(t)])^\top]
\end{equation}
by applying \eqref{eq:tilde_ve} and exploiting the assumptions ${\rm E}[\nu_k\nu_k^\top]=I_n\otimes R_{\nu_s}$, ${\rm E}[\nu_k\tilde{u}_k^\top]=0$, ${\rm E}[\nu_k\tilde{u}_e^\top(t)]=0$, ${\rm E}[\tilde{u}_k\tilde{u}_k^\top]=0$, ${\rm E}[\tilde{u}_k\tilde{u}_e^\top(t)]=0$, ${\rm E}[\tilde{v}_e(t)] = 0$ and ${\rm E}[\nu_k]=0$.
After some algebra, it results to be 
\begin{equation}
	\label{eq:Rv} 
		R_v(t)=	\Phi(t-k\tau) K_s( I_n\otimes R_{\nu_s}) K_s^\top \Phi^\top(t-k\tau).
\end{equation}
In conclusion, it is worth noting that 
$R_v$ is bounded because, 
since $\Phi(t-k\tau)$ is a polynomial function of $t$, there exists a finite $\overline{\Phi}(\tau) > 0$ such that $\|\Phi(t-k\tau)\| \le \overline{\Phi}(\tau)$ for all $t \in [k\tau,\,(k+1)\tau)$ and for any $k \in \mathbb{N}$.
In particular, as a conservative approach, the upper bound $\bar{R}_{vk}$ of  the covariance can be chosen accordingly to
\begin{equation}
	\label{eq:bar_R_vk}
	\bar{R}_{vk}=R_v((k+1)\tau).
\end{equation}
 
\subsubsection{Analysis of $\tilde{\xi}$ and $\tilde{\dot{\xi}}$} 

Use \eqref{eq:fhatxie} and  $\tilde{v}_e = \hat{v}_e-v_e$ to compute
\begin{equation}
	\label{eq:tildexie}
\begin{aligned}
\tilde{\xi}_e =  &\,
f_{\xi_e}(\hat{v}_e) - f_{\xi_e}(\hat{v}_e-\tilde{v}_e) \approx  J_{\xi_e}(\hat{v}_e) \tilde{v}_e \nonumber
\end{aligned}
\end{equation}
where $J_{\xi_e}(\hat{v}_e) := \partial f_{\xi_e}(\hat{v}_e)/\partial \hat{v}_e$.  
Then, the expected value is ${\rm E}[\tilde{\xi}_e] \approx J_{\xi_e}(\hat{v}_e) {\rm E}[\tilde{v}_e]$. To conclude,
 we exploit \eqref{eq:Rv} to compute
\begin{equation}
	\begin{aligned}
R_{\xi_{e}}(t)=&\,{\rm E}[(\tilde{\xi}_{e}(t)-{\rm E}[\tilde{\xi}_e(t)])(\tilde{\xi}_{e}(t)-{\rm E}[\tilde{\xi}_e(t)])^\top]\\
\approx &\, J_{\xi_e}(\hat{v}_e(t)) R_{v}(t) J_{\xi_e}^\top(\hat{v}_e(t))  -{\rm E}[\tilde{\xi}_{e}(t)] {\rm E}[\tilde{\xi}_e^\top(t)]
\end{aligned}
\nonumber
\end{equation}
where it is worth observing that $\|J_{{\xi}_e}(\cdot)\|$ is bounded under the following Assumption. 

	\begin{assumption}[Motorbike Speed]
	\label{hyp:Speed}
	The inertial speed $v(t)$ is a Lipschitz continuous function. Moreover, there exist $\underline{v}, \overline{v} > 0$ such that $\sqrt{v_x^2(t)+v_y^2(t)}>\underline{v}$ and $\mathtt{v}(t) < \overline{v}$ for all $t \ge 0$. \hfill $\square$
\end{assumption}
As for Assumption \ref{hyp:Speed}, Eq. \eqref{eq:v_aero} becomes bijective if $\mathtt{v}(t) = \sqrt{v^\top(t) v(t)} >\sqrt{v_x^2(t)+v_y^2(t)}$ is strictly positive. 
Moreover, the Lipschitz continuity of $v(t)$ ensures that $\dot{v}(t)$ and $y_a(t)$ are bounded. It is worth noting that Assumption \ref{hyp:Speed} does not represents a constraint because, in practice, motorbikes are power- and force-limited systems for which the assumption of a Lipschitz continuous speed represents a matter of fact. Indeed, engines deliver bounded powers and torques while tires transfer bounded traction/braking forces to the ground, thus limiting accelerations.

\subsubsection{Analysis of $\nu_\nu$}

Exploit \eqref{eq:hat_theta_e_av} and \eqref{eq:Theta_e_av} to calculate
\begin{equation}
	\tilde{\Theta}_{\text{av}}:= \hat{\Theta}_{\text{av}}-\Theta_{\text{av}} \approx J_{\Theta}(t) \mathtt{col}(\nu_a, \tilde{\xi}_e,\tilde{\mathtt{g}}) \nonumber
\end{equation}
in which $ J_{\Theta}(t) = \partial f_{\Theta}(a,\xi_e,\mathtt{g}) / \partial \mathtt{col}(a,\xi_e,\mathtt{g})$ evaluated at $a=y_a(t)$, $\xi_e=\hat{\xi}_e(t)$, and $\mathtt{g} = \hat{\mathtt{g}}$. The expected value of $\tilde{\Theta}_{\text{av}}$ is
\begin{equation}
	\begin{aligned}
			{\rm E}[\tilde{\Theta}_{\text{av}}(t)] \approx &\, J_{\Theta}(t) \mathtt{col}({\rm E}[\nu_a], {\rm E}[\tilde{\xi}_e(t)], {\rm E}[\tilde{\mathtt{g}}])   \\
	\end{aligned}\nonumber
\end{equation}
 while the covariance of $\tilde{\Theta}_{\text{av}}$ is
\begin{equation}
\begin{aligned}
R_{\Theta}(t) :=&\, {\rm E}[(\tilde{\Theta}_{\text{av}}(t)-{\rm E}[\tilde{\Theta}_{\text{av}}(t)])(\tilde{\Theta}_{\text{av}}(t)-{\rm E}[\tilde{\Theta}_{\text{av}}(t)])^\top] \\
 \approx &\, J_{\Theta}(t) \mathtt{blkdiag}(R_{\nu_a}, R_{{\xi}_e}(t), 0) J^\top_{\Theta}(t) - {\rm E}^2[\tilde{\Theta}_{\text{av}}(t)]
\end{aligned}
\nonumber
\end{equation}
in which we exploited ${\rm E}[\tilde{\mathtt{g}}^2] = 0$. It is worth noting that $\|J_\Theta(t)\|_\infty$ is bounded. Finally, with the same steps outlined in Eq. (37)-(43) of \cite{Markley2004QuatEstim}, we obtain
\begin{equation}
	\label{eq:Rnu_singular}
		\begin{aligned}
			{\rm E}[\nu_\nu\,\nu_\nu^\top] =  M(q_1(t))R_\Theta (t)M^\top(q_1(t)).		
		\end{aligned}
\end{equation}

\subsection{EKF Tuning Guidelines}
\label{sec:EKFTuning}
The tunable quantities of \eqref{eq:BasicEKF_CF} are $\lambda$, $Q$, and $R(t)$. In this section, we exploit results of Section \ref{sec:PreFilterPerformance} to design $R(t)$. Besides, we propose to exploit \eqref{eq:SensorSuite} and \eqref{eq:dotTheta}, and the stability arguments detailed in the proof of Theorem \ref{Theorem:DREBounds} to define $Q$. The tuning of EKF matrices can be divided in two parts: the selection of intra-matrix weights and the selection of inter-matrix weights. The idea is that the hardest part of the tuning, \textit{i.e.} the selection of intra-weights for $Q$ and $R(t)$ (which are not necessarily purely diagonal), is made through an automatic computational procedure. Consequently, we leave $\lambda$, representing the inter-matrix weight, as single scalar hand-tuned parameter.

Markley and Pittelkau demonstrated in \cite{Markley2004QuatEstim} and \cite{Pittelkau2003QuatAttitude} that ${\rm E}[\nu_\nu(t)\,\nu_\nu^\top(t)]$ computed in \eqref{eq:Rnu_singular} is singular because of $\bar{h}(\cdot)$, which enforces the constraint of unitary norm. Consequently, in agreement with \cite{Markley2004QuatEstim} and \cite{Pittelkau2003QuatAttitude}, we avoid singularities by introducing $\beta_R (\cdot)\,:\,\mathbb{R}\to \mathbb{R}$ and take
\begin{equation}
	\label{eq:Rnu}
			R(t) := \beta_R^2(t) q_1(t) q_1^\top(t)+  M(q_1(t))\bar{R}_\Theta (t)M^\top(q_1(t))		
\end{equation}
where $\bar{R}_\Theta (t)$ represents ${R}_\Theta (t)$ in the worst condition \eqref{eq:bar_R_vk}. 
As for $\beta_R(t)$, we set $\beta_R(t) = (\underline{\sigma}(R_\Theta(t))+\overline{\sigma}(R_\Theta(t)))/2$, where $\underline{\sigma}(R_\Theta)$ and $\overline{\sigma}(R_\Theta)$ are the smallest and largest singular value of $R_\Theta$. This choice assures both the non-singularity and the well-conditioning of $R(t)$, important for the computation of its inverse.

As for $Q$, we propose the following setting. In agreement with \eqref{eq:GenSensorModel}, let $w_b$, $w_z$, and $w_\nu \in \mathbb{R}^3$ such that $w_g= \mathtt{col}(w_b,w_z,w_\nu)$. Define ${\rm E}[w_b w_b^\top]$, ${\rm E}[w_z w_z^\top]$, and ${\rm E}[w_\nu w_\nu^\top]$ such that
$R_g = \mathtt{blkdiag}({\rm E}[w_b w_b^\top],\,{\rm E}[w_z w_z^\top],\,{\rm E}[w_\nu w_\nu^\top])$.  Introduce ${\rm E}[z_g z_g^\top] :=  (0.4365)^2 \tau_g^2 {\rm E}[w_z w_z^\top]$ corresponding to the covariance of $z_g$ evaluated at the worst Allan power spectral density [\cite{Farrell2022IMU}, Eq. (37)]. Then, introduce $ \beta_Q, \epsilon > 0$ and write
\[
Q = \mathtt{blkdiag}( {\rm E}[w_b w_b^\top],\epsilon I,{\rm E}[w_\nu w_\nu^\top]+{\rm E}[z_g z_g^\top], \beta_Q^2 ).
\]
Inspired by the same arguments adopted for $\beta_R$, the degree of freedom $\beta_Q$ is set as the average of the smallest and largest singular values of ${\rm E}[w_\nu w_\nu^\top]+{\rm E}[z_g z_g^\top]$, \textit{i.e.},
\[
\resizebox{\columnwidth}{!}{$
\beta_Q = \dfrac{1}{2}\left(\underline{\sigma}({\rm E}[w_\nu w_\nu^\top]+{\rm E}[z_g z_g^\top]) + \overline{\sigma}({\rm E}[w_\nu w_\nu^\top]+{\rm E}[z_g z_g^\top])\right).$}
\]

Since the parameter $\epsilon$ is only necessary for the observer stability, \textit{i.e.}, to guarantee that $\hat{q}=0$ does not belong to any forward invariant set for \eqref{eq:BasicEKF_CF},  we set $\epsilon \ll \underline{\sigma}({\rm E}[w_\nu w_\nu^\top]+{\rm E}[z_g z_g^\top])$ to let its contribution be relevant only for $\|\hat{q}\|\approx 0$.
				
Finally, $\lambda$ regulates the magnitude $Q$ relative to $R(t)$ and, consequently, the behaviour of $S(t)$ and $K(t):= S^{-1}(t)H^\top R^{-1}(t)$. It is a common fact that increasing the ratio between $Q$ and $R(t)$ increases $\overline{\sigma}(K(t))$ and thus makes EKFs have a shorter transient but a higher noise sensitivity. Therefore, in practice, the best compromise is found through an on-field trial and error procedure on $\lambda$.
		
\subsection{Field test}
\label{sec:FieldTest}
\begin{figure}
	\centering
	\includegraphics[clip, trim= 0cm 0cm 0cm 0cm, width=\columnwidth]{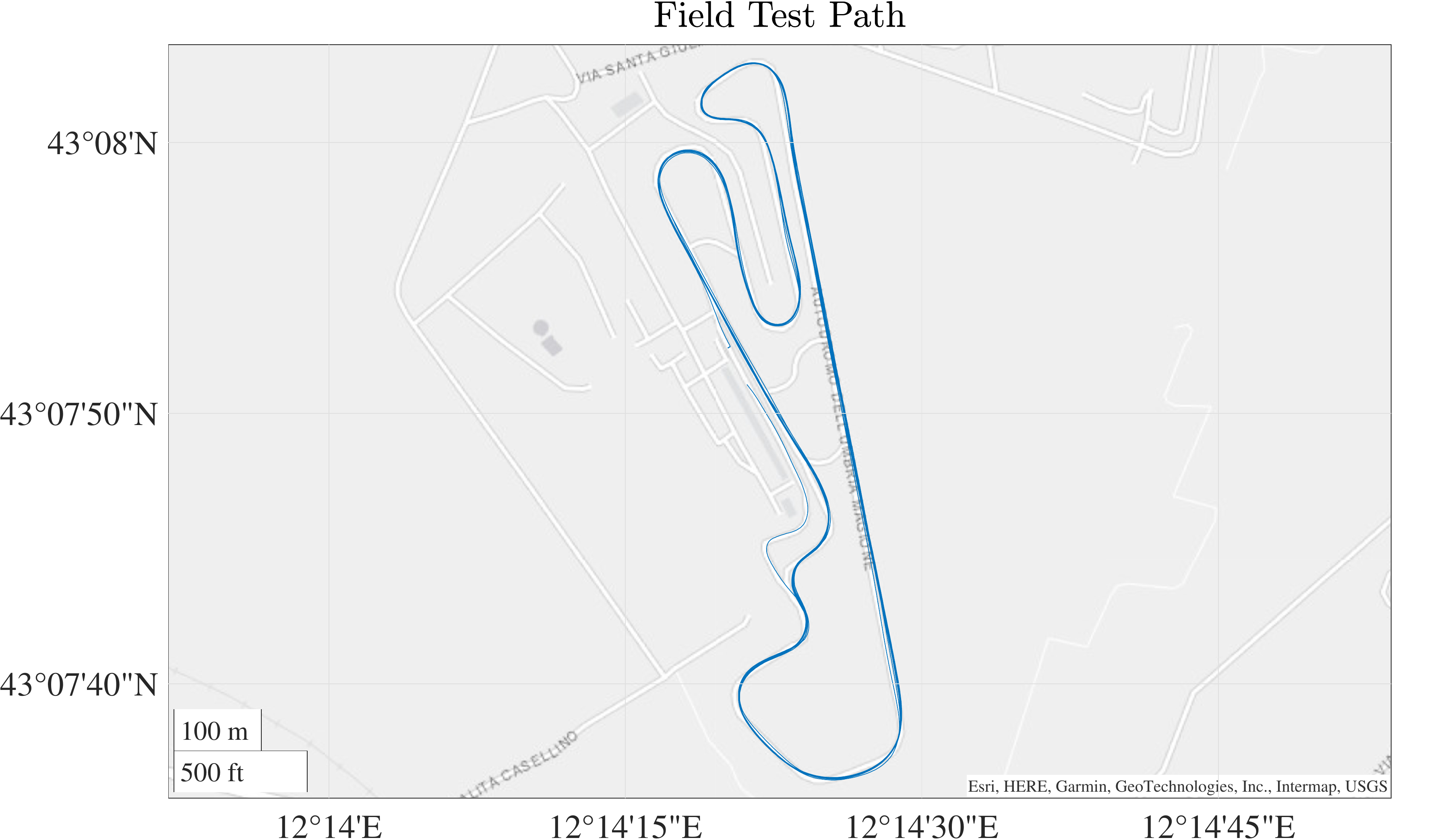}
	\caption{Autodromo Nazionale dell'Umbria "Mario Umberto Borzacchini". Circuit and race path sensed by the on-board GNSS receiver. The path was travelled clockwise.}
	\label{fig:Umbria}
\end{figure}

\begin{figure}
	\centering
	\begin{subfigure}[t]{0.55\columnwidth}
		\centering
			\includegraphics[clip, trim= 0.2cm 0 0.2cm 0, width=\textwidth]{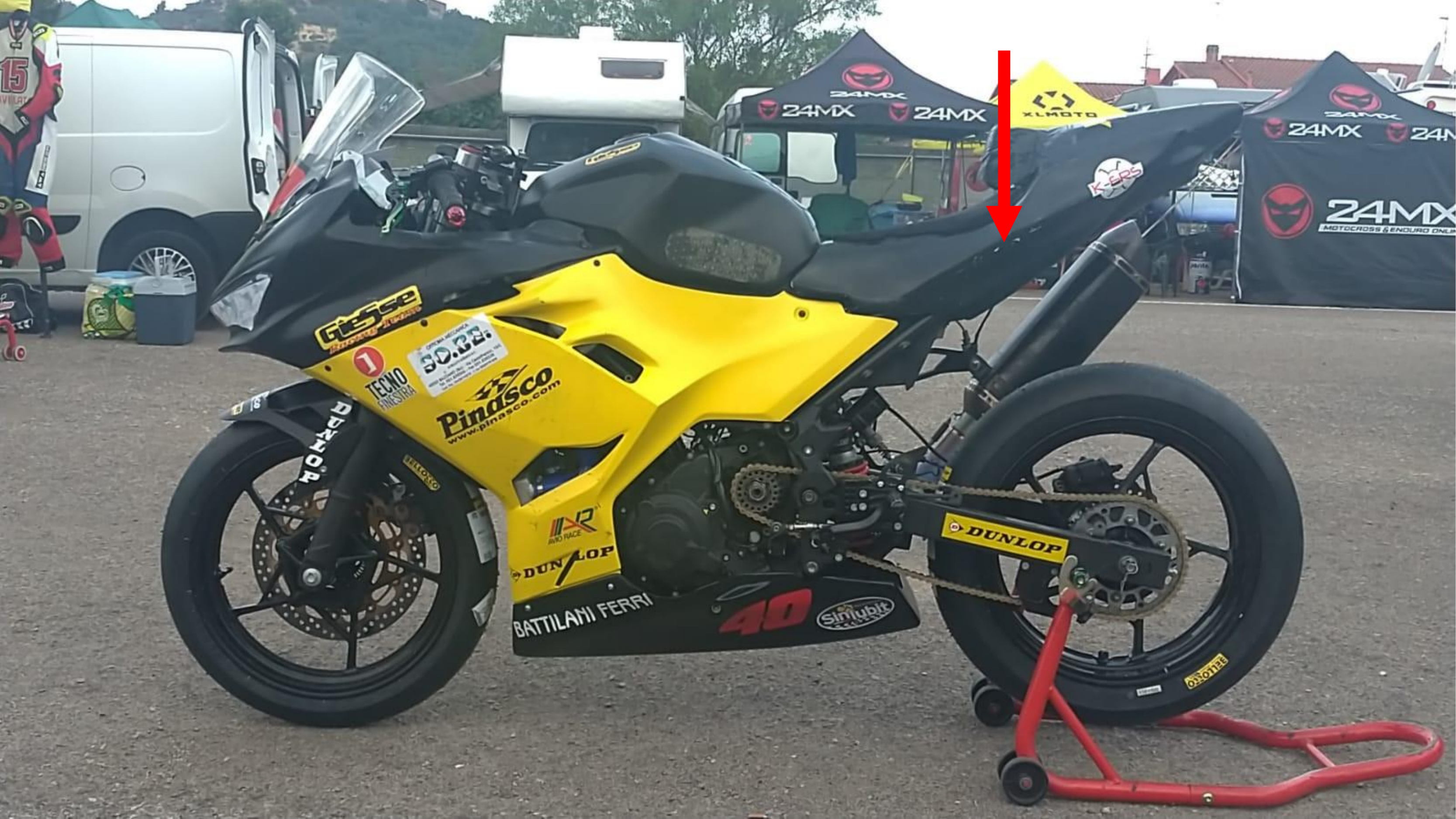}
			\caption{}
			\label{fig:Motorbike}
	\end{subfigure}
\begin{subfigure}[t]{0.42\columnwidth}
	\centering
	\includegraphics[clip, trim= 0cm 1.5cm 0cm 0, width=\textwidth]{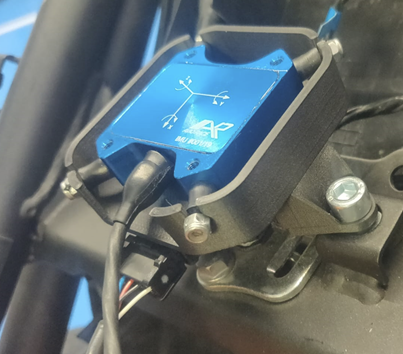}
	\caption{}
\label{fig:IMU}
\end{subfigure}
\caption{(a) The Kawasaki Ninja 400 used for the field test. The yellow arrow points to the location at which the sensor suite is installed. (b) The sensor suite the motorbike was equipped with. It embeds a 9-DOF IMU and a GNSS receiver.}
\label{fig:ExperimentalSetup}
\end{figure}

The field tests have been executed on a Kawasaki Ninja 400 driven by a professional rider in Autodromo Nazionale dell'Umbria "Mario Umberto Borzacchini", see Figure \ref{fig:Umbria}. The motorbike was equipped with a combined IMU+GNSS receiver sensor suite, designed and built by AvioRace, see Figure \ref{fig:ExperimentalSetup}. The IMU consists of a 3-axis accelerometer and a 3-axis gyroscope aligned with the sensor suite symmetry axes, see Figure \ref{fig:IMU}. The sensor suite is installed under the saddle at the location the arrow displayed in Figure \ref{fig:Motorbike} is pointing to.

\begin{figure}
	\centering
	\begin{subfigure}[t]{0.49\columnwidth}
	\centering
	\includegraphics[clip, trim= 1.3cm 0cm 3cm 0cm,width = \textwidth]{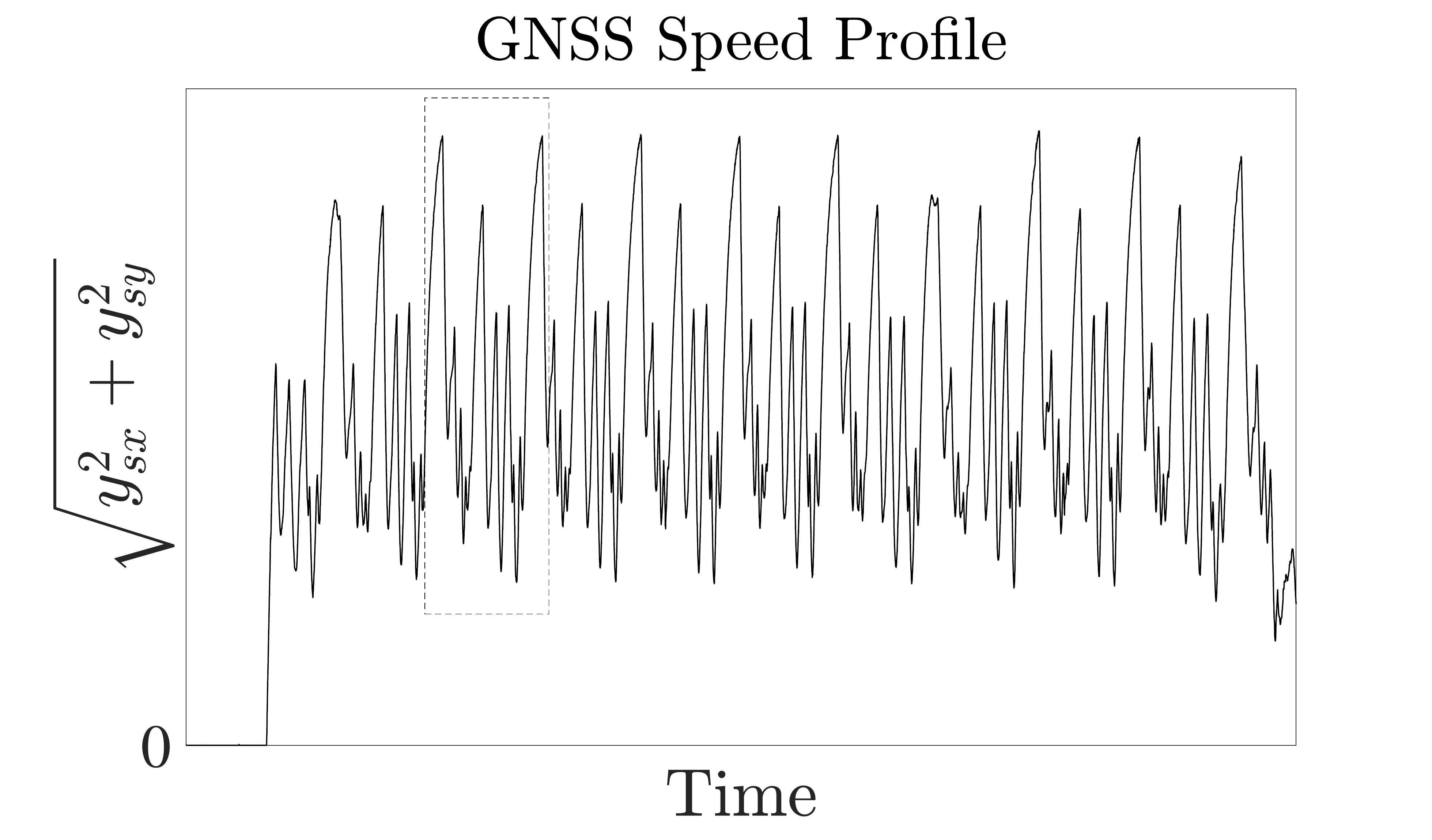}\\
	\includegraphics[clip, trim= 1.34cm 0.5cm 3cm 0cm,width = \textwidth]{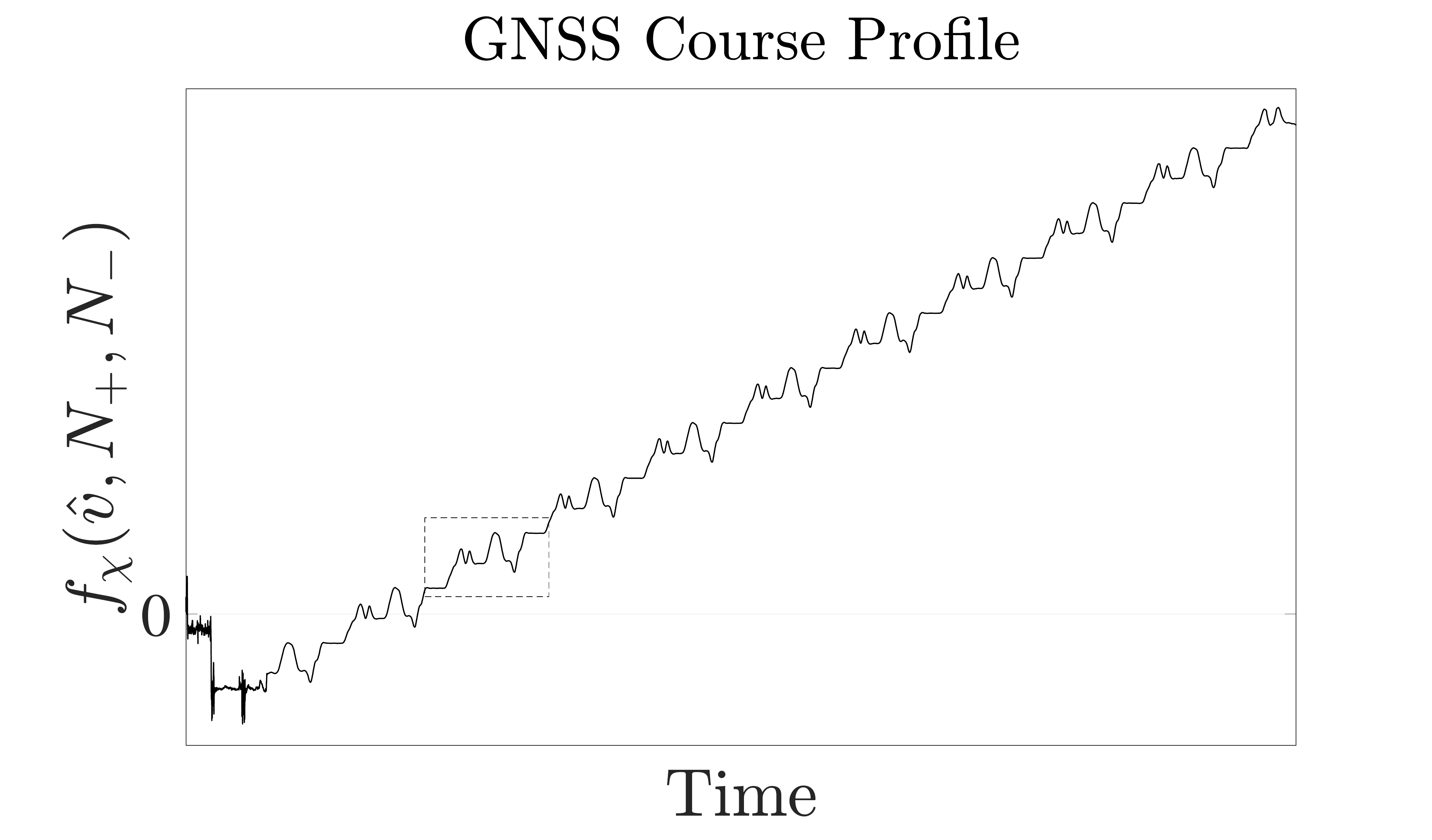}\\
	\caption{}
	\label{fig:SpeedUmbria}
\end{subfigure}
	\begin{subfigure}[t]{0.49\columnwidth}
	\centering
	\includegraphics[clip, trim= 1.3cm 0cm 3cm 0cm, width = \textwidth]{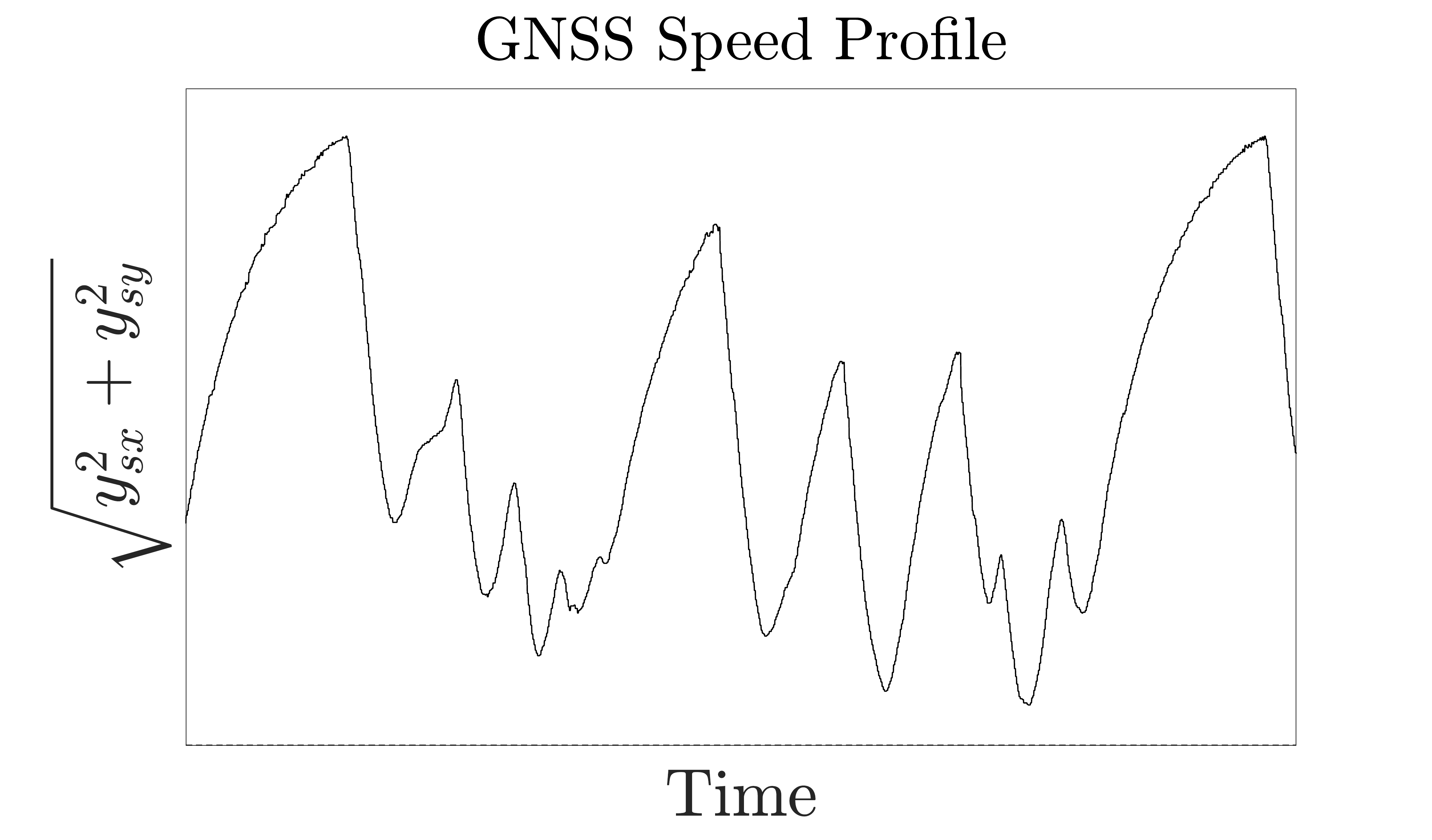}\\
		\includegraphics[clip, trim= 1.3cm 0.5cm 3cm 0cm,width = \textwidth]{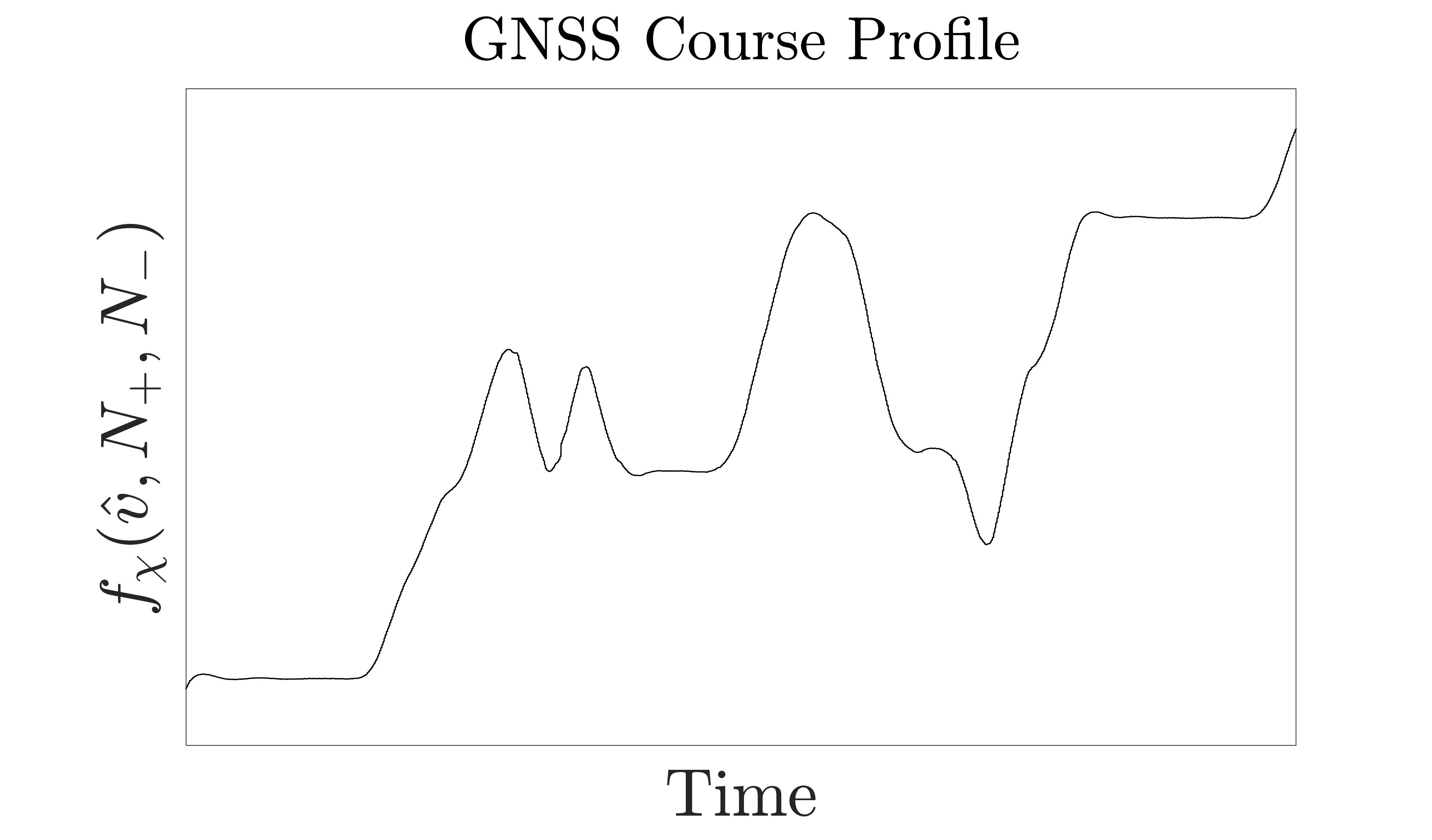}\\
	\caption{}
	\label{fig:SpeedUmbria_zoom}
\end{subfigure}
\caption{(a) GNSS speed profile of the field test. The experiment lasted 9 laps. The 2$^{\text{nd}}$ lap is boxed and magnified in Figure \ref{fig:SpeedUmbria_zoom}. (b) Magnification of the 2$^{\text{nd}}$ lap speed profile. A comparison with Figure \ref{fig:dsdt} highlights the realism of the synthetic data produced through the simulator.}
\label{fig:GNSSUmbria}
\end{figure}

The GNSS data associated to the path illustrated in Figure \ref{fig:Umbria} are reported in Figure \ref{fig:SpeedUmbria} while Figure \ref{fig:SpeedUmbria_zoom} magnifies the second lap (used for the assessment of the realism of the synthetic data produced with the simulator, see Section \ref{sec:Simulation}). As for the GNSS course angle, Figures \ref{fig:SpeedUmbria} and \ref{fig:SpeedUmbria_zoom} show the progressive course made incremental through the lap counter \eqref{eq:LapCounter}. For the presented test, the GNSS vertical speed was not available and the algorithm was evaluated with $\gamma\equiv 0$. 

\begin{figure}
	\centering
	\begin{subfigure}[t]{0.49\columnwidth}
		\centering
		\includegraphics[clip, trim= 1.5cm 0cm 3cm 0cm,width = \textwidth]{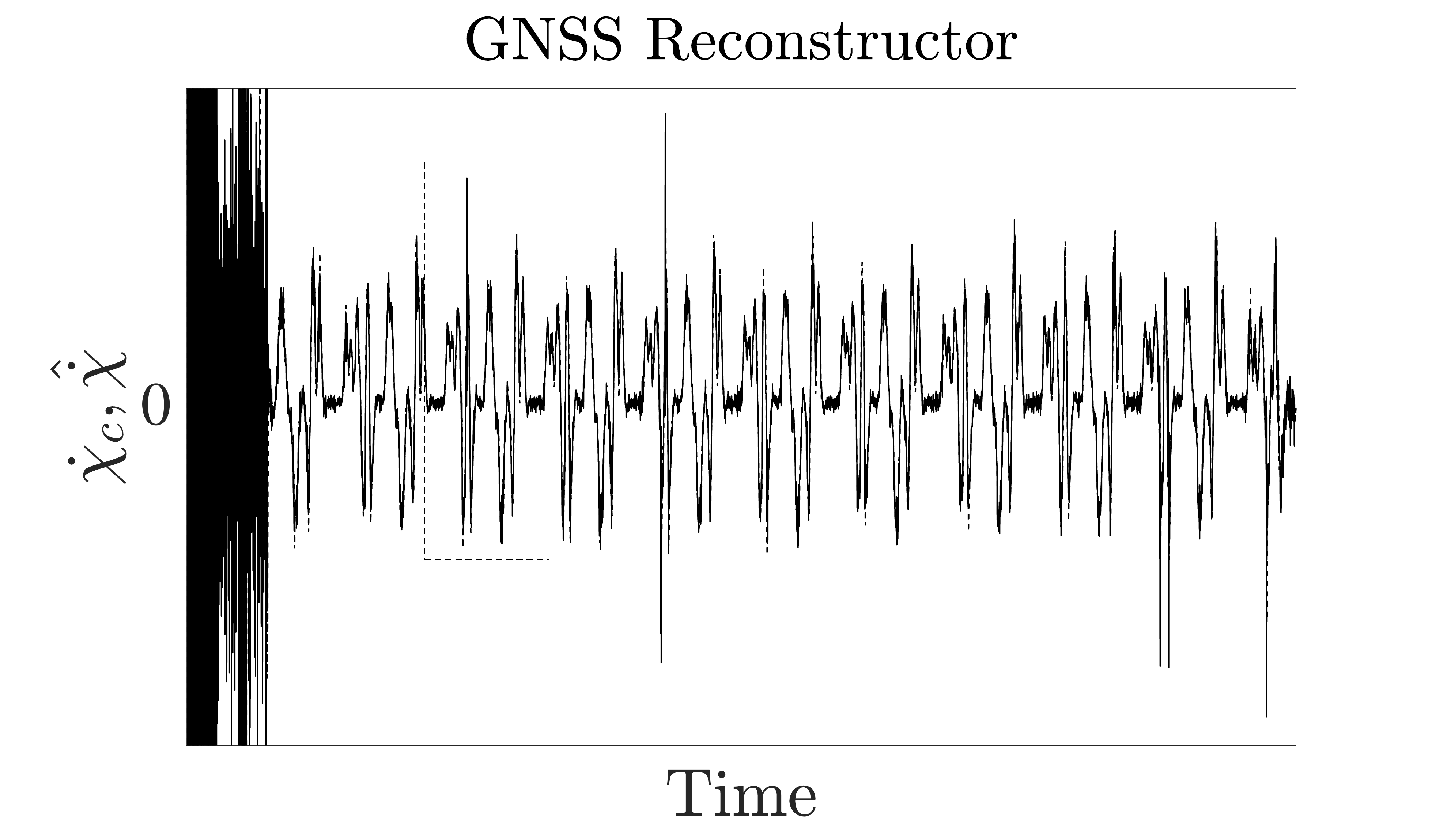}\\
		\caption{}
		\label{fig:hatdotchiUmbria_full}
	\end{subfigure}
	\begin{subfigure}[t]{0.49\columnwidth}
		\centering
		\includegraphics[clip, trim= 1.5cm 0cm 3cm 0cm, width = \textwidth]{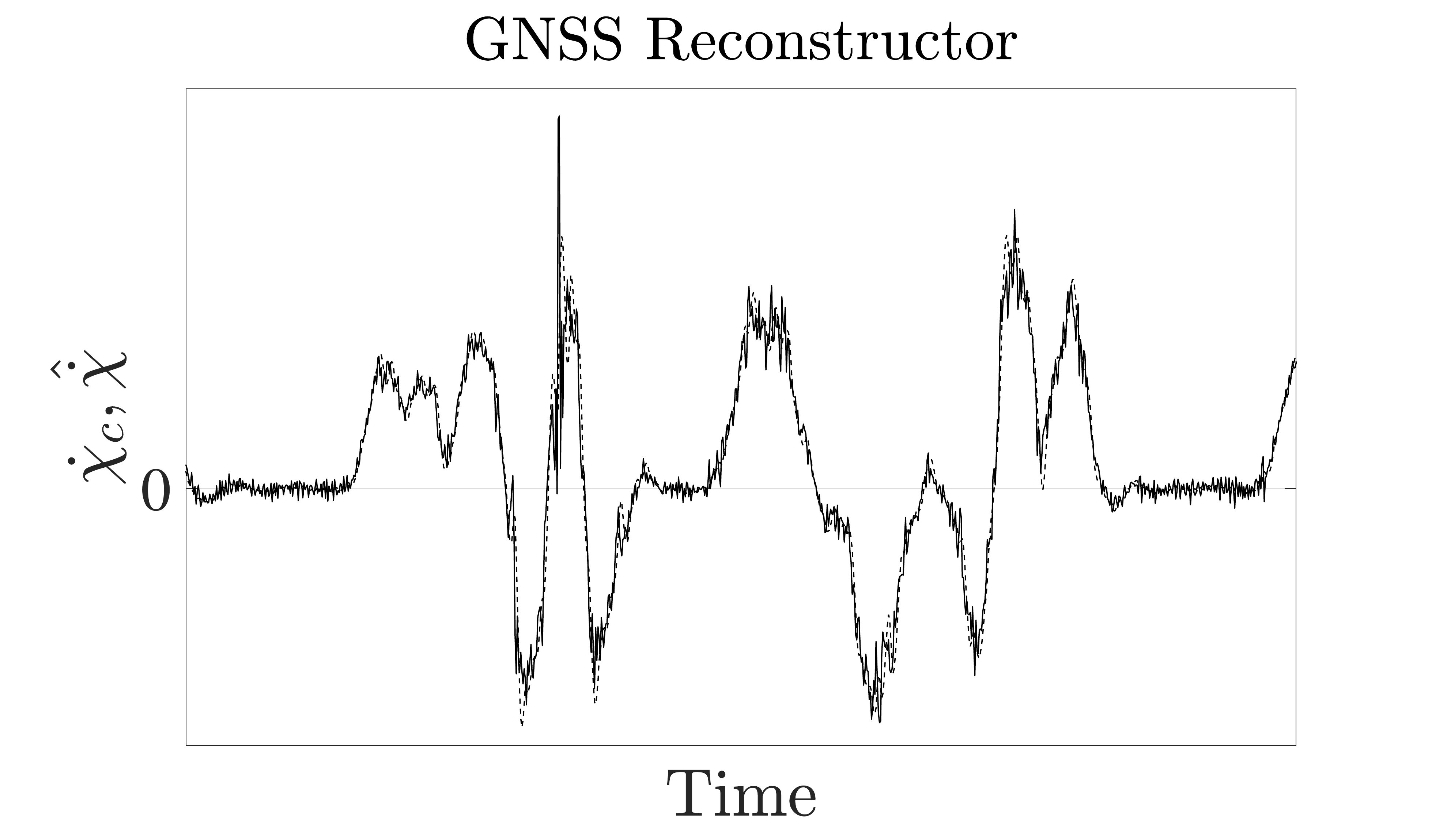}\\
		\caption{}
		\label{fig:hatdotchiUmbria_detail}
	\end{subfigure}
	\caption{(a) GNSS course angle reconstruction. Comparison of the estimations obtained via a batch crude numerical computation (continuous line) and the estimate given by the reconstructor (dashed line). (b) Magnification of the 2$^{\text{nd}}$ lap. The reconstructor  (dashed line) well tracks the reference course angle derivative obtained via crude numerical batch computations (continuous line).}
	\label{fig:hatdotchiUmbria}
\end{figure}

The application of the GNSS reconstructor, described in Section \ref{sec:Pre}, on the course angle of Figure \ref{fig:SpeedUmbria} lead to the estimation of Figure \ref{fig:hatdotchiUmbria}. To appreciate the performance of the reconstructor \eqref{eq:fhatxie} and \eqref{eq:hat_v_vdot}, we compared the estimated $\hat{\dot{\chi}}$ with a batch crude numerical computation. This latter, founded on the central finite difference method of order 8, provides the estimation $\dot{\chi}_c$, see Figure \ref{fig:hatdotchiUmbria}.

\begin{figure}[t]
	\centering
	\begin{subfigure}[t]{0.49\columnwidth}
		\centering
		\includegraphics[clip, trim= 1.5cm 0.5cm 3cm 0cm,width = \textwidth]{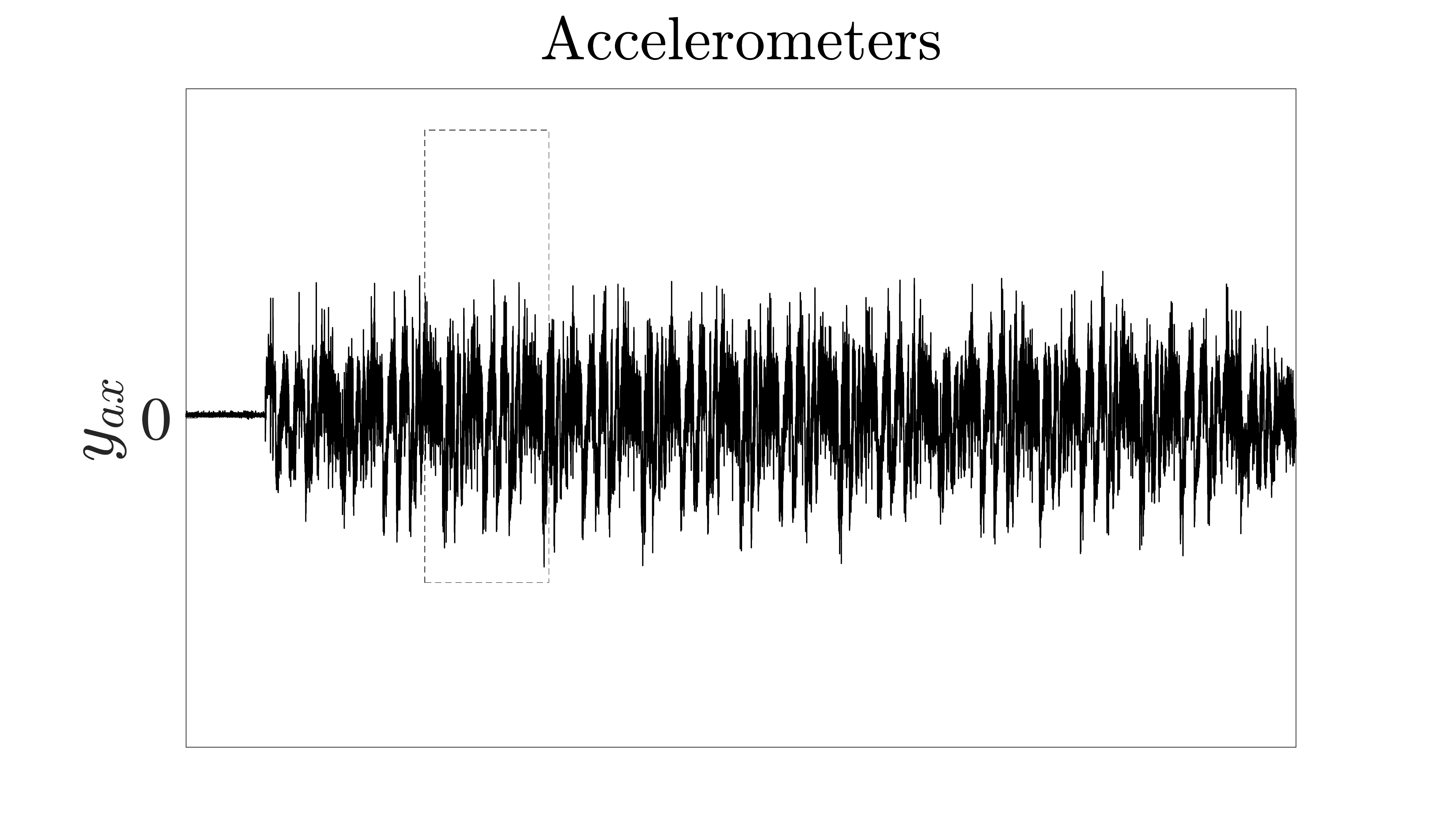}\\
		\includegraphics[clip, trim= 1.5cm 0.5cm 3cm 0cm,width = \textwidth]{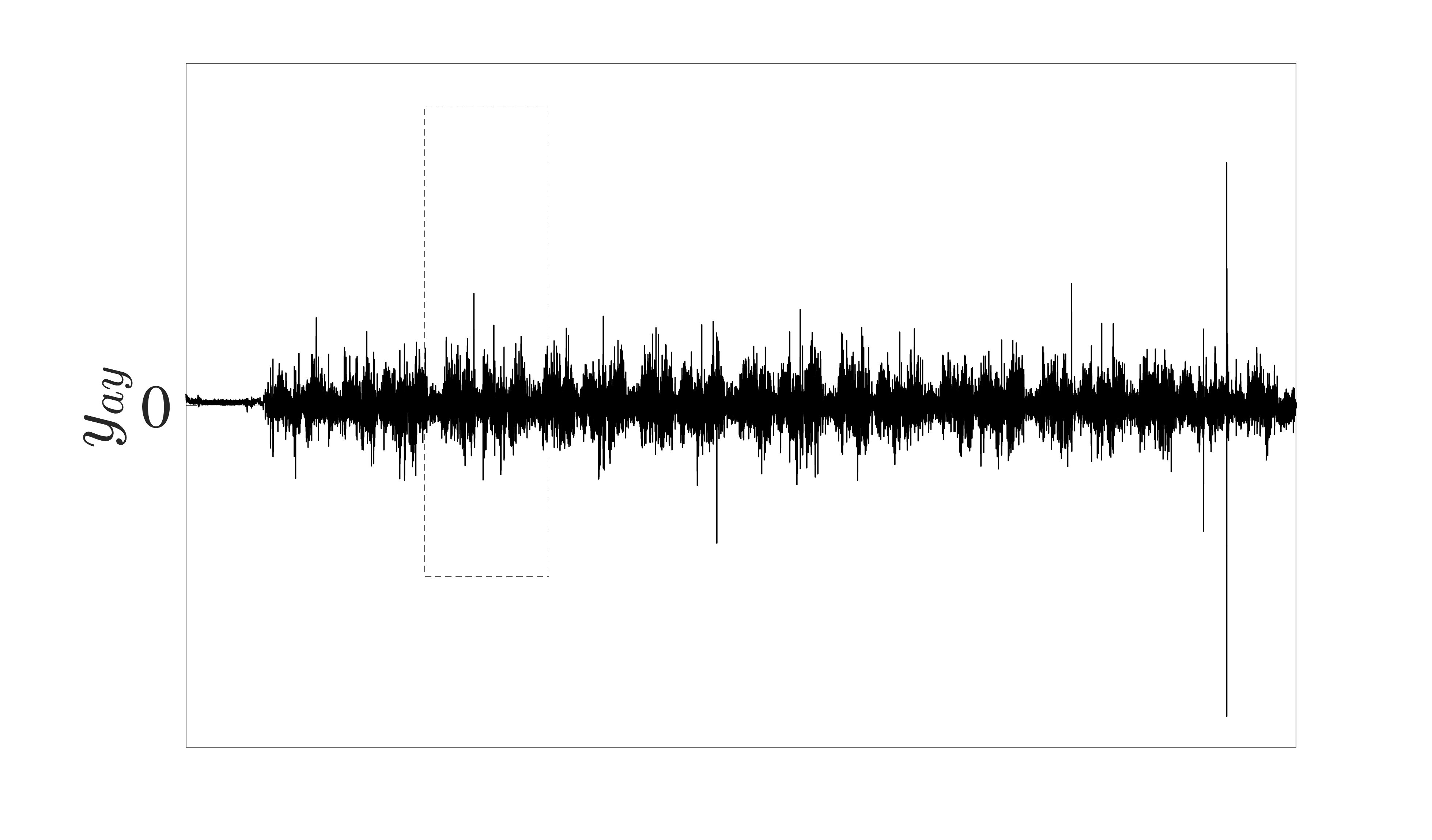}\\
		\includegraphics[clip, trim= 1.5cm 0.5cm 3cm 0cm,width = \textwidth]{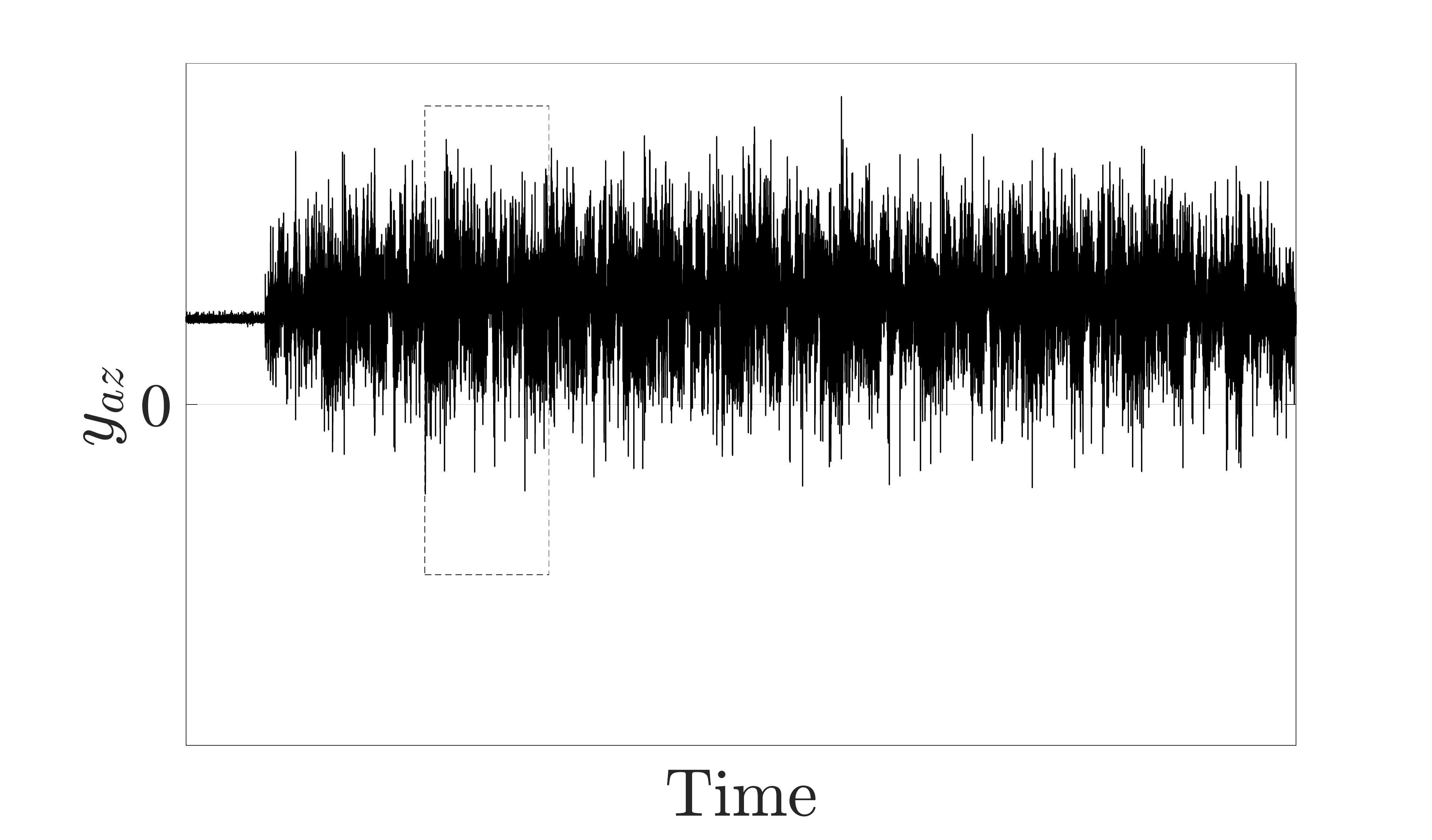}
		\caption{}
		\label{fig:AccUmbria}
	\end{subfigure}
	\begin{subfigure}[t]{0.49\columnwidth}
		\centering
		\includegraphics[clip, trim= 1.5cm 0.5cm 3cm 0cm,width = \textwidth]{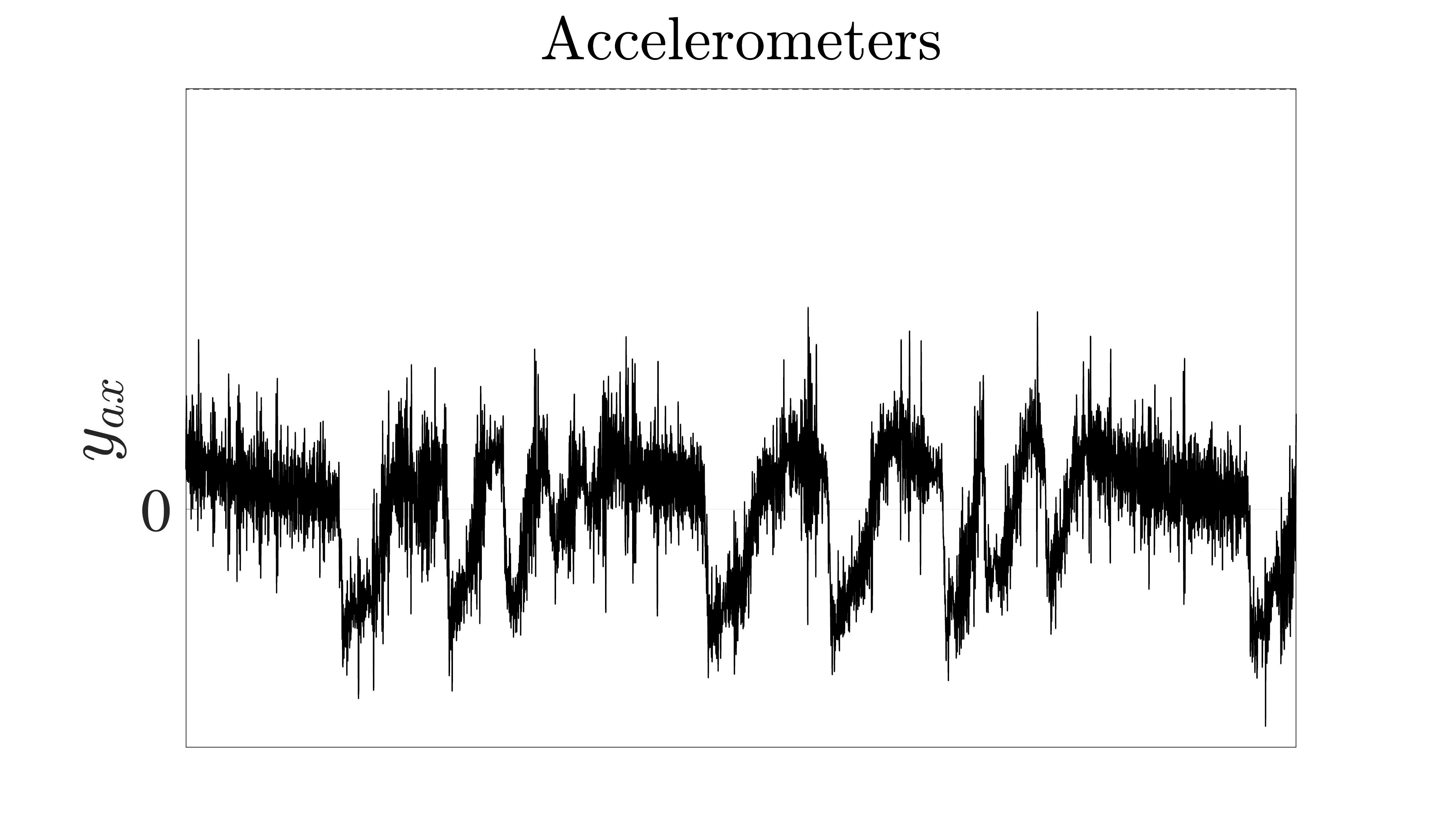}\\
\includegraphics[clip, trim= 1.5cm 0.5cm 3cm 0cm,width = \textwidth]{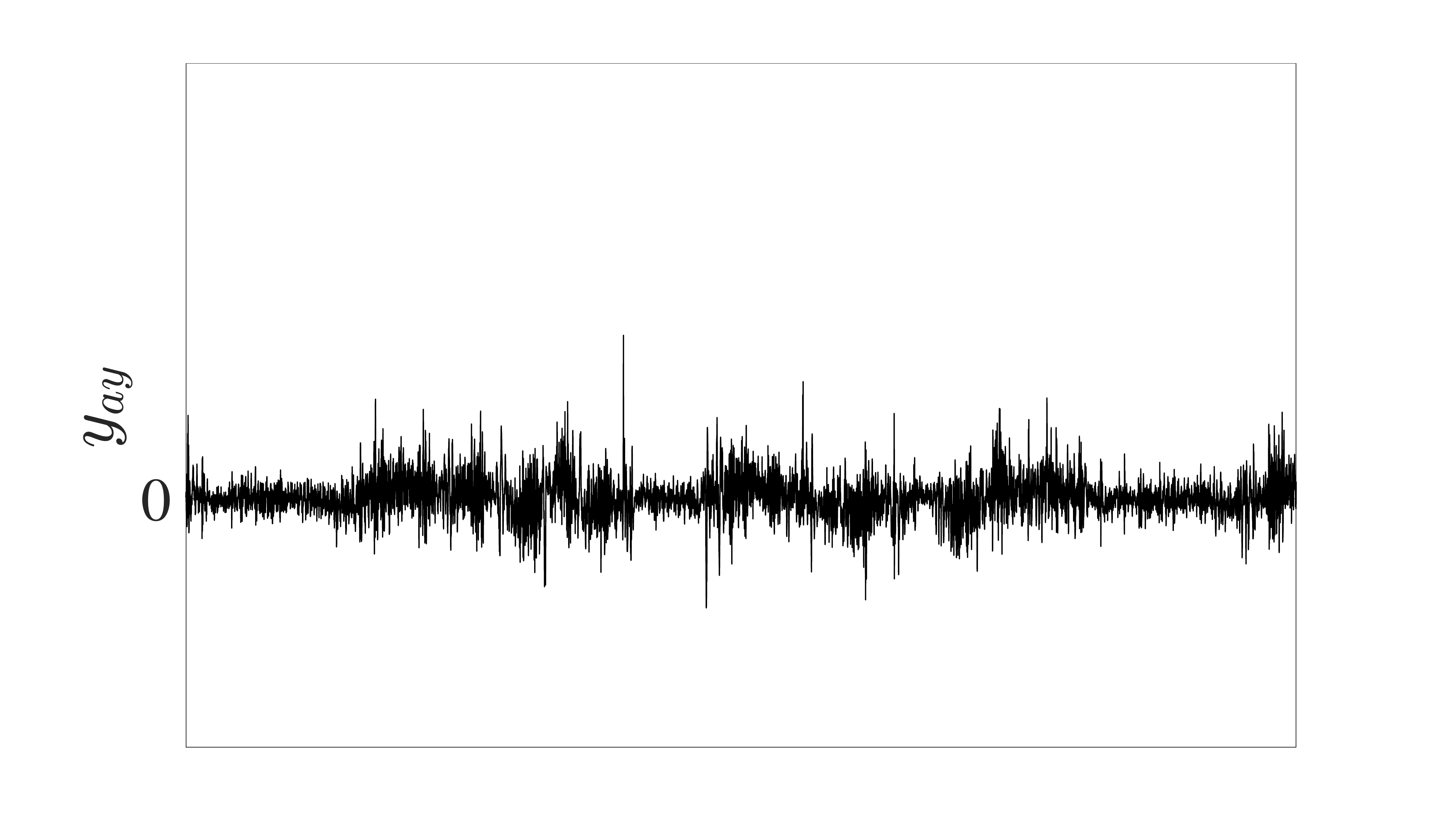}\\
\includegraphics[clip, trim= 1.5cm 0.5cm 3cm 0cm,width = \textwidth]{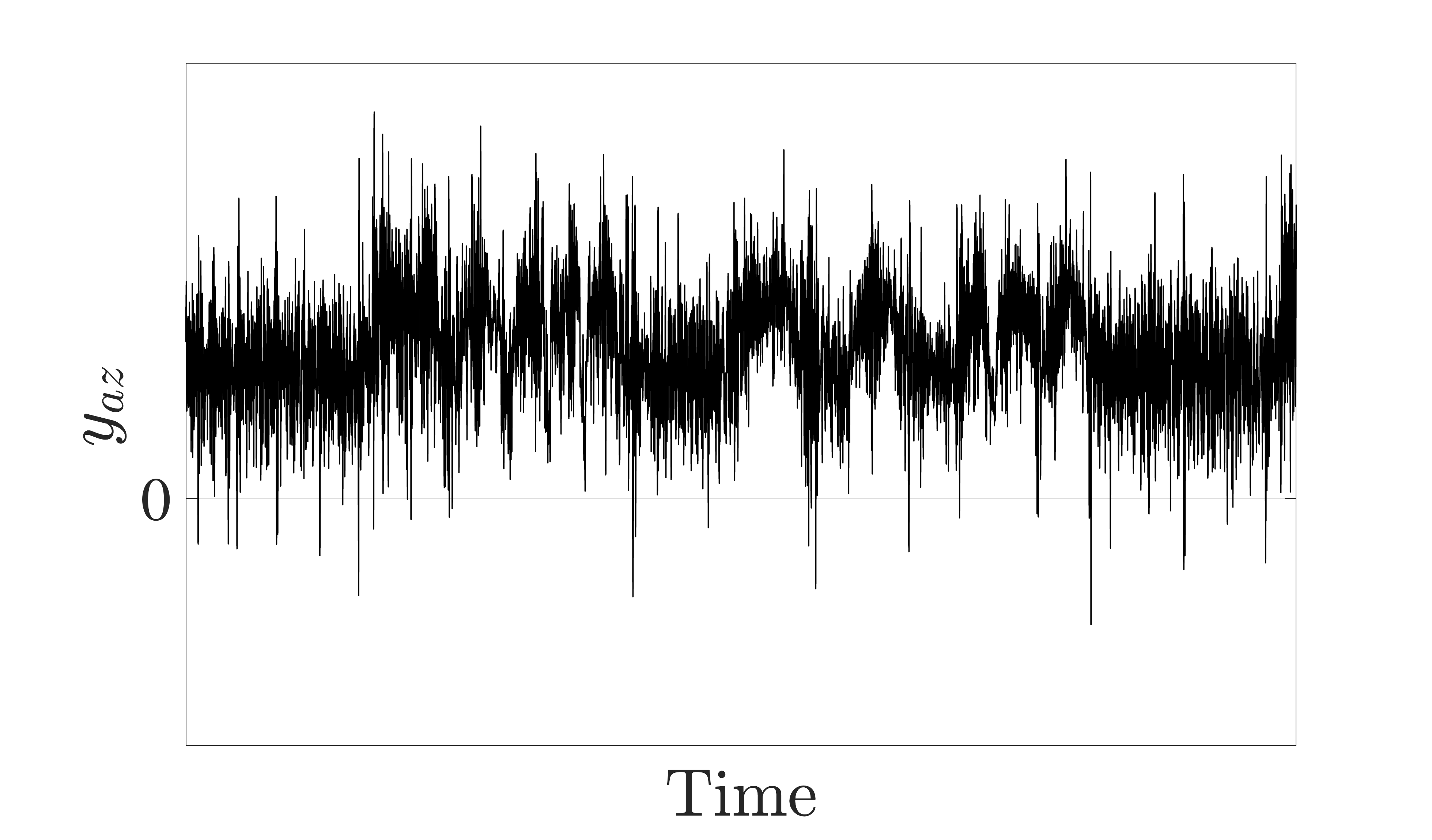}
		\caption{}
		\label{fig:AccUmbriaZoom}
	\end{subfigure}
	\caption{(a) Accelerometer outputs. The 2$^{\text{nd}}$ lap is boxed and magnified in Figure \ref{fig:AccUmbriaZoom}. (b) Magnification of the 2$^{\text{nd}}$ lap. }
	\label{fig:Y1Umbria}
\end{figure}

\begin{figure}[t]
	\centering
	\begin{subfigure}[t]{0.49\columnwidth}
		\centering
		\includegraphics[clip, trim= 1.5cm 0.5cm 3cm 0cm,width = \textwidth]{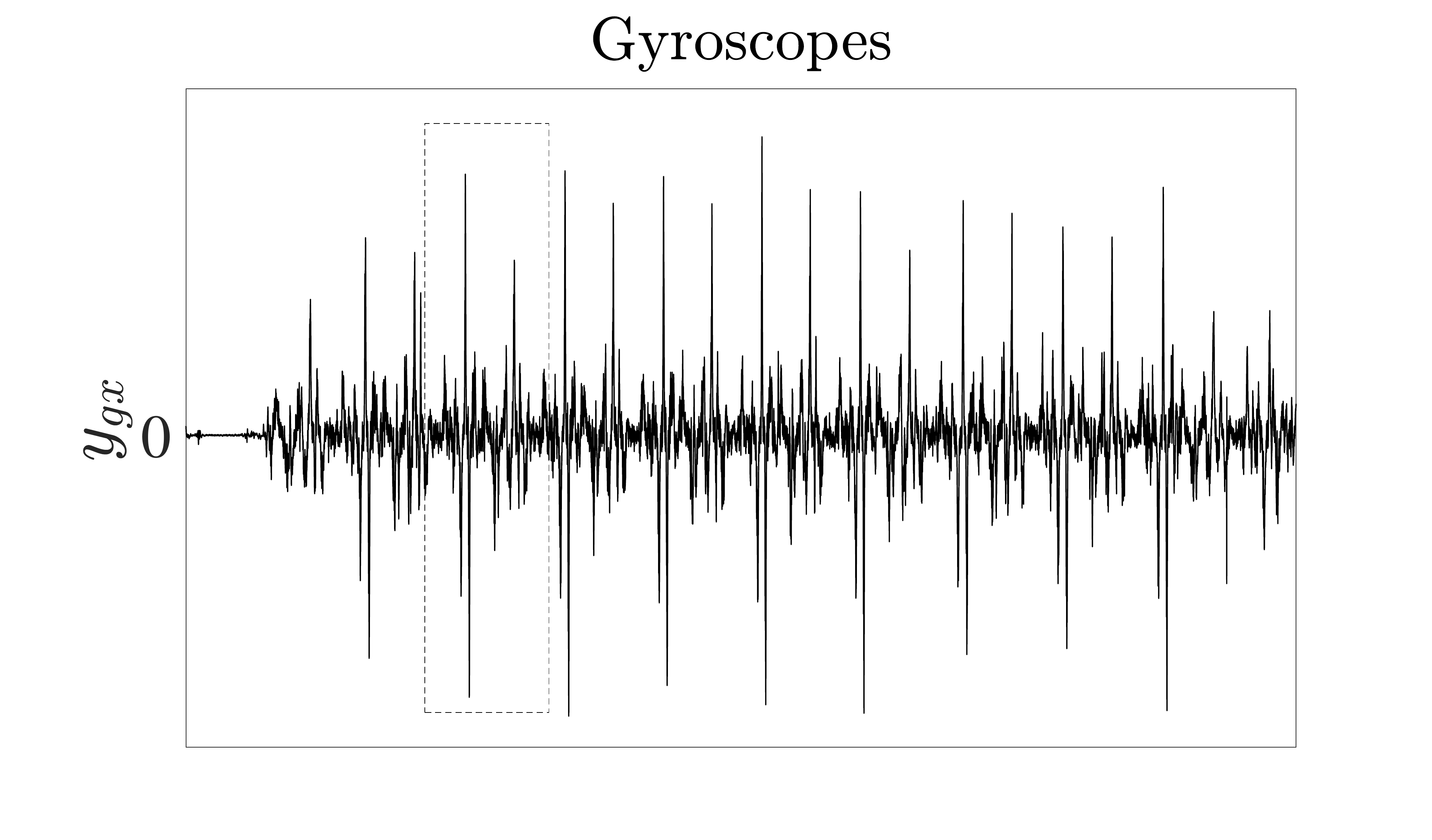}\\
		\includegraphics[clip, trim= 1.5cm 0.5cm 3cm 0cm,width = \textwidth]{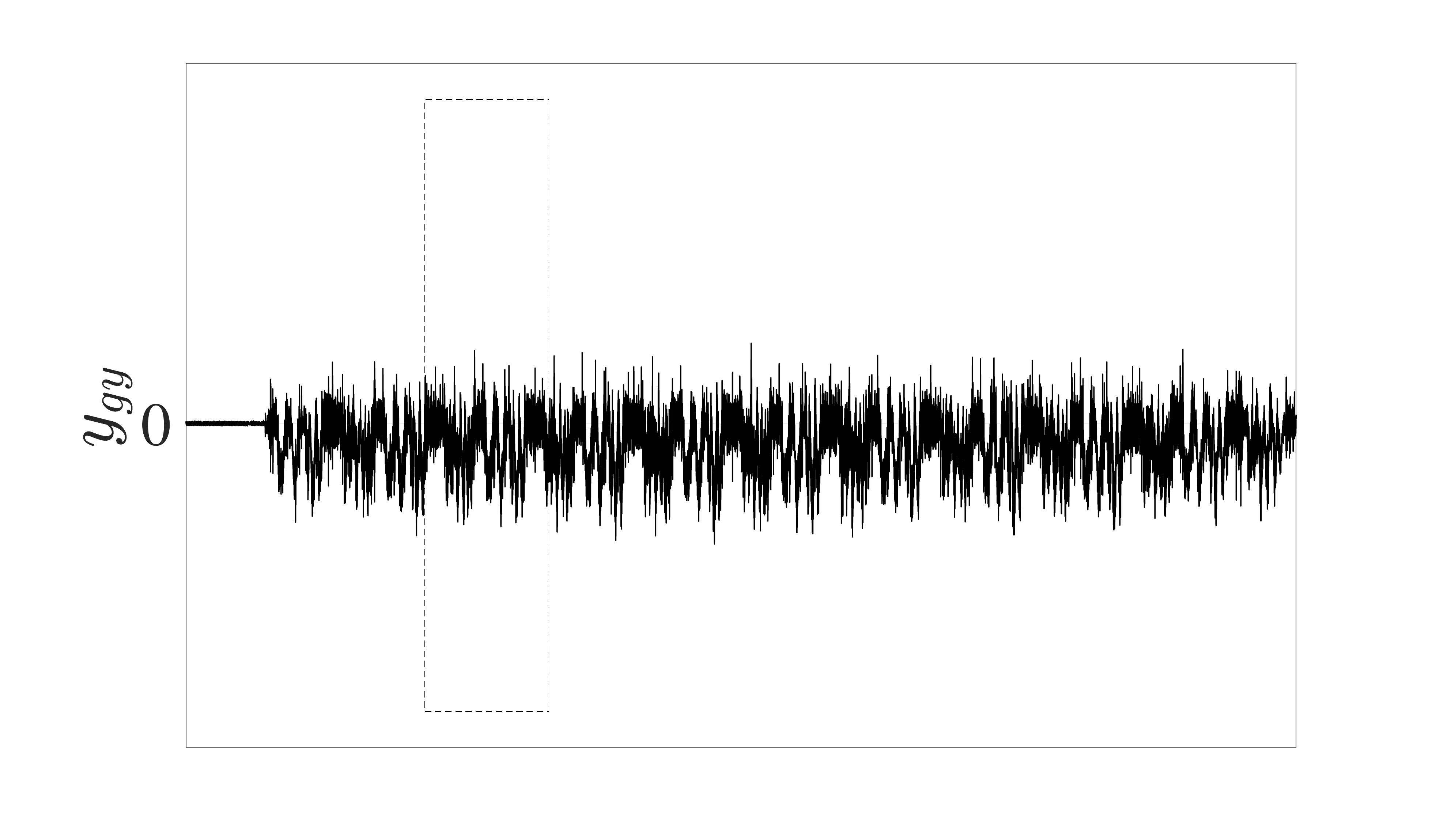}\\
		\includegraphics[clip, trim= 1.5cm 0.5cm 3cm 0cm,width = \textwidth]{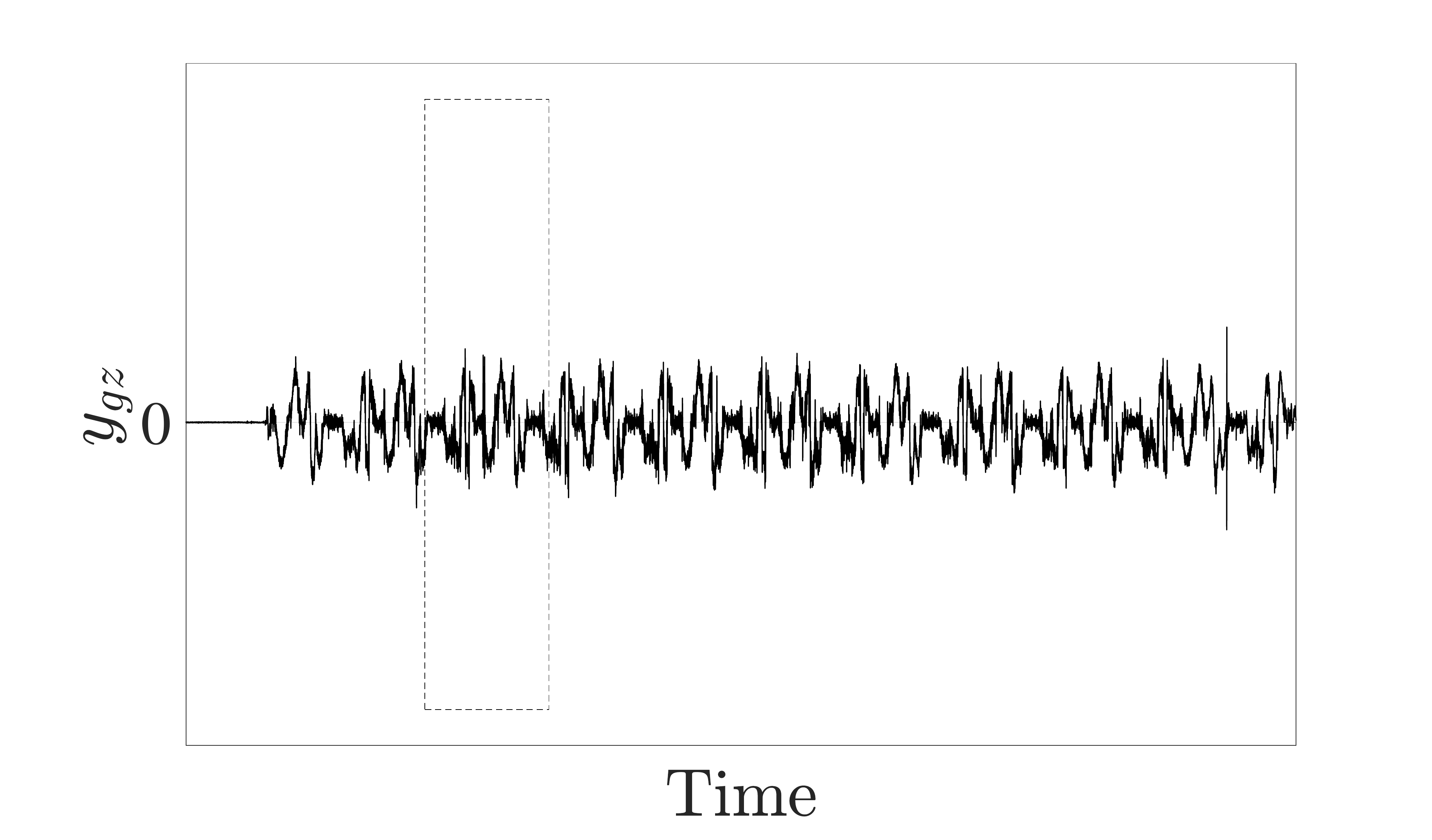}
		\caption{}
		\label{fig:GyroUmbria}
	\end{subfigure}
	\begin{subfigure}[t]{0.49\columnwidth}
		\centering
		\includegraphics[clip, trim= 1.5cm 0.5cm 3cm 0cm,width = \textwidth]{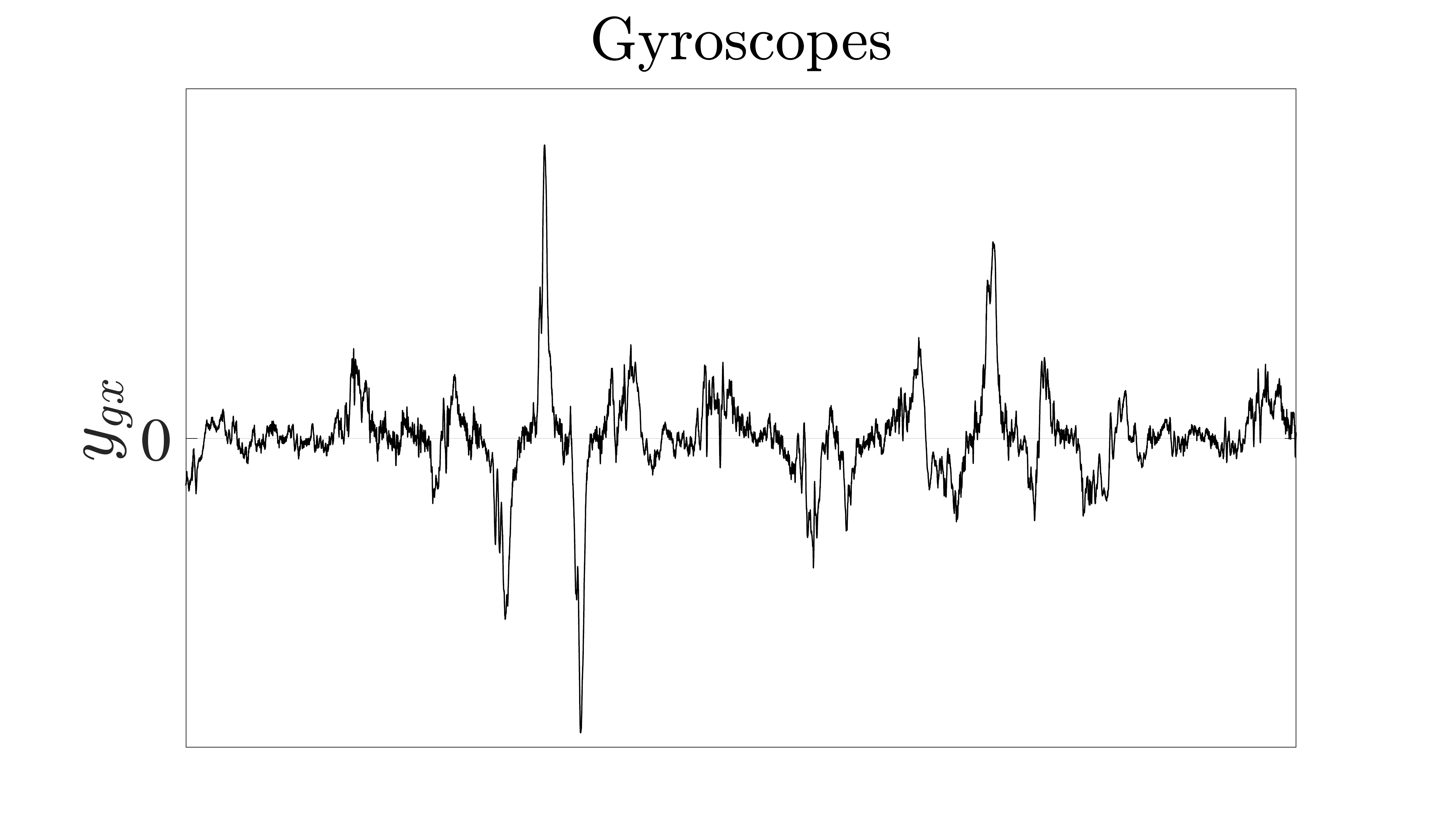}\\
		\includegraphics[clip, trim= 1.5cm 0.5cm 3cm 0cm,width = \textwidth]{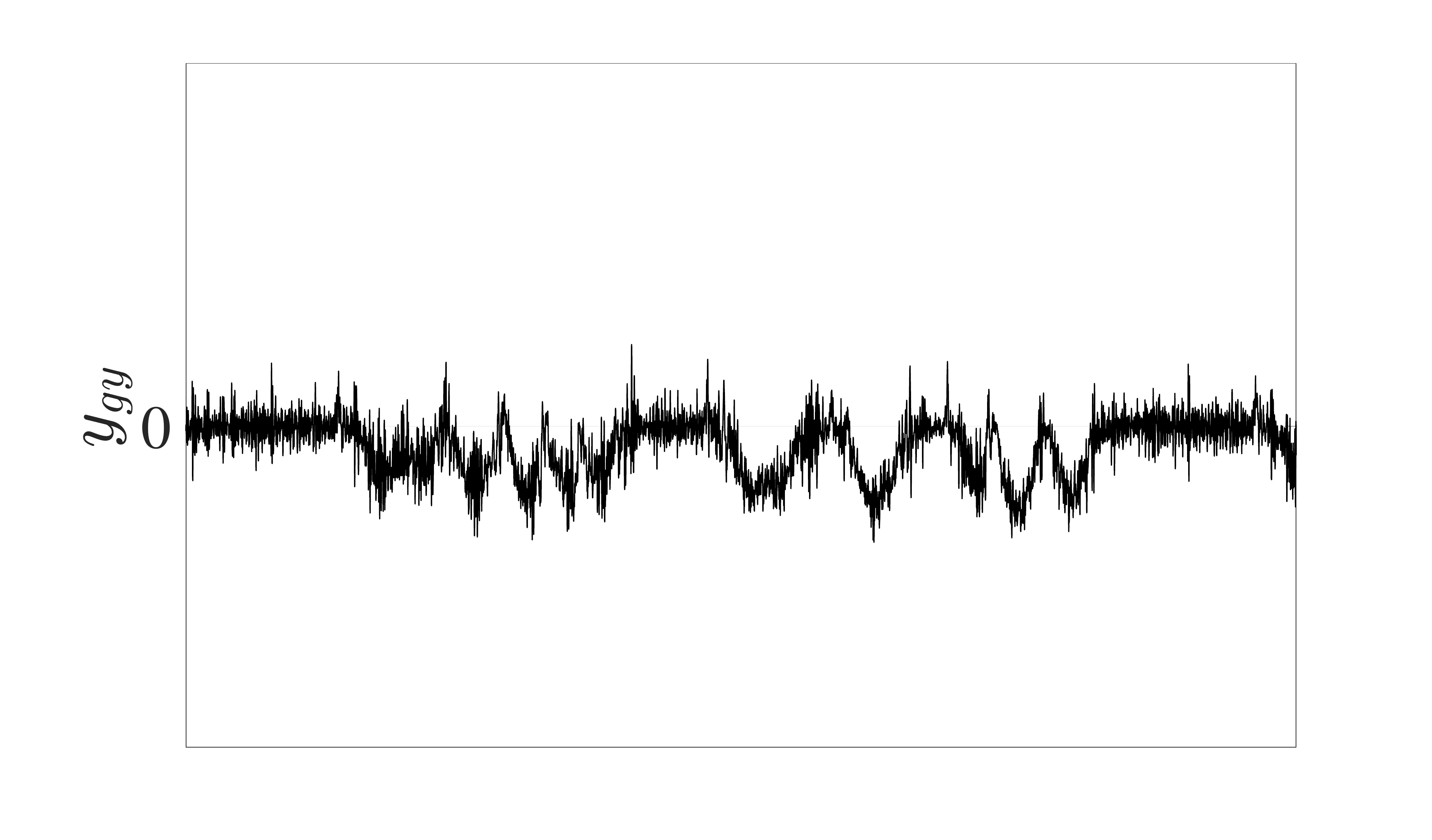}\\
		\includegraphics[clip, trim= 1.5cm 0.5cm 3cm 0cm,width = \textwidth]{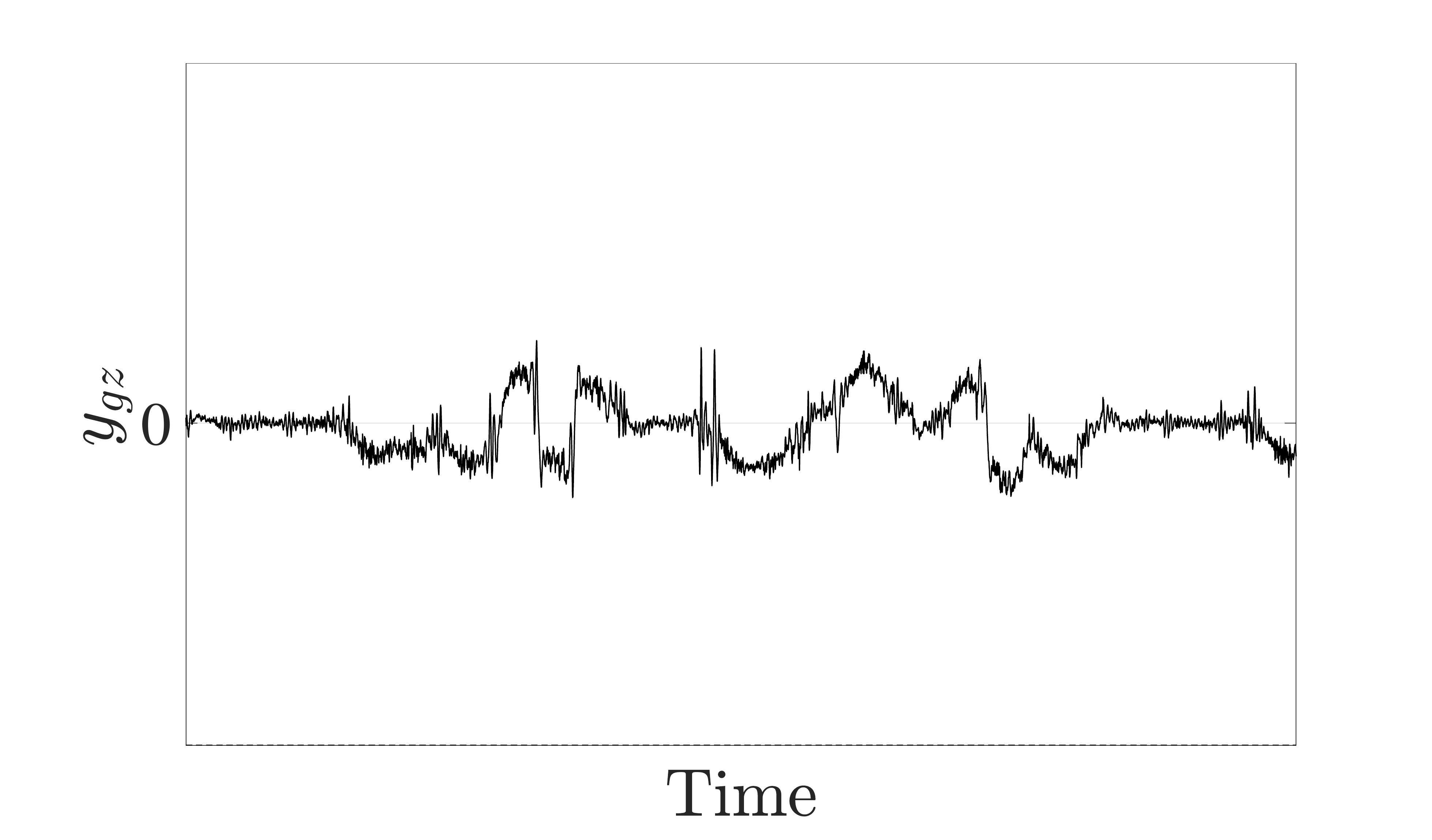}
		\caption{}
		\label{fig:GuroUmbriaZoom}
	\end{subfigure}
	\caption{(a) Gyroscope outputs. The 2$^{\text{nd}}$ lap is boxed and magnified in Figure \ref{fig:AccUmbriaZoom}. (b) Magnification of the 2$^{\text{nd}}$ lap. }
	\label{fig:Y2Umbria}
\end{figure}

The data produced by accelerometers and gyroscopes are shown in Figures \ref{fig:Y1Umbria} and \ref{fig:Y2Umbria}, both calibrated for compensating the installation misalignment. 

\begin{figure}[t]
	\centering
	\begin{subfigure}[t]{0.49\columnwidth}
		\centering
		\includegraphics[clip, trim= 0cm 0.5cm 3cm 0cm,width = \textwidth]{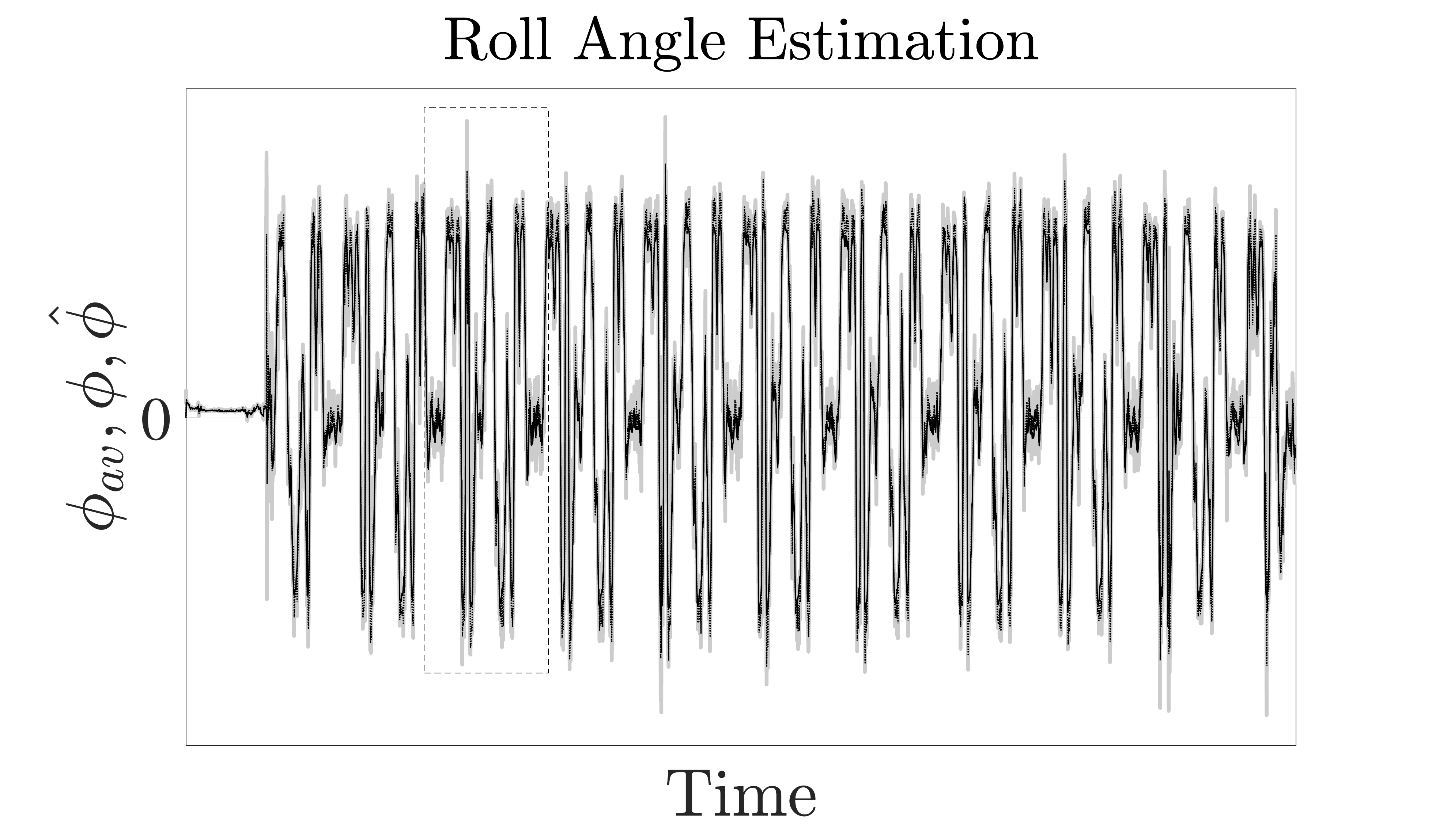}\\
		\caption{}
	\end{subfigure}
	\begin{subfigure}[t]{0.49\columnwidth}
		\centering
		\includegraphics[clip, trim= 0cm 0.5cm 3cm 0cm,width = \textwidth]{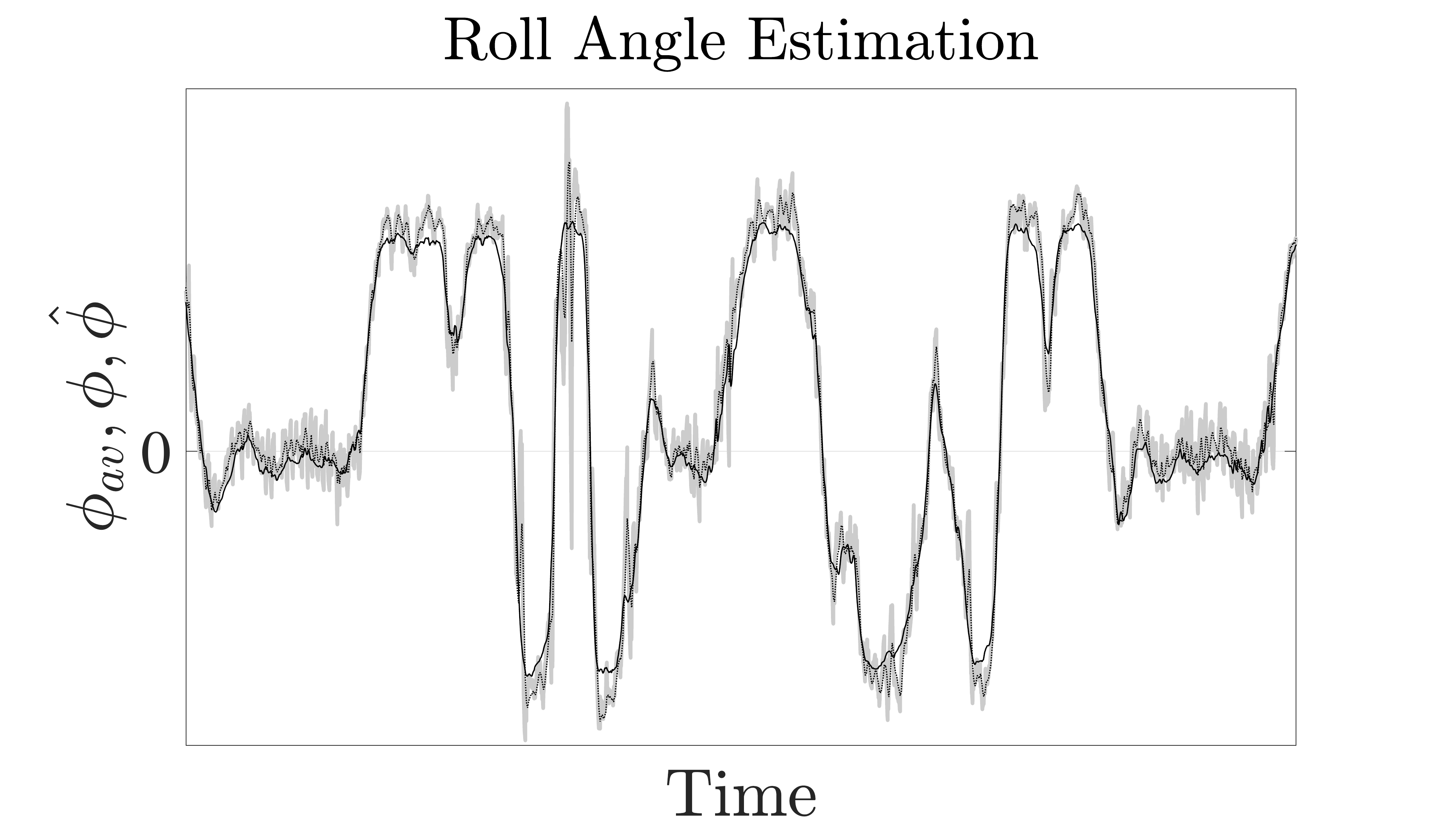}\\
		\caption{}
	\end{subfigure}
	\begin{subfigure}[t]{0.49\columnwidth}
	\centering
	\includegraphics[clip, trim= 0cm 0.5cm 3cm 0cm,width = \textwidth]{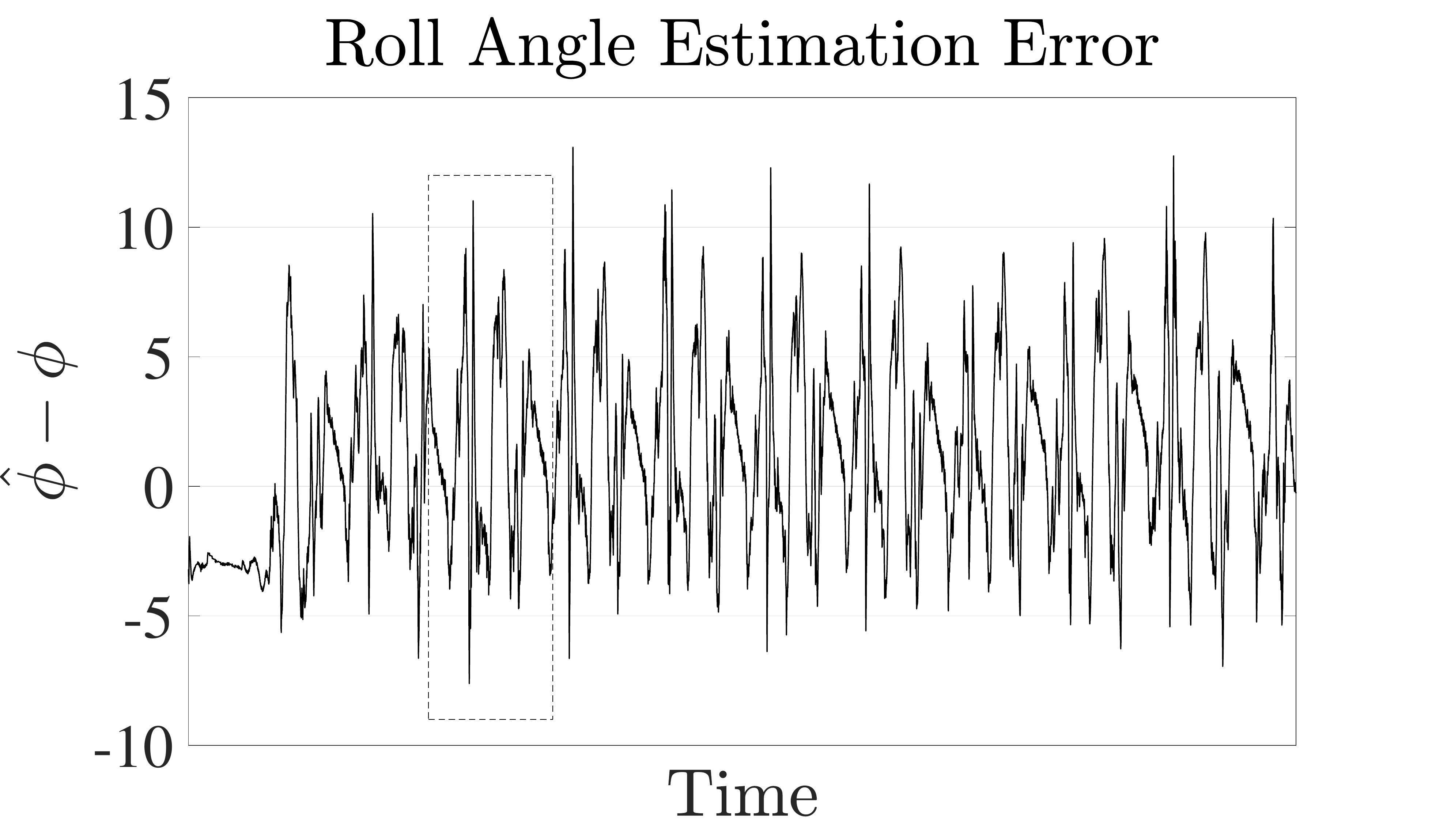}
	\caption{}
\end{subfigure}
	\begin{subfigure}[t]{0.49\columnwidth}
	\centering
	\includegraphics[clip, trim= 0cm 0.5cm 3cm 0cm,width = \textwidth]{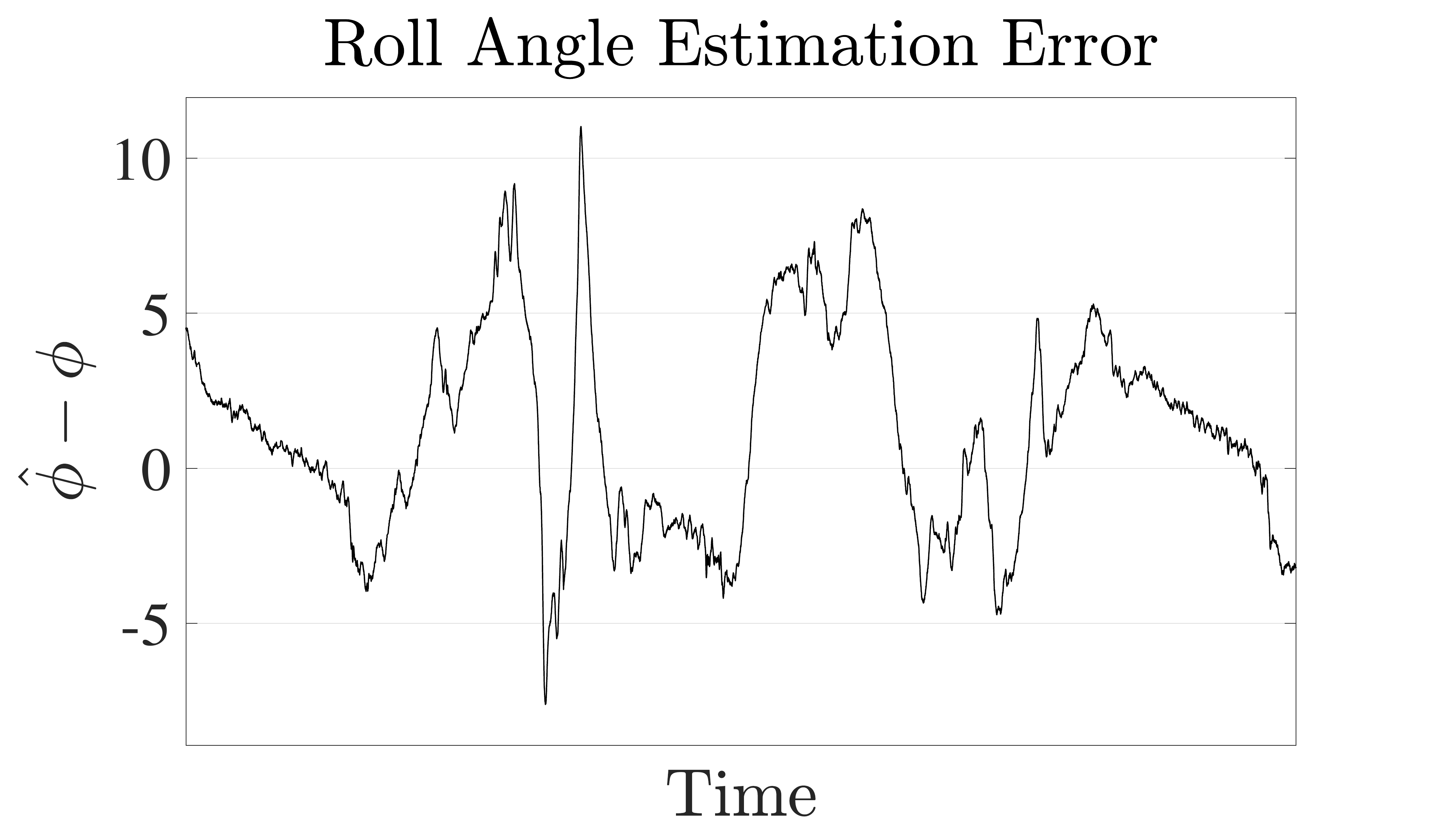}
	\caption{}
\end{subfigure}
	\caption{(a) Roll angle estimation. (b) Magnification of the 2$^{\text{nd}}$ lap. The estimation $\phi_{\text{av}}$ generated accordingly to the assumption of coordinated manoeuvre is reported in grey. The dotted line represents the estimation $\hat{\phi}$ provided by the observer. The continuous line denotes the reference angle $\phi$. (c) Roll angle estimation error. (d) Magnification of the 2$^{\text{nd}}$ lap. The estimation is satisfactory even in those few moments (isolated picks) in which the assumption of coordinated manoeuvre is violated by a drifting condition.}
	\label{fig:phi}
\end{figure}

These data, together with $\hat{\xi}_e$ provided by \eqref{eq:fhatxie}, \eqref{eq:hat_v_vdot}, and \eqref{eq:fchi_new2}, are exploited through \eqref{eq:phi_av} to provide $\phi_{\text{av}}$. This latter is then compared with the estimation of $\phi$ elaborated by a proprietary algorithm, also using a 3-axis magnetometer, and considered as a reliable reference, see Figure \ref{fig:phi}. It is worth noting that, the assumption of coordinated manoeuvre, on which the observer proposed in this paper is based, accurately models the actual dynamics of system motorbike+biker. Indeed, $\phi_\text{av}$ is highly coherent with the reference $\phi$ even in those circumstances in which the rider made the bike skidding. For the test described in this section, we computed the mean and the standard deviation of the roll estimation error as key performance indices, with ${\rm E}[\hat{\phi}(t)-\phi(t)] \approx 1.2$ deg and $\sqrt{E[((\hat{\phi}(t)-\phi(t))-{\rm E}[\hat{\phi}(t)-\phi(t)])^2]} \approx 4.6$ deg. 

%

\subsection{Simulations}
\label{sec:Simulation}

\begin{figure}
	\centering
	\includegraphics[clip, trim= 3cm 0 3cm 0, width=\columnwidth]{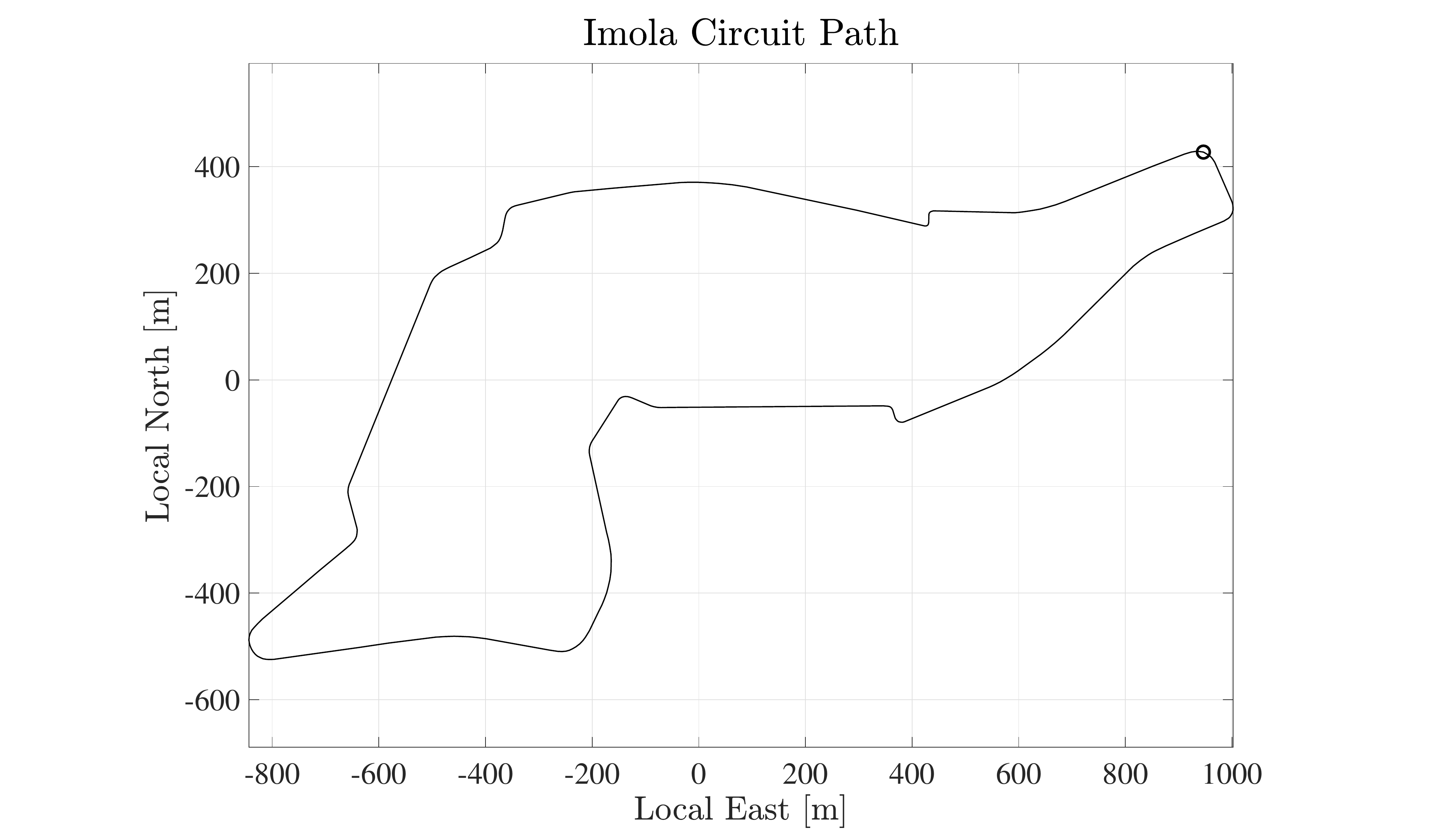}
	\caption{Imola circuit, for motorbikes, represented in the NED reference frame. The circle denotes the starting point of the curvilinear abscissa (counter-clockwise direction).}
	\label{fig:Imola}
\end{figure}

The aim of this section is that of comparing, in clear, the performance of our algorithm to those of the most relevant algorithms found in the literature. In particular, since one of the fundamental elements of our work is represented by the introduction of $\phi_{\text{av}}$ through the definition of the coordinated manoeuvres (see Section \ref{sec:Coordinated}), we evaluated the algorithms proposed in \cite{Boniolo2010Motorcycle,Schlipsing2012Roll,Maceira2021Roll} because, under the assumption of a flat-coordinated turn, they provide possible alternatives to $\phi_{\text{av}}$. More in details, let $v_x^B \in \mathbb{R}$ be projection of the inertial speed on the motorcycle x-axis. Then, we investigated
\begin{itemize}
	\item \phantom{.}[\cite{Boniolo2010Motorcycle}, Eq. (9)] which adopts the z-axis gyroscope measurement to compute
		\[
		\phi_1 = -\tan^{-1}\left(\dfrac{v_x^B y_{2z}}{\mathtt{g}}\right)
		\] 
	\item \phantom{.}[\cite{Boniolo2010Motorcycle}, Eq. (16)] that exploits the y-axis gyroscope to elaborate
		\[
		\begin{aligned}
		\phi_2 =&\, -\mathtt{sign}(y_{2z})\cos^{-1}\left(\sqrt{1+\Phi ^2}-\Phi \right)\\
		\Phi =&\, \dfrac{v_x^B|y_{2y}|}{2\mathtt{g}}
	\end{aligned}
		\]
	\item \phantom{.}[\cite{Boniolo2010Motorcycle}, Eq.s (21), (22)] using both the y- and z-axis gyroscopes to determine
		\[
		\phi_3 = -\tan^{-1}\left(\dfrac{v_x^B}{\mathtt{g}}\mathtt{sign}(y_{2z})\sqrt{y_{2y}^2+y_{2z}^2}\right),
		\]
	\item \phantom{.}[\cite{Schlipsing2012Roll}, Eq. (7)] which, in implicit form, represents a heuristic improvement of $\phi_1$ 
		\[
		\phi_4\in \mathbb{R}\,:\,  \tan(0.9\phi_4) \cos(\phi_4) = -\left(\dfrac{v_x^B y_{2z}}{\mathtt{g}}\right).
		\]
	\item \phantom{.}[\cite{Maceira2021Roll}, Eq. (54), (56)-(58)] that introduces an heuristic weight function to mix $\phi_1$ and a proxy of $\phi_3$ as
		\[
		\begin{aligned}
			\phi_5 =&\, W \phi_1 -\left(1-W\right)\mathtt{sign}(y_{2z})\sin^{-1}\left(\dfrac{y_{2y}}{\sqrt{y_{2y}^2+y_{2z}^2}}\right)\\
			W =&\, \exp(-25\phi_1^2).
		\end{aligned}
		\]
\end{itemize}

To create a synthetic but realistic dataset we modelled a track lap. More in detail, we added to the Imola circuit path (see Figure \ref{fig:Imola}) a time law whose generation, made exploiting \cite{Hauser2006Motorcycle}, takes into account  lateral and longitudinal maximum tire forces, wheelie conditions, engine power and efficiency, circuit slope and aerodynamic drag. 

\begin{figure}
	\centering
	\includegraphics[clip, trim= 1cm 0 1cm 0, width=\columnwidth]{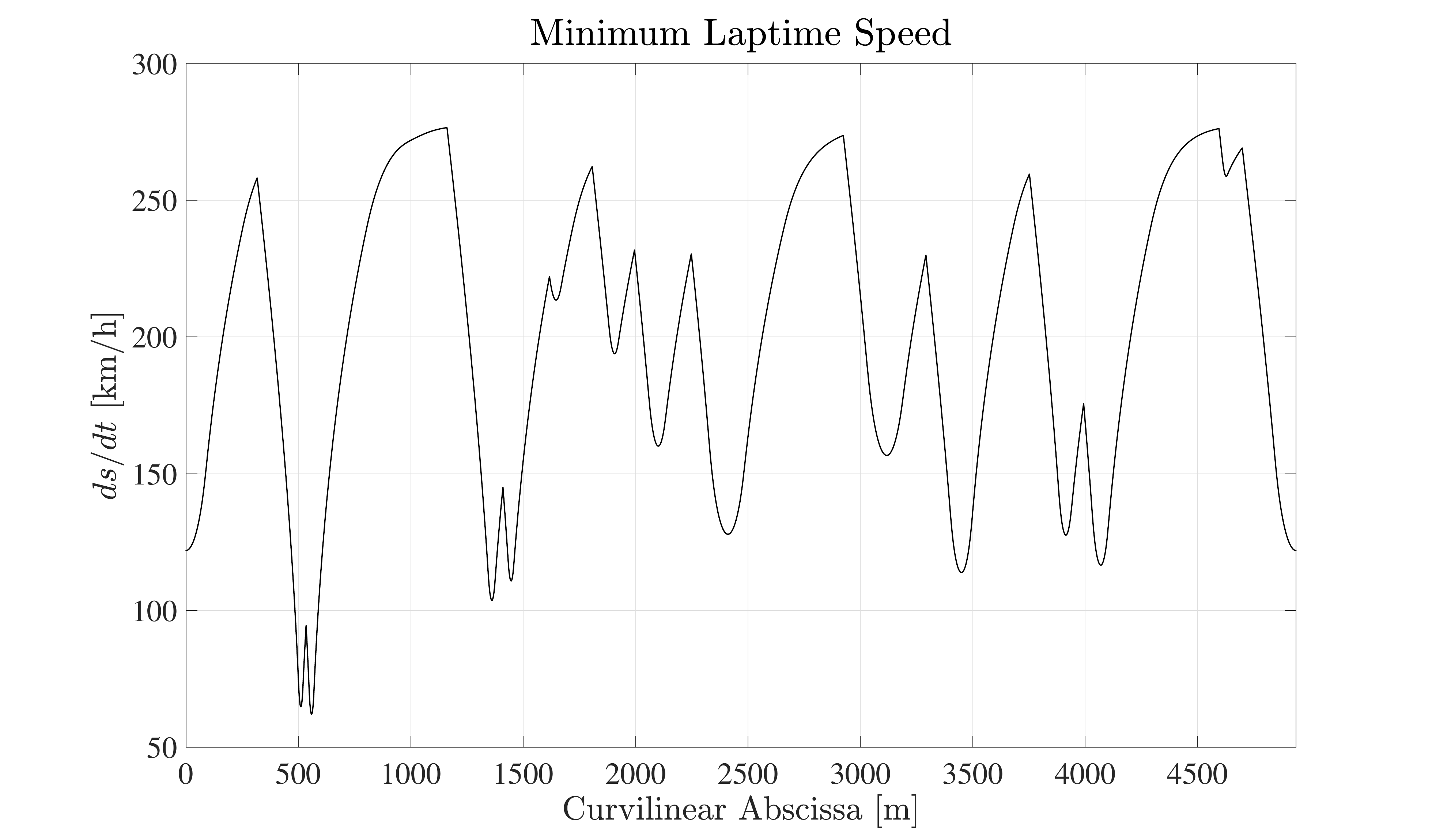}
	\caption{The minimum laptime speed is generated accordingly to the road slope, the aerodynamic drag, the engine performance, and the tire maximum cohesion coefficient.}
	\label{fig:dsdt}
\end{figure}

Let $L>0$ be the circuit length, then this procedure leads to the definition of the curvilinear speed and acceleration, namely $ds/dt(\cdot),\,d^2s/dt^2(\cdot)\,:\,[0\,\,L] \to \mathbb{R}^3$ (see Figure \ref{fig:dsdt}), and the heading and the slope $\chi(\cdot),\gamma(\cdot)\,:\,[0\,\,L]\to \mathbb{R}$. From these quantities,  $v(\cdot),\,a(\cdot),\,\dot{\chi}(\cdot),\,\dot{\gamma}(\cdot)$ are computed by standard geometric arguments. 

\begin{figure}
	\centering
	\includegraphics[clip, trim= 1cm 0 1cm 0, width=\columnwidth]{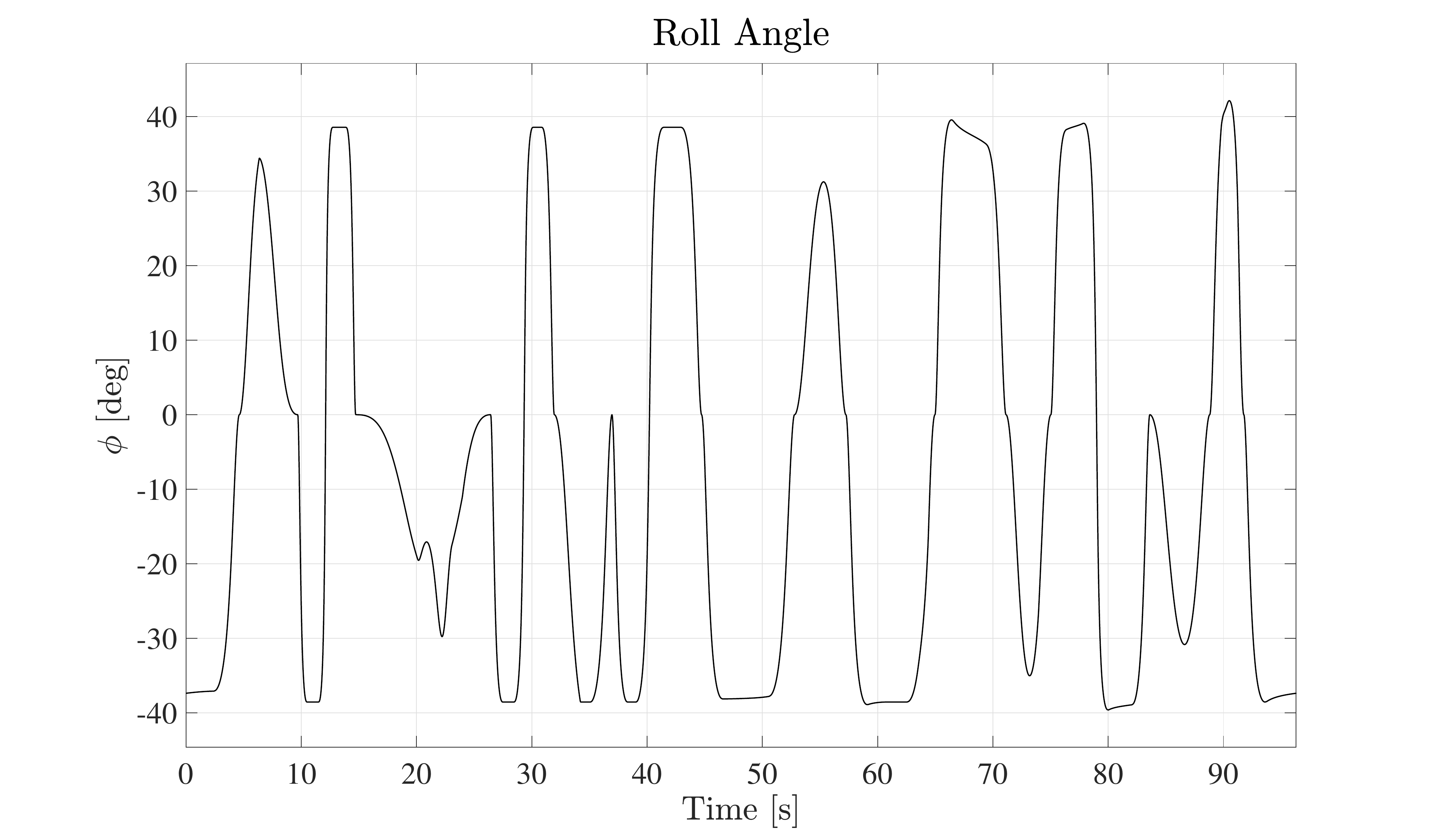}
	\caption{Motorbike lean angle generated accordingly with the simulation procedure described in Section \ref{sec:Simulation}. }
	\label{fig:ImolaPhi}
\end{figure}

As for the generation of the Euler angles and the body rotational speeds, in agreement with [\cite{Maceira2021Roll}, Eq. (32)], we assumed 
a linear law $\Delta\phi = k \phi$ with $k >0 $ (tuned to keep the pilot gravity centre above the road surface) and we computed $\phi\,:\,[0\,\,L]\to\mathbb{R}$ accordingly to a coordinated turn, see Section \ref{sec:Coordinated}. The motorbike lean angle obtained with this procedure is shown in Figure \ref{fig:ImolaPhi}.

The body angular speeds, $\omega\,:\,[0,\,L]\to \mathbb{R}^3$ were obtained from $\phi(\cdot),\gamma(\cdot),\chi(\cdot)$ by using standard derivative arguments and under the assumption of coordinated manoeuvres. To let the reader able to replicate the simulations detailed in this Section, the parameters modelling the sensor suite are listed in Table \ref{table:SimParamters}. More in detail, accordingly to \cite{Farrell2021IMU}, the model \eqref{eq:GenSensorModel} is completed by taking $b_{i0} = \mathtt{col}(0,\sim \mathcal{N}(0,\sigma^2))$, and
\[
\resizebox{\columnwidth}{!}{$
w_{\#} = \left[\begin{array}{c}
	1\\1\\1
\end{array}\right] \otimes \left[\begin{array}{c}
\sim \mathcal{N}\left(0,\dfrac{2\log(2)}{\pi 0.4365^2} \dfrac{B^2}{\tau_{\#}^2}\right)\\ \sim \mathcal{N}\left(0,K^2\right)\\ \sim \mathcal{N}\left(0,N^2\right) 
\end{array}\right],\, \# \in\{a,g,s\}.
$}
\] 
	\begin{table}
		\begin{center}
				{\scriptsize
	\begin{tabular}{ |c  | c  ||c | c| c| }

\hline
    Symbol & \begin{tabular}{@{}c@{}}Unit \\ {[Gyros, Acc.s, GNSS]} \end{tabular}
 & Gyros & Acc.s & GNSS \\
    \hline
    \hline
     $B$ & [rad/s, m/s$^2$, m/s]& 3.3e-3 & 2.0e-4 & 0\\  \hline
     $\tau_{\#}$ & s &20 & 30 & 0 \\   \hline	
     $N$ &  [rad/s, m/s$^2$, m/s] & 8.5e-3 & 3.3e-3& 
     \begin{tabular}{@{} cc@{}} 1.1e-2 & $x,y$\\
     	1.0e-1 &$z$ \end{tabular}
     \\  \hline
     $K$ & [rad/s$^2$, m/s$^3$, m/s$^2$] & 0 & 3.3e-3 & 0\\  \hline
     $\sigma$ & [rad/s, m/s$^2$, m/s] & 5.0e-2 & 1.0e-1 & 0\\   \hline
     $T_s$ & s & 1.0e-2 & 1.0e-2 & 1.0e-1\\  \hline
	\end{tabular}}
\end{center}
\caption{List of parameters used for sensor noise generation.} 
\label{table:SimParamters}
\end{table}

\begin{figure}
	\centering
		\begin{subfigure}[t]{0.49\columnwidth}
		\centering
		\includegraphics[clip, trim= 0cm 0.5cm 3cm 0cm,width = \textwidth]{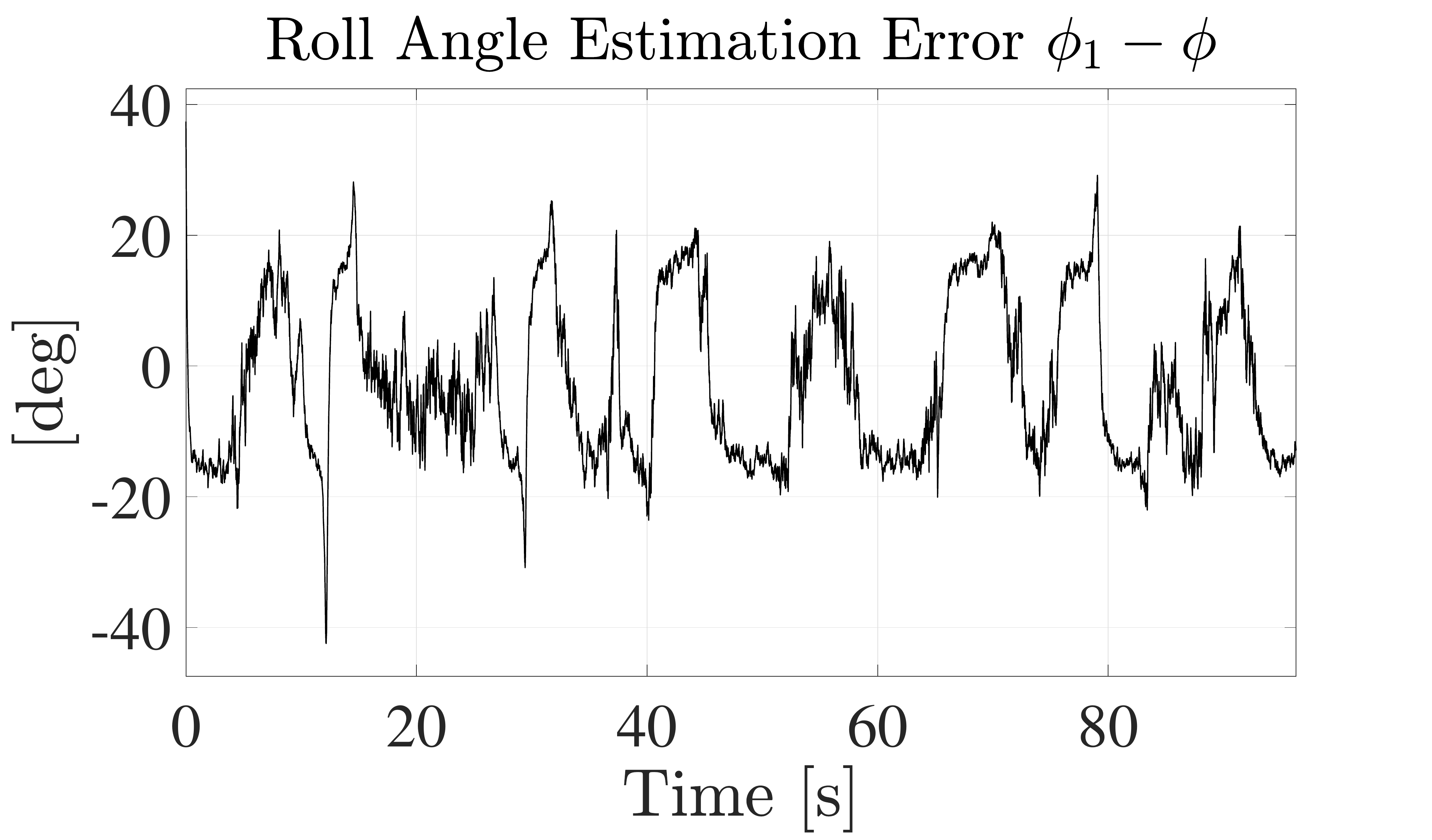}\\
		\caption{}
	\end{subfigure}
		\begin{subfigure}[t]{0.49\columnwidth}
	\centering
	\includegraphics[clip, trim= 0cm 0.5cm 3cm 0cm,width = \textwidth]{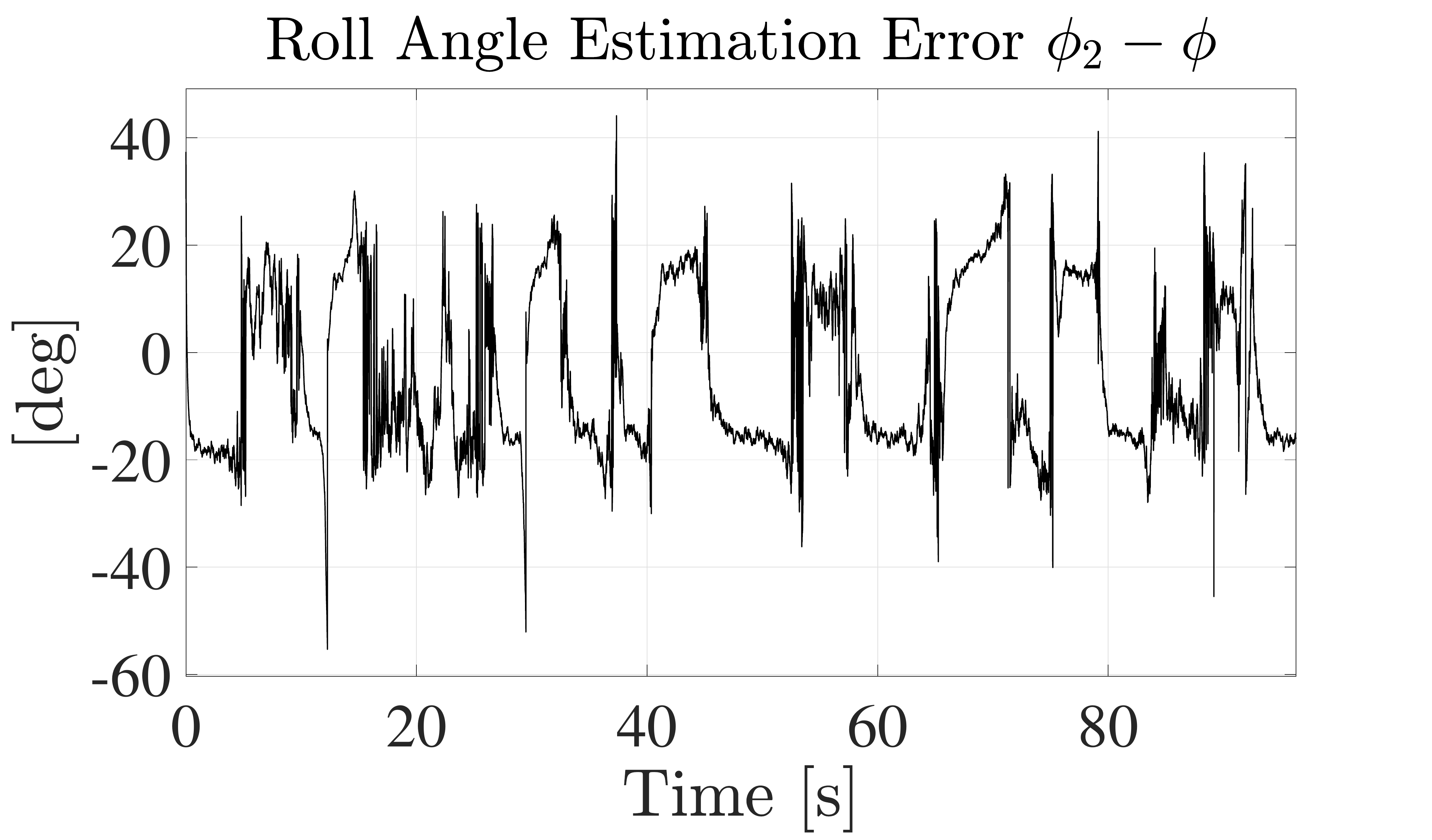}\\
	\caption{}
\end{subfigure}
		\begin{subfigure}[t]{0.49\columnwidth}
	\centering
	\includegraphics[clip, trim= 0cm 0.5cm 3cm 0cm,width = \textwidth]{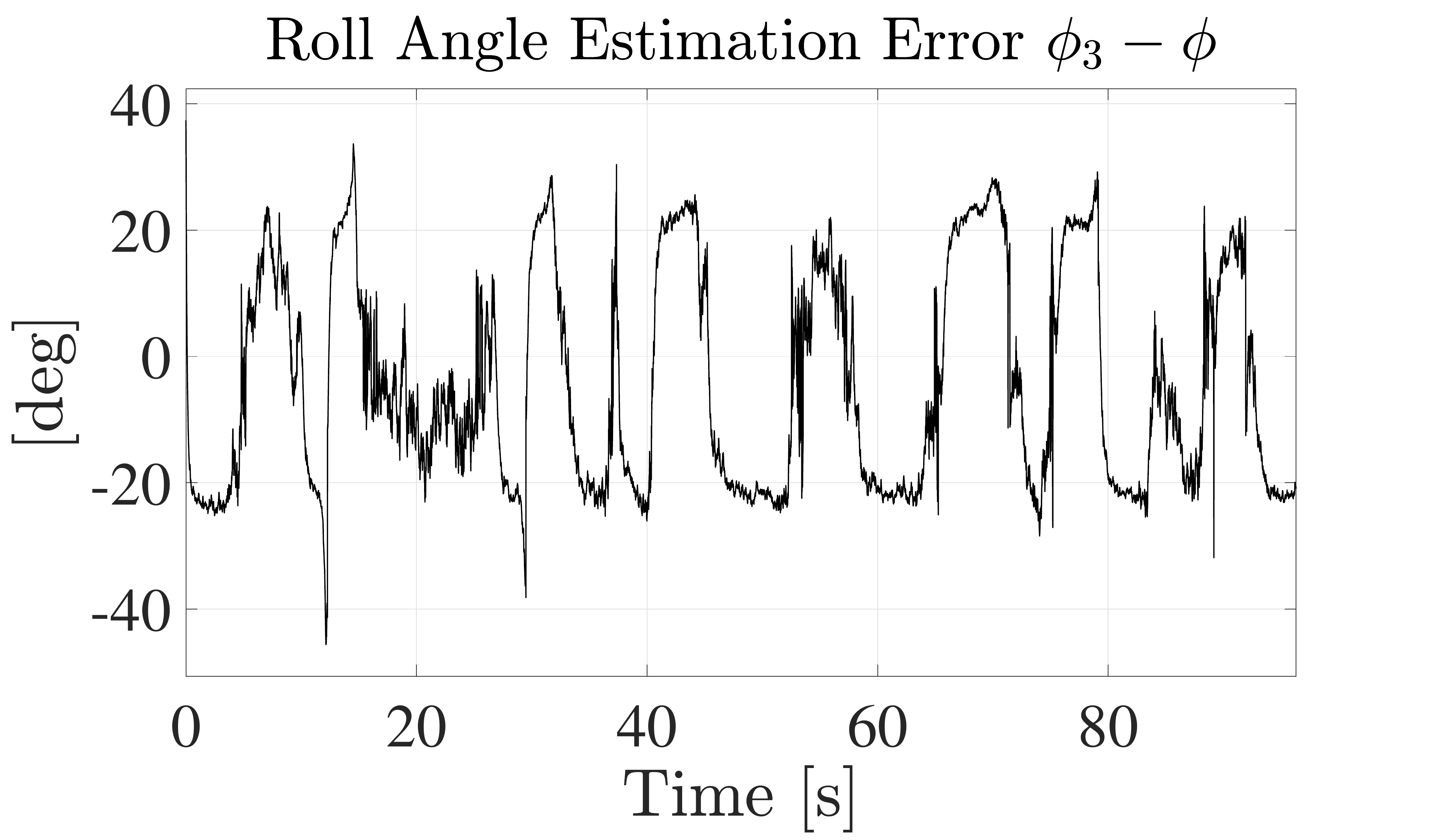}\\
	\caption{}
\end{subfigure}
\begin{subfigure}[t]{0.49\columnwidth}
	\centering
	\includegraphics[clip, trim= 0cm 0.5cm 3cm 0cm,width = \textwidth]{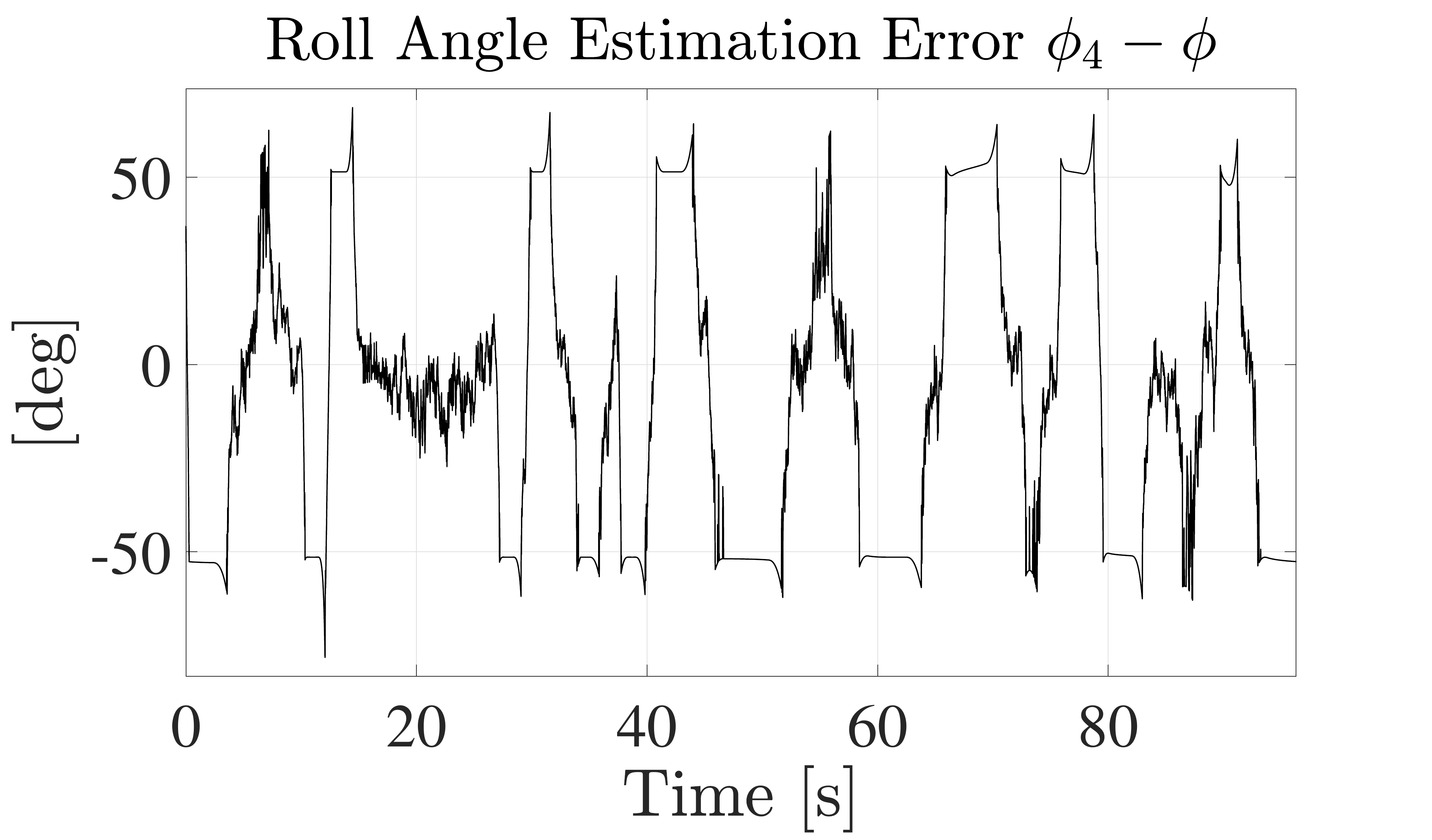}\\
	\caption{}
\end{subfigure}
		\begin{subfigure}[t]{0.49\columnwidth}
	\centering
	\includegraphics[clip, trim= 0cm 0.5cm 3cm 0cm,width = \textwidth]{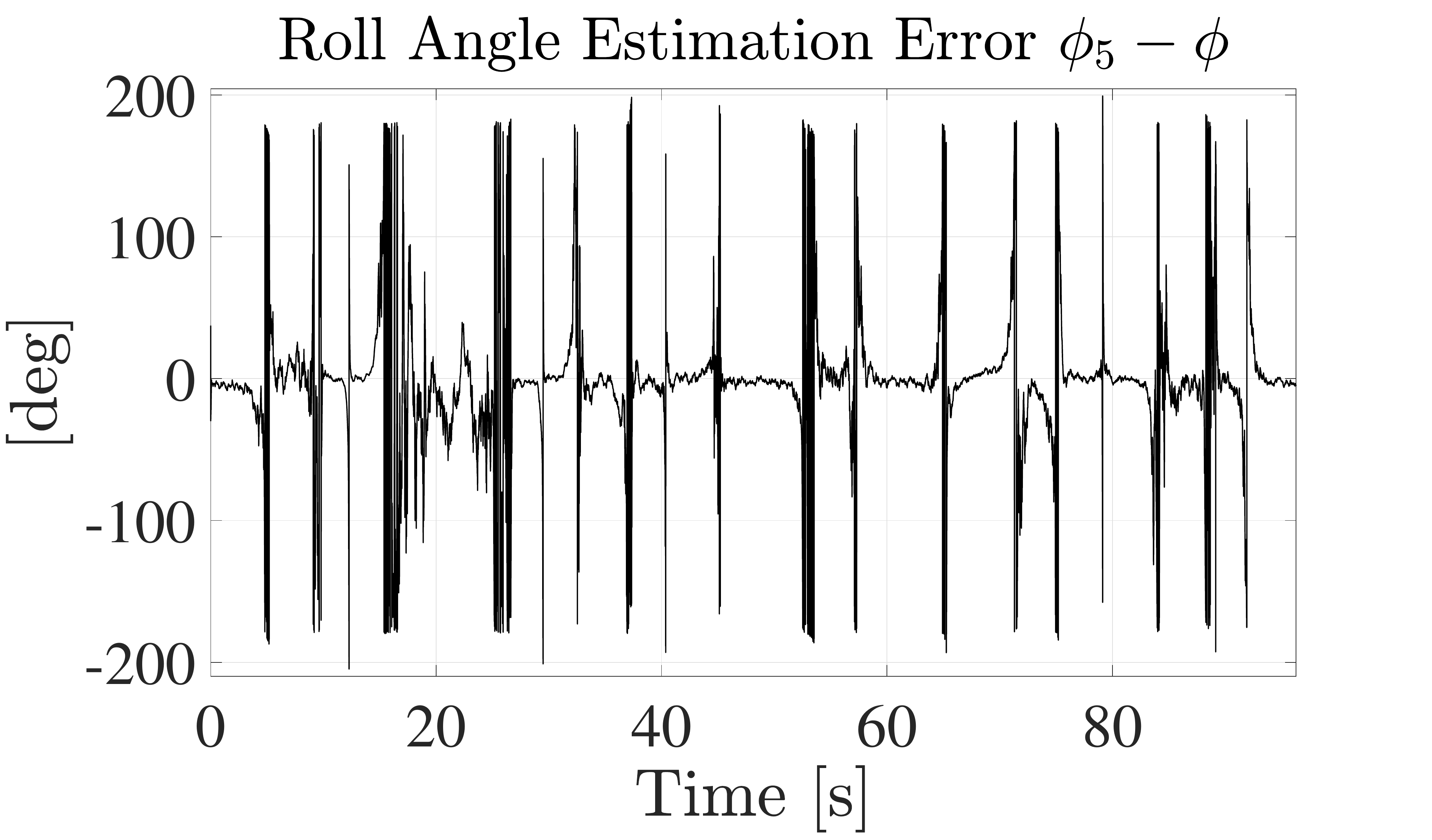}\\
	\caption{}
\end{subfigure}
\begin{subfigure}[t]{0.49\columnwidth}
	\centering
	\includegraphics[clip, trim= 0cm 0.5cm 3cm 0cm,width = \textwidth]{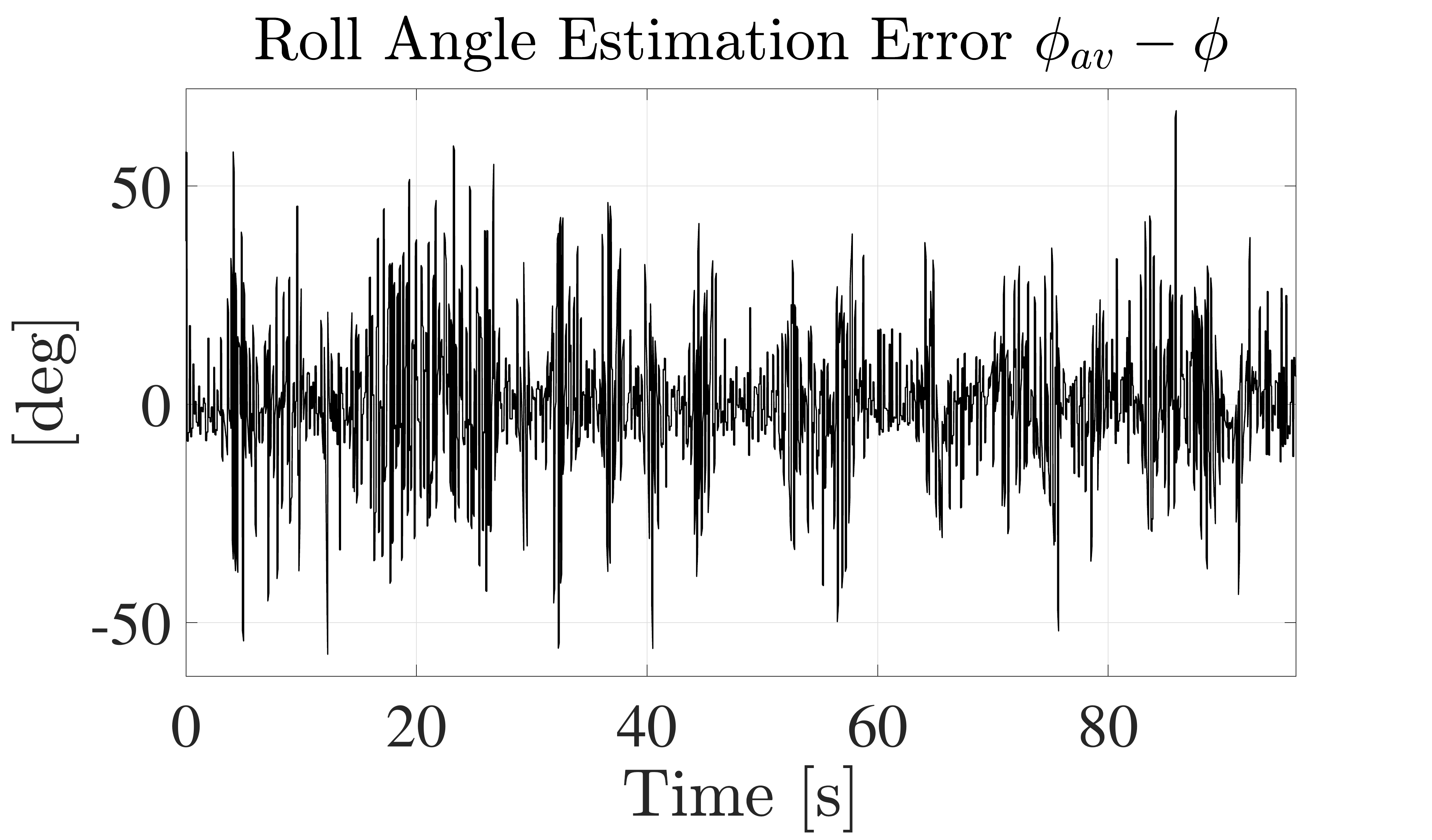}\\
	\caption{}
\end{subfigure}
\caption{Estimation errors $\tilde{\phi}_i$, for $i=1,\dots,5$, and $\tilde{\phi}_{\text{av}}$. These results show that the reference roll angles proposed by the cited literature suffer of low-frequency errors, see Figures from (a) to (e). At the opposite, the error $\tilde{\phi}_{\text{av}}$ seems to be more robust to low frequency errors.}
\label{fig:Comparison}
\end{figure}

\begin{figure}
	\centering
	\includegraphics[clip, trim= 2cm 0cm 2cm 0cm, width=\columnwidth]{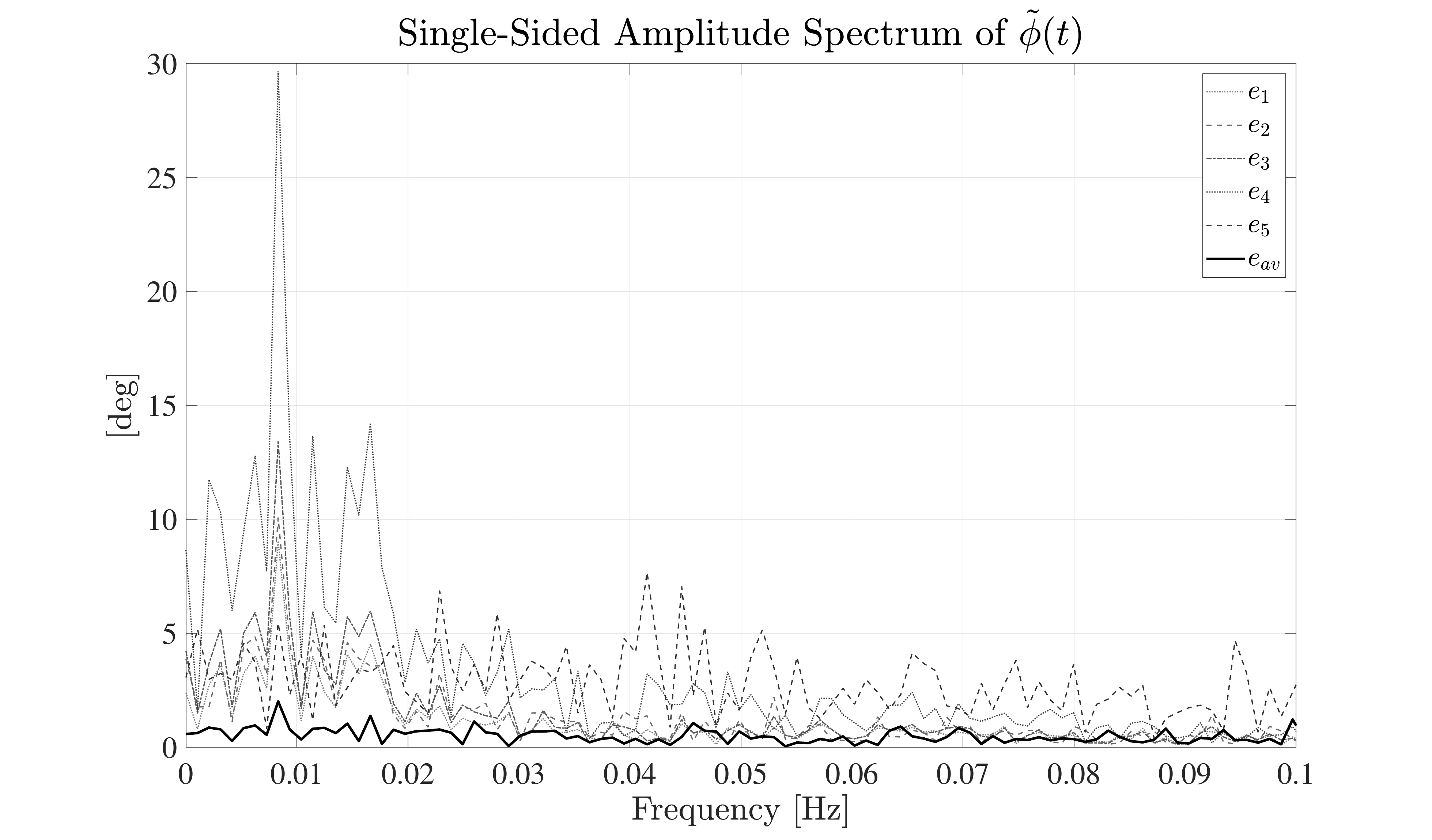}
	\caption{Spectra of errors $\tilde{\phi}_i-\phi$, for $i=1,\dots,5$, and $\tilde{\phi}_{\text{av}}$ obtained via Fast Fourier Transform. This picture highlights that the error $\tilde{\phi}_{\text{av}}$ is characterised by a more uniform spectrum that makes $\phi_{av}$ suitable for a observer implementation.}
	\label{fig:ImolaSpectrum}
\end{figure}

With reference to the roll angle profile of Figure \ref{fig:ImolaPhi}, Figures \ref{fig:Comparison}(a)-(e) report the errors $\tilde{\phi}_i := \phi_i-\phi$, for $i = 1, \dots,5$, while Figure \ref{fig:Comparison}(f) shows $\tilde{\phi}_{\text{av}} := \phi_\text{av}-\phi$. As can be seen, the estimation errors associated to the algorithms proposed by the cited literature have low frequency errors which seems to be not present in $\tilde{\phi}_{\text{av}}$. To confirm this visual intuition, Figure \ref{fig:ImolaSpectrum} displays the single-sided spectra of $\tilde{\phi}_i$,  for $i = 1,\dots,5$, and $\tilde{\phi}_{\text{av}}$.

\section{Conclusions}
The lean angle estimate proposed in this work has been formulated by implementing two key ideas: the definition of a two-stage-structure observer and the coordinated manoeuvre assumption. The first stage elaborates GNSS and accelerometer data to estimate an attitude quaternion accordingly with the assumption of coordinated manoeuvre. Then, the second stage, consisting of an EKF, integrates gyroscope data to improve the estimation accuracy during fast rolling manoeuvres. 
Compared to other lean angle estimation schemes, the coordinated manoeuvre assumption, particularly suitable for high-performance motorbikes, assures good estimation accuracy, as the plots in the experimental section show. Furthermore, theoretical proofs assure the observer estimation error is uniform with all the actual potential trajectories.

\bibliographystyle{ieeetr}
\bibliography{Bib,BibIntro}

\appendices

\section*{Acknowledgment}

We thank engineers Dario Zinelli and Mirko Sthylla, from AvioRace, for their crucial support in executing field tests and the precious discussions on IMU performance in motorsport applications.  We thanks also Giesse Racing Team, which allowed us to collect the data presented in this paper.

\section{Proofs and Computations}

%

\subsection{Proof of Theorem \ref{Theorem:DREBounds}}
\label{sec:ProofOfTheorem1}

Since the algorithm under investigation is a standard EKF, we rely on already available results to prove the observer stability. Hereafter, we report only those steps which are not conventionally found in the cited literature.

This proof is made of four parts: a) lower and upper bound for $S$, b) local asymptotic stability of $\tilde{x}=\hat{x}-x$,   c) demonstration that \eqref{eq:BasicEKF_CF} is not at the equilibrium for $\hat{q}=0$, and d) upper bound for $\|\hat{\phi}-\phi\|$.
 
As for a), we adopt \cite{Silverman1967Controllability} to show that, if $\hat{q}(t) \ne 0$ for all $t \ge 0$, the couple $(A(\hat{x}(t),t),H)$ is Uniformly Completely Observable (UCO). Then, \cite{Bristeau2010Design} shows that
 $(A(\hat{x}(t),t),H)$ UCO implies that the reconstructability Gramian is lower and upper bounded by positive matrices.  Then, exploit the reachability of $(A(\hat{x},t),I)$, the positiveness of $R(t)$ and $Q$, and  use \cite{Pengov2001} to demonstrate the existence of $\underline{s}$ and $\overline{s}$ verifying the first claim of Theorem \ref{Theorem:DREBounds}.

As for b), we rely on results of point a) and use standard EKF arguments (see \cite{besanccon2007nonlinear} and citations therein).

We prove point c) by contradiction. Define $x^\star = \mathtt{col}(\bar{b}_g^\star(t),0)$ and $S^\star(t)$ as en equilibrium manifold for \eqref{eq:BasicEKF_CF}, introduce $P_{ij}$, with $i,j = 1,2$, as matrices of suitable dimensions such that $$ P(t): = \left[\begin{array}{cc}
	P_{11}(t) & P_{12}(t)\\
	P_{21}(t) & P_{22}(t)
\end{array}\right] = S^{-1}(t).$$ Evaluate \eqref{eq:EKF_CF} at $x^\star(t)$ and $S^\star(t)$ to find that $\dot{\hat{q}} = 0$ if and only if $P_{22}(t)$ is singular for all $t \ge 0$. Let ${X}(t) \in \mathbb{R}^{4\times 4}$ such that $A(x^\star,t) = \mathtt{blkdiag}(0, X(t))$, 
and use \eqref{eq:DRE_CF} to compute the dynamics of
$P_{22}$ as
\[
\dot{P}_{22} = 
{X}(t)P_{22} + P_{22} {X}^\top(t) - P_{22}
R^{-1}(t)P_{22} + \epsilon I.
\]		
Finally, since this represents a differential Riccati equation associated with the fully observable couple $({X}(t),I)$ and $\epsilon I, R(t) \succ 0$, $P_{22}(t)$ cannot be singular $\forall\, t \ge 0$.

As for d), we rely on b) and perform algebraic computations to show that 
$$\|\hat{\phi}-\phi\| \le \|\left[\begin{array}{ccc}
	1 & 0 & 0
\end{array}\right] \partial  \bar{h}^{-1}(q) / \partial q\| \|\tilde{q}\|$$ where Assumption \ref{hyp:Pitch} guarantees $$\|\left[\begin{array}{ccc}
1 & 0 & 0
\end{array}\right] \partial  \bar{h}^{-1}(q) / \partial q\| \le \sqrt{2}/(2\cos{\overline{\theta}}).$$  

\subsection{Proof of Lemma \ref{lemma:PhiAv}}
\label{sec:ProofPhiAv}
We define three reference frames to prove \eqref{eq:phi_av2}: Inertial, Body, and Navigation. Let $v^N$ and $v^B$ 
be the projection on the Navigation and Body frames of the inertial speed $v$. Then, we assume that $v^N = v^B 
 = \mathtt{col}(\mathtt{v},0,0)$.
 We define 			$R_{BI} = R_1(\phi) R_2(\theta)R_3(\psi)$ and $R_{NI} = R_2(\gamma)R_3(\chi)$ where $R_i(s)$ denotes the  matrix associated with a rotation, of magnitude $s$, around the $i$-th axis.
Exploit these matrices and the definition of $v^N$ to write
$
v^B = R_{BI} R^\top_{NI} v^N.
$ 
 Explicit the matrices in the latter equation to obtain
\[
v^B= R_1(\phi) R_2(\theta)R_3(\psi)  R_3^\top(\chi) R_2^\top(\gamma) v^N.
\]
Since $v^N 
 = \mathtt{col}(\mathtt{v},0,0)$ by definition, the previous equality is verified if $\psi = \chi$ and $\theta=\gamma$. 
We exploit $a = R_{BI}(\dot{v}-g)$ to compute $\phi$ as follows.  
Explicit $R_{BI} = R_1(\phi) R_2(\gamma)R_3(\chi) = R_1(\phi) R_{NI} $ and define $a^N:=\mathtt{col}(a_x^N,a_y^N,a_z^N) = R_{NI} (\dot{v}-g)$. Then $a = R_1(\phi)a^N$ 
from which, with Assumption \ref{hyp:Non-Ballistic},
\[
\phi = \tan^{-1}\left(\dfrac{a_z^N a_y -a_y^N a_z}{a_y^Na_y + a_z^Na_z}\right).
\]

\subsection{Definition of $f_{\xi_e}$}
\label{sec:fxie}
Let
\begin{equation*}
		\begin{aligned}
			\mathcal{T}^{+}:=&\, \{t \ge 0\,:\,{v}_y(t) = 0,\,{\dot{v}}_y(t) < 0,\,{v}_x(t) < 0\}\\
			\mathcal{T}^{-}:=&\, \{t \ge 0\,:\,{v}_y(t) = 0,\,{\dot{v}}_y(t) > 0,\,{v}_x(t) < 0\},\\
		\end{aligned}	
\end{equation*} 
and define two lap counters $N_{+}$ and $N_{-} \in \mathbb{N}$, whose dynamics is
	\begin{equation}
	\label{eq:LapCounter}
		\begin{aligned}
			\dot{N}_{\$} = &\, 0 && t \ge 0\\
			N_{\$}^+ =&\, N_{\$}+1 &&t \in \mathcal{T}^{\$},  
		\end{aligned}
	\end{equation}
	with $\$ \in \{+, -\}$, and where $N_{\$}(0)=0$. Then, take
\[
f_{\chi}(v,N_{+},N_{-}):=  \mathtt{atan}_2(v_y,v_x) + 2\pi(N_{+}-N_{-})
\]
and use it and \eqref{eq:v_aero} to build
\[
f_{\xi}(v)  =  \mathtt{col}(
\sqrt{v^\top v},
-\sin^{-1}(v_z/\sqrt{v^\top v}),
f_\chi(v,N_+,N_-)
)
\]
and
\[
 f_{\dot{\xi}}(v_e) = \left[\left.\dfrac{\partial h_v(\xi)}{\partial \xi}\right|_{\xi = f_{\xi}(v)}\right]^{-1} \dot{v}.
\]
Finally, define
\[
f_{\xi_e}(v_e) = \mathtt{col}(	f_{\xi}(v),f_{\dot{\xi}}(v_e)).
\]

\subsection{Proof of Lemma \ref{lemma:ContinuityHatChi}}
\label{sec:ProofContinuityHatChi}
Let $t_\star \in \mathcal{T}_k^{-}\cup \mathcal{T}_k^{+}$, then the continuity $\mathcal{C}^1$ of $\hat{\chi}(t)$ is implied by
\[
\lim_{t \to t_\star^-} \hat{\chi}(t) = \lim_{t \to t_\star^+} \hat{\chi}(t)\qquad \forall\, t_\star \in \mathcal{T}_k^{-}\cup \mathcal{T}_k^{+},
\]
continuity $\mathcal{C}^1$ of $\mathtt{atan}_2(\hat{v}_y(t),\hat{v}_x(t))$ for all $t \not\in \mathcal{T}_k^{-}\cup \mathcal{T}_k^{+}$, and
\[
\lim_{t \to t_\star^-} \dfrac{d}{dt}\mathtt{atan}_2(\hat{v}_y(t),\hat{v}_x(t)) = \lim_{t \to t_\star^+} \dfrac{d}{dt}\mathtt{atan}_2(\hat{v}_y(t),\hat{v}_x(t))
\]
for all $t_\star \in \mathcal{T}_k^{-}\cup \mathcal{T}_k^{+}$.

\end{document}